\documentclass[a4paper]{amsart}
\usepackage{fullpage}
\usepackage{amsmath,amssymb,mathtools,booktabs,tikz}
\usetikzlibrary{arrows.meta,calc,decorations.markings,positioning}
\usepackage[T1]{fontenc}
\usepackage{lmodern}
\usepackage{microtype}
\usepackage{needspace}
\usepackage{tabularx,array}
\usepackage[hypertexnames=false]{hyperref}
\usepackage{doi}
\usepackage{float}
\hypersetup{colorlinks=true,linkcolor=blue,citecolor=blue,urlcolor=blue}
\newtheorem{thm}{Theorem}[section]
\numberwithin{equation}{section}
\newtheorem{lem}[thm]{Lemma}
\newtheorem{prop}[thm]{Proposition}
\newtheorem{cor}[thm]{Corollary}
\newtheorem{conjecture}[thm]{Conjecture}
\newtheorem{rem}[thm]{Remark}
\newtheorem{eg}[thm]{Example}

\newcommand{\C}{\mathbb C}
\newcommand{\R}{\mathbb R}

\newcommand{\Z}{\mathbb Z}

\title[The Erd\H{o}s--Herzog--Piranian lemniscate conjecture]
{A proof of the Erd\H{o}s--Herzog--Piranian lemniscate conjecture}
\author{Teng Zhang}
\address{School of Mathematics and Statistics, Xi'an Jiaotong University,
Xi'an 710049, P. R. China}
\email{teng.zhang@stu.xjtu.edu.cn}
\date{\today}
\begin{document}
\begin{abstract}
Let \(p:\mathbb{C}\to\mathbb{C}\) be a monic polynomial of
degree \(n\geq 1\). In 1958, Erd\H{o}s, Herzog, and Piranian conjectured
that the length of the lemniscate
\(\{z\in\mathbb C:|p(z)|=1\}\) is maximized by \(p(z)=z^n-1\).
We prove this conjecture in every degree and show that equality
holds precisely for \(p(z)=(z-z_0)^n-e^{i\theta}\), where
\(z_0\in\mathbb C\) and \(\theta\in\mathbb R\).
The proof uses Möbius variations of the critical values and estimates
for the inverse branches of \(p\), supplemented by rigorous
computer-assisted estimates in a finite range of degrees.
\end{abstract}
\maketitle

\section{Introduction}
\subsection{Erd\H{o}s--Herzog--Piranian lemniscate conjecture}
Let \(p:\C\to\C\) be a monic polynomial of degree \(n\ge1\), and
write \(\mathbb T=\{w\in\C:|w|=1\}\). For \(R>0\), let
\[
E_R=E_R(p):=\{z\in\C:|p(z)|<R\},\qquad
\partial E_R(p)=\{z\in\C:|p(z)|=R\}.
\]
We call \(\partial E_1(p)=p^{-1}(\mathbb T)\) the
\emph{Erd\H{o}s--Herzog--Piranian lemniscate} of \(p\) and
denote its arclength by \(\ell(\partial E_1(p))\).

In 1958, Erd\H{o}s, Herzog, and Piranian~\cite[Problem~12]{EHP58} (see also Erd\H{o}s~\cite[Problem~1]{Erd83})
conjectured
that the length $\ell(\partial E_1(p))$ of the lemniscate
\(\partial E_1(p)\) is maximized by \(p(z)=z^n-1\).
\begin{conjecture}[Erd\H{o}s--Herzog--Piranian]\label{conj}
Let \(p:\mathbb{C}\to\mathbb{C}\) be a monic complex polynomial of
degree \(n\geq 1\),
and set \(p_0(z):=z^n-1\). Then
\[
\ell\bigl(\partial E_1(p)\bigr)
\leq
\ell\bigl(\partial E_1(p_0)\bigr).
\]
\end{conjecture}

\begin{figure}[H]
\centering
\includegraphics[width=0.95\textwidth]{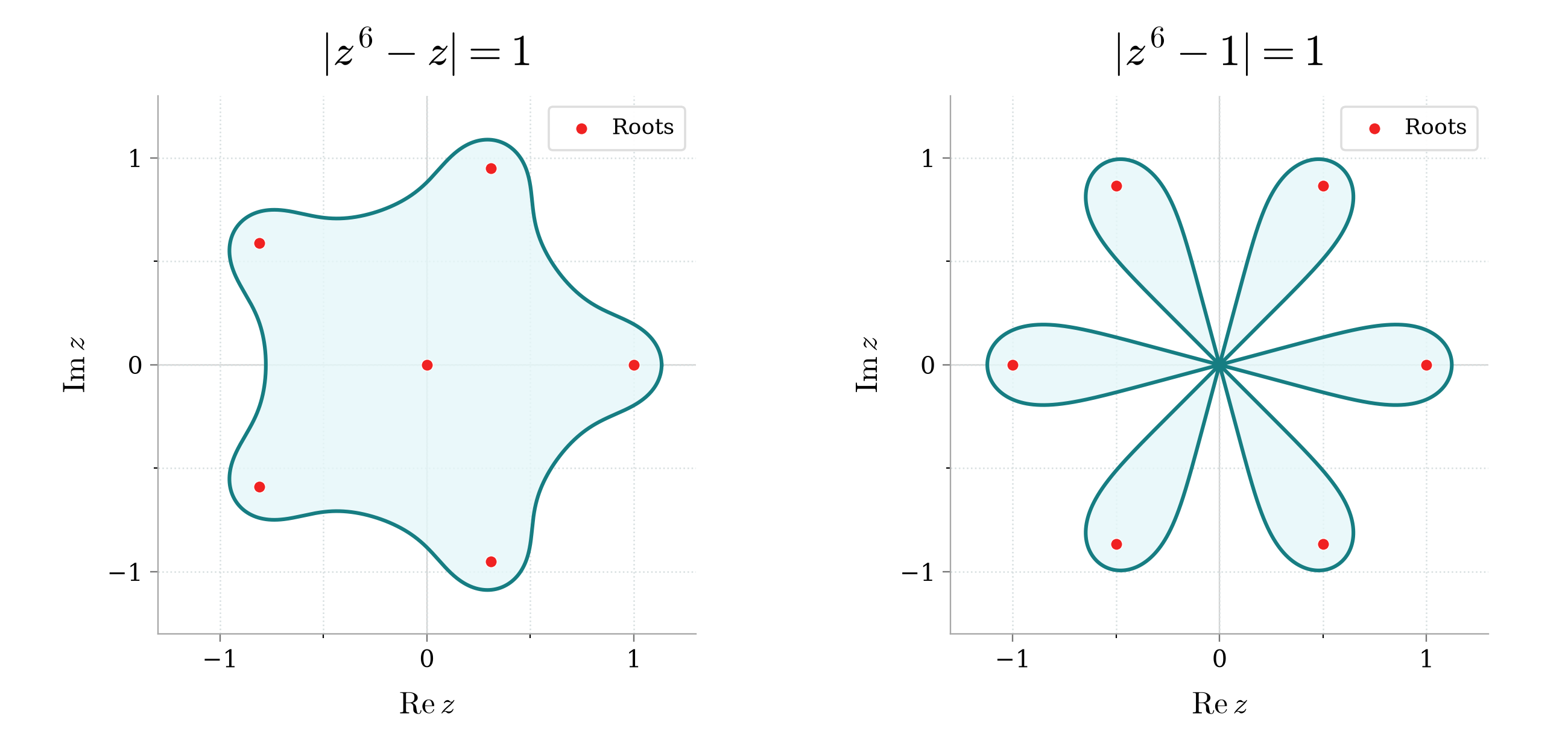}
\caption{The lemniscates \(\partial E_1(z^6-z)\) (left) and
\(\partial E_1(z^6-1)\) (right). The shaded regions indicate
the corresponding sets \(E_1(p)\), and the red dots mark
the zeros of the polynomials.}
\label{fig:lemniscates-n6}
\end{figure}

Figure~\ref{fig:lemniscates-n6} compares the lemniscates associated
with \(z^6-z\) and \(z^6-1\). The latter consists of six petals
meeting at the origin.

\begin{rem}\label{rem}
As suggested by Tao~\cite[Remark~1.1]{Tao25}, it is natural
to conjecture that equality in Conjecture~\ref{conj} holds
if and only if
$
p(z)=(z-z_0)^n-e^{i\theta}
$
for some \(z_0\in\mathbb{C}\) and \(\theta\in\mathbb{R}\).
In other words, the conjectured maximizers are precisely
the translates and rotations of \(p_0\). 

The model length is
\[
\ell\bigl(\partial E_1(p_0)\bigr)
 =2^{1/n}B\!\left(\frac12,\frac1{2n}\right),
\]
where \(B\) is Euler's beta function; see Danchenko~\cite{Dan07}. An equivalent
formula in terms of the gamma function see also Elia and Galizia Angeli~\cite[Eq.~(13)]{EG84}.
\end{rem}
Several upper bounds for \(\ell(\partial E_1(p))\) have been obtained
towards Conjecture~\ref{conj}; see Table~\ref{tab:lemniscate-progress}\footnote{This table updates the list of upper bounds compiled by Tao~\cite[Table~1]{Tao25}.}.
In particular, the case \(n=1\) is immediate, and
Eremenko and Hayman~\cite{EH99} settled the case \(n=2\).
Dahlke~\cite{Dah26} proved the case \(n=3\) and the correponding
characterized equality see
Plasma AI~\cite{Pla26}.
Tao~\cite{Tao25} proved the conjecture and
characterized equality for all sufficiently large \(n\).
\par\nopagebreak[4]
\begin{table}[H]
\centering
\small
\renewcommand{\arraystretch}{1.18}
\setlength{\tabcolsep}{5pt}
\begin{tabularx}{0.8\textwidth}{
    |c|c|>{\centering\arraybackslash}X|
}
\hline
Upper bound & Hypotheses & Reference \\
\hline
\(74n^2\)
&
& Pommerenke~\cite{Pom61} \\

\(8e\pi n\)
&
& Borwein~\cite{Bor95} \\

\(4\pi n\)
&
& Dolzhenko~\cite{Dol63} \\

\(4\pi n-2\pi\)
&
& Fryntov--Nazarov~\cite{FN09} \\

\(9.173n\)
&
& Eremenko--Hayman~\cite{EH99} \\

\(2\pi n+2\pi\sqrt n-2\pi\)
&
& Fryntov--Nazarov~\cite{FN09} \\

\(2\pi n\)
&
& Danchenko~\cite{Dan07} \\

\(\pi n+O(\sqrt{n\log n})\)
&
& Kosukhin~\cite{Kos12} \\
\hline
\(2n+O(n^{7/8})\)
&
& Fryntov--Nazarov~\cite{FN09} \\

\(\ell(\partial E_1(p_0))\)
& \(n=1\)
& trivial \\

\(\ell(\partial E_1(p_0))\)
& \(n=2\)
& Eremenko--Hayman~\cite{EH99} \\

\(\ell(\partial E_1(p_0))\)
& \(p(z)=z^n+c\)
& Wang--Peng~\cite{WP06} \\

\(\ell(\partial E_1(p_0))\)
& \(p\) close to \(p_0\)
& Fryntov--Nazarov~\cite{FN09} \\

\(\ell(\partial E_1(p_0))\)
& \(n=3\)
& Dahlke~\cite{Dah26} \\
\hline
\(2n+O(\sqrt n)\)
&
& Tao~\cite{Tao25} \\

\(2n+O(1)\)
&
& Tao~\cite{Tao25} \\

\(2n+4\log 2+o(1)\)
&
& Tao~\cite{Tao25} \\

\(\ell(\partial E_1(p_0))\)
& \(n\) sufficiently large
& Tao~\cite{Tao25} \\
\hline
\(\ell(\partial E_1(p_0))\)
& \(n\ge1\)
& Theorem~\ref{thm:main}  \\
\hline
\end{tabularx}
\par\smallskip
\caption{A summary of upper bounds on the lemniscate length
\(\ell(\partial E_1(p))\) for monic polynomials \(p\) of degree \(n\), including the main result of the present paper,
where \(p_0(z)=z^n-1\). Blank entries in the second column indicate
\(n\ge2\) with no additional hypothesis. All asymptotic estimates
are uniform in \(p\).}
\label{tab:lemniscate-progress}
\end{table}

In this paper, we prove Conjecture~\ref{conj} in every degree and
show that equality holds precisely in the cases described in
Remark~\ref{rem}.
\begin{thm}\label{thm:main}
Let \(p:\mathbb{C}\to\mathbb{C}\) be a monic complex polynomial of
degree \(n\geq 1\),
and set \(p_0(z):=z^n-1\). Then
\[
\ell\bigl(\partial E_1(p)\bigr)
\leq
\ell\bigl(\partial E_1(p_0)\bigr).
\]
Equality holds if and only if \(p(z)=(z-z_0)^n-e^{i\theta}\)
for some \(z_0\in\C\) and \(\theta\in\R\).
\end{thm}

\subsection{High-level description of the proof}
Fryntov and Nazarov~\cite{FN09} compare intersections after
normalizing the target line. For \(z\mapsto z^n\), the inverse
image of \(\mathbb R\) meets almost every source line in \(n\)
points. A general monic polynomial may have fewer intersections.
We deform the target circle to a line and estimate how much of
this loss survives on returning to the circle. We assume
\(n\ge2\) throughout the outline. For $n=1$, the lemniscate
is a unit circle.

\medskip
\noindent
\(\bullet\)\enspace\textbf{Step I. Choosing a maximizer with critical values on the unit circle.}
A point \(\zeta\in\mathbb{C}\) is called a \emph{critical point}
of \(p\) if \(p'(\zeta)=0\). A number \(v\in\mathbb{C}\) is called
a \emph{critical value} of \(p\) if \(v=p(\zeta)\) for some
critical point \(\zeta\) of \(p\).

First, we observe that the extremal polynomial \(p_0(z)=z^n-1\)  has all its critical points at the
origin and a single finite critical value, namely \(-1\).
Conversely, a monic polynomial with a single finite critical
value \(v\) must have the form
\[
p(z)=(z-z_0)^n+v.
\]
This motivates us to prove that 
\begin{center}
    If \(n\ge2\) and \(p\) maximizes \(\ell(\partial E_1(p))\)
among all monic polynomials of degree \(n\), then \(p\) has
a single finite critical value of modulus one.
\end{center}

Eremenko and Hayman~\cite{EH99} show that a maximizing polynomial
can be chosen with all its critical values on the unit circle. Thus, as an immediate consequence, Conjecture~\ref{conj} holds for \(n=2\).
To prove Theorem~\ref{thm:main}, we may therefore consider such
a maximizer, Tao~\cite{Tao25}
also starts with this choice.  For the length inequality it therefore
suffices to prove that a maximizer chosen in this way has only
one finite critical value. Step~VI treats an arbitrary maximizer.

\medskip
\noindent
\(\bullet\)\enspace\textbf{Step II. Lifting a M\"obius deformation of the target circle.}
Eremenko and Hayman~\cite{EH99}
lift a local motion of one critical value to a deformation of the
polynomial. Here we move all the critical values by the same
M\"obius transformation. Since \(\partial E_1(p)=p^{-1}(\mathbb T)\),
this also gives us a family of target circles.

For \(|\beta|<1\), consider the Möbius transformation
\[
M_\beta(w)=\frac{w}{1-\beta w}.
\]
We apply \(M_\beta\) both to the value \(w=p(z)\) and to the unit
circle in the target plane. This replaces \(p\) by the rational
map \(M_\beta\circ p=p/(1-\beta p)\) and the circle by \(M_\beta(\mathbb T)\).
Since \(M_\beta\) is injective,
\begin{equation}\label{eq:mobius-circle}
|p(z)|=1
\quad\Longleftrightarrow\quad
M_\beta(p(z))\in M_\beta(\mathbb T).
\end{equation}
Thus the two conditions select the same points \(z\).

\begin{rem}
In the coordinate \(u=M_\beta(w)\), the image circle has equation
\[
(1-|\beta|^2)|u|^2-2\operatorname{Re}(\beta u)=1.
\]
As \(\beta\to e^{i\theta}\) from within the unit disk,
this circle converges locally to the line
\begin{equation}\label{eq:limiting-target-line}
\{u\in\mathbb C:\operatorname{Re}(e^{i\theta}u)=-1/2\}.
\end{equation}
\end{rem}

To lift this particular motion, we divide \(F^2\) by \(F'\) at the
current monic polynomial \(F\). The remainder prescribes the change
in the polynomial, and the quotient prescribes the change of source
coordinate. The resulting polynomial \(F_\beta\) is monic of degree
\(n\), with \(F_0=p\). The source map \(\Phi_\beta\) is univalent
near \(\partial E_1(p)\), with \(\Phi_0(z)=z\), and satisfies
\begin{equation}\label{eq:mobius-lift}
F_\beta(\Phi_\beta(z))=M_\beta(p(z)).
\end{equation}

By \eqref{eq:mobius-circle} and \eqref{eq:mobius-lift}, the
source flow carries each fibre over \(\mathbb T\) onto the
corresponding fibre over \(M_\beta(\mathbb T)\). Indeed, the flow
is injective and preserves the local multiplicities, whose sum in
each fibre is \(n\). Hence
\begin{equation}\label{eq:deformed-lemniscate}
\Phi_\beta\bigl(\partial E_1(p)\bigr)
  =F_\beta^{-1}\bigl(M_\beta(\mathbb T)\bigr).
\end{equation}
Integrating along the regular arcs in \eqref{eq:deformed-lemniscate} gives
\[
\ell\bigl(F_\beta^{-1}(M_\beta(\mathbb T))\bigr)
 =\int_{\partial E_1(p)}|\Phi_\beta'(z)|\,|dz|.
\]
Write \(\beta=r\xi\), where \(0<r<1\) and \(\xi\in\mathbb T\).
If \(1-\xi v\ne0\) for every critical value \(v\) of \(p\), the
polynomials \(F_{r\xi}\) converge in coefficient space to a
polynomial \(F_\xi\). Its critical values lie on the line in
\eqref{eq:limiting-target-line}, with \(\xi=e^{i\theta}\).
We average the length integral in \(\xi\) before letting
\(r\uparrow1\).

\begin{figure}[H]
\centering
\includegraphics[width=\linewidth]{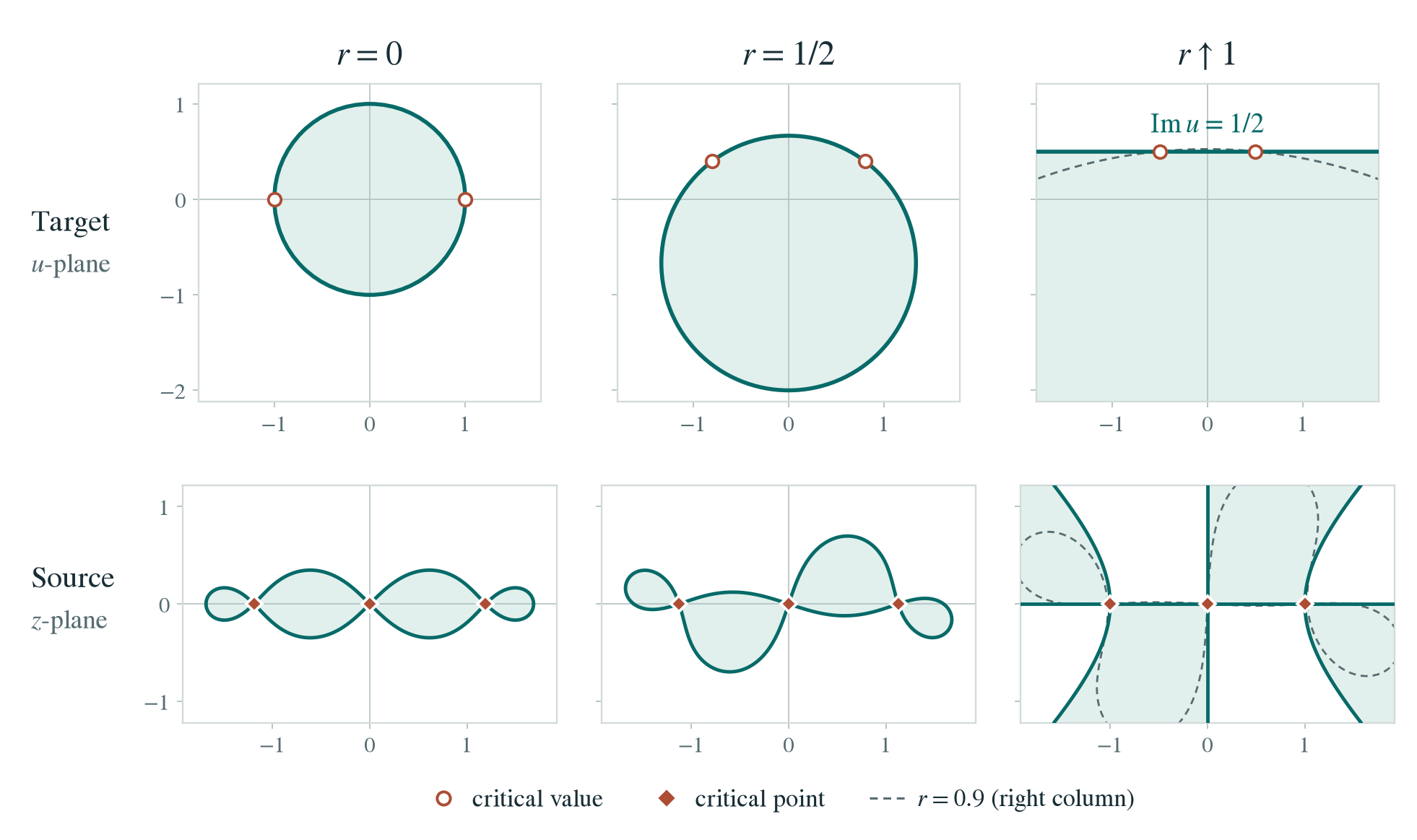}
\caption{The deformation of \(p(z)=(z^2-\sqrt{2})^2-1\) along
\(\beta=ir\). The top row shows \(M_{ir}(\mathbb T)\), and the bottom
row shows its inverse image under
\(F_{ir}(z)=(z^2-\sqrt{2/(1+r^2)})^2+M_{ir}(-1)\).
Shading marks the image of the open unit disk and its inverse image,
respectively. In the last column the solid curves show the limit on
a fixed window. The target line is \(\operatorname{Im}u=1/2\), and its
inverse image is \(\{z:z^4-2z^2\in\mathbb R\}\).}
\label{fig:ehp-moving-circle}
\end{figure}

Figure~\ref{fig:ehp-moving-circle} shows the moving circle and its
full inverse image for \(\xi=i\), when the limiting target line is
horizontal.
\par\nopagebreak[4]

\medskip
\noindent
\(\bullet\)\enspace\textbf{Step III. Counting missing intersections on a line.}
Fryntov and Nazarov~\cite{FN09} estimated the length of
\(\{\Re F=0\}\) in a disk by counting intersections with lines.
We use their counting argument with a common truncation in the
target plane. A translation and rotation take the limiting line
in \eqref{eq:limiting-target-line} to \(\R\), and a rotation of
the source preserves monicity.

\begin{eg}\label{eg:power-map}
For the model map \(z\mapsto z^n\), the inverse image of
\(\mathbb R\) consists of \(n\) lines through the origin, given by
\[
\{z\in\mathbb C:z^n\in\mathbb R\}
=
\bigcup_{k=0}^{n-1}e^{ik\pi/n}\mathbb R.
\]
For \(R>0\), the inverse image of the interval
\([-R^n,R^n]\) is
\[
\{z\in\mathbb C:z^n\in[-R^n,R^n]\}
=
\bigcup_{k=0}^{n-1}e^{ik\pi/n}[-R,R].
\]
These \(n\) segments have length \(2R\) each and meet only
at the origin. Hence
\[
\ell\bigl(\{z\in\mathbb C:z^n\in[-R^n,R^n]\}\bigr)=2nR.
\]
\end{eg}

The inverse image of \([-R^n,R^n]\) under \(z\mapsto z^n\)
has length \(2nR\). For a general monic polynomial \(F\), we
compare the two lengths through the limiting difference
\[
\Delta(F,\mathbb R)
=
\lim_{R\to\infty}
\left(
2nR-\ell\bigl(F^{-1}([-R^n,R^n])\bigr)
\right).
\]
Applying the planar Cauchy--Crofton formula in the normalization
of Dolzhenko~\cite{Dol63} and passing to the limit
in the truncation gives
\begin{equation}\label{eq:crofton-deficit}
\Delta(F,\mathbb R)
=
\frac12\int_0^\pi\int_{\mathbb R}
\bigl(n-N(\theta,h)\bigr)\,dh\,d\theta,
\end{equation}
where \(N(\theta,h)\) counts the intersections of
\(F^{-1}(\mathbb R)\) with the line
\(e^{i\theta}(\mathbb R+ih)\).

For almost every \((\theta,h)\), these intersections correspond
to the distinct real solutions of
\[
\operatorname{Im}F\bigl(e^{i\theta}(x+ih)\bigr)=0.
\]
The left-hand side is a nonzero real polynomial in \(x\)
of degree at most \(n\), so \(N(\theta,h)\le n\).
For the model map \(z\mapsto z^n\), the inverse image of
\(\mathbb R\) meets almost every source line in exactly
\(n\) points. Thus \(n-N(\theta,h)\) counts the intersections
missing relative to the model. In particular,
\eqref{eq:crofton-deficit} gives
\(\Delta(F,\mathbb R)\ge0\).

\begin{figure}[H]
\centering
\includegraphics[width=\linewidth]{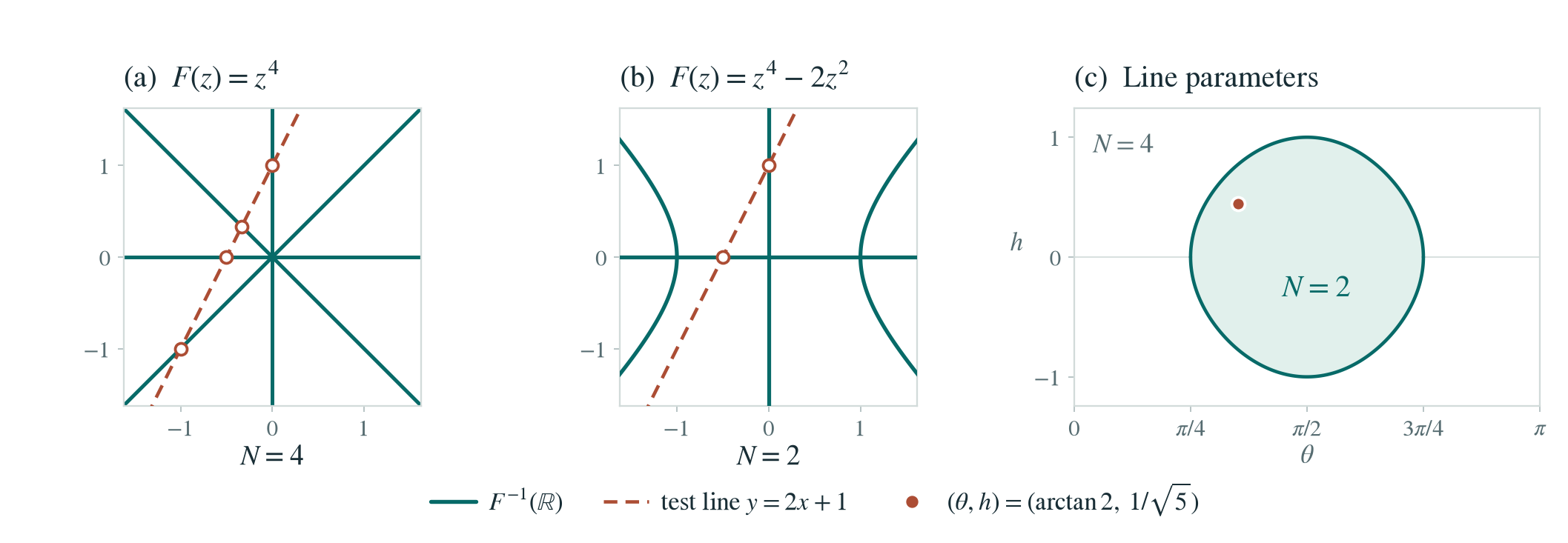}
\caption{The line \(y=2x+1\) meets the inverse image of the real axis
in four points for \(F(z)=z^4\) and in two for \(F(z)=z^4-2z^2\).
For the latter polynomial the inverse image consists of the two
coordinate axes and the hyperbola \(x^2-y^2=1\).
On writing \(z=e^{i\theta}(X+ih)\), the equation of the hyperbola
is quadratic in \(X\) for almost every \(\theta\), with discriminant
\(4(h^2+\cos(2\theta))\). Thus \(N(\theta,h)=2\) for almost every
parameter pair in the shaded region \(h^2+\cos(2\theta)<0\), and
\(N(\theta,h)=4\) almost everywhere outside it. By
\eqref{eq:crofton-deficit}, the shaded area equals
\(\Delta(F,\mathbb R)\).}
\label{fig:ehp-crofton-example}
\end{figure}

Translating the target in the last column of
Figure~\ref{fig:ehp-moving-circle} gives \(F(z)=z^4-2z^2\) and the
target line \(\mathbb R\). Figure~\ref{fig:ehp-crofton-example} shows
which source lines miss two intersections for this quartic.
\par\nopagebreak[4]

For fixed \(n\), the scaled form of Lemma~1 in
Fryntov and Nazarov~\cite[\S5]{FN09}
already gives \(\Delta(F,\mathbb R)\ge c_n d^{1/n}\), where
\(d\) is the diameter of the finite critical values and \(c_n>0\)
depends only on \(n\).\footnote{In their normalization the target
line is the imaginary axis and \(F(z)=z^n+\sum_{j=2}^n a_jz^{n-j}\)
with \(a_n\) real. Put \(s=\max_{2\le j\le n}|a_j|^{1/j}\).
Their lower bound for the loss is a positive multiple of \(s\),
with a constant depending only on \(n\). The critical points lie in \(|z|\le2s\),
so \(d\le2\max_{|z|\le2s}|F(z)|\le(3\cdot2^n-2)s^n\).
On each end of the inverse image, \(f_F(w)=w+O_F(w^{-1})\)
and \(f_F'(w)=1+O_F(w^{-2})\) at infinity. Thus truncating at
\(|F|=R^n\) or at \(|z|=R\) changes the length by \(O_F(R^{-1})\).}
To return to the circle we need an explicit constant that is large
enough for the comparison. For \(n\ge256\), we prove
\[
\Delta(F,\mathbb R)>\frac45(d/2)^{1/n},
\]
whenever \(d>0\).
The same estimate holds for every affine target line.

\medskip
\noindent
\(\bullet\)\enspace\textbf{Step IV. From the line deficit to the circle comparison.}
Kuznetsova and Tkachev~\cite[Theorem~4]{KT03} express the length
of an exterior level curve in terms of the Laurent coefficients of
the square root of the derivative of the normalized inverse map.
Their formula also holds at the largest critical level.
We use the weighted identity \eqref{eq:weighted-exterior-parseval}.
Applying Parseval's formula in the deformation parameter \(\beta\)
then yields \eqref{eq:outline-circle-bound}.

For the angular average at fixed \(|\beta|\), expand
\(\sqrt{\Phi_\beta'(z)}=\sum_{k\ge0} A_k(z)\beta^k\).
Parseval's formula expresses the average length as a sum of
squared coefficient norms. Each \(A_k\) is a polynomial of degree
\(nk\). Writing \([z^{nk}]A_k\) for its leading coefficient, the
exterior coefficient estimate gives
\[
\mathcal E_k
 :=\int_{\partial E_1(p)}|A_k(z)|^2\,|dz|
     -2\pi\bigl|[z^{nk}]A_k\bigr|^2\ge0.
\]
To pass to the boundary in \(\beta\), we use the local maximality
of \(p_0\) proved by Fryntov and Nazarov~\cite{FN09}.
For \(r=|\beta|\) near one, normalization of the moving circle
back to \(\mathbb T\) puts the critical values close enough to
one another for the local comparison to apply, except on arcs
of directions of total length \(O_n(1-r^2)\). A uniform bound
on length controls the contribution from these arcs. Writing
\(\mathcal L_\xi=\{w\in\C:\operatorname{Re}(\xi w)=-1/2\}\)
for \(\xi\in\mathbb T\), we obtain
\begin{equation}\label{eq:outline-circle-bound}
\ell(\partial E_1(p))
+\sum_{k\ge1}\mathcal E_k
+\frac1{2\pi}\int_0^{2\pi}
\Delta(F_{e^{i\theta}},\mathcal L_{e^{i\theta}})\,d\theta\le
2\pi-
2\pi\frac{\Gamma(-1/n)}
{\Gamma((n-1)/(2n))^2},
\end{equation}
where \(\Gamma\) is Euler's gamma function and the finitely
many exceptional boundary directions are omitted.

To obtain strict inequality, it suffices to show that the
coefficient sum and the averaged deficit in
\eqref{eq:outline-circle-bound} together exceed the difference
between its right-hand side and \(\ell(\partial E_1(p_0))\).
When the critical-value diameter is large, the averaged
line deficit alone is sufficient. At intermediate diameters,
the coefficient sum supplies the remaining contribution.
For the smaller degrees, we also use the vanishing first
variation at a maximizer: after normalizing the moving circle
back to the unit circle, maximality gives
\(\int_{\partial E_1(p)}A_1(z)\,|dz|=0\).
This imposes further constraints on the coefficients and
strengthens the lower bounds for their contribution.

For critical values lying on a very short arc, we instead
compare the circle and the line directly. Rotate a farthest
pair of critical values to \(e^{\pm i\alpha}\), and put
\(\tau=\tan(\alpha/2)\).
Straightening the arc gives a normalized monic polynomial
\(F\) whose critical values lie in \([-1,1]\) and include
both endpoints. We prove
\begin{equation}\label{eq:outline-short-arc}
\begin{aligned}
\ell(\partial E_1(p_0))-\ell(\partial E_1(p))
\ge
(2\tau)^{1/n}
\bigl[\Delta(F,\mathbb R)-\operatorname{Err}_n(\tau)\bigr].
\end{aligned}
\end{equation}
For \(n\ge3\) and \(0<\tau<1/6\), one may take
\(\operatorname{Err}_n(\tau)\) to be the right-hand side of
\eqref{eq:curvature-general} with \(R=3\), \(b_0=1/3\), and
\(h\) replaced by \(H=(U_n/(2\pi)-1)^{1/2}\), where \(U_n\) is
the right-hand side of \eqref{eq:outline-circle-bound}. The bound
\(h\le H\) follows from \eqref{eq:exterior-parseval} and
\eqref{eq:universal-length-bound}. This error bound is uniform
in \(F\) and tends to zero as \(\tau\downarrow0\) for fixed \(n\).
The required line estimate gives a positive lower bound for
\(\Delta(F,\mathbb R)\) depending only on \(n\).
Choosing the arc sufficiently short makes the error smaller
than this lower bound. Equation~\eqref{eq:outline-short-arc}
then gives strict comparison with the model.

Once the quantitative line estimates are supplied, these
comparisons exclude a maximizing polynomial with two distinct
critical values.

\medskip
\noindent
\(\bullet\)\enspace\textbf{Step V. Bounding the deficit by counting real zeros.}
Tao~\cite[Lemma~7.2 and Corollary~7.2]{Tao25} used the integrated
defect in the triangle inequality to measure the dispersion of
critical points. We identify the corresponding integral for the
zeros of a polynomial \(G\) with the mean of its line deficits
over target lines through the origin, as in
\eqref{eq:angular-root-energy}. Fryntov and Nazarov~\cite{FN09}
had used higher derivatives of real slices to detect missing real
zeros. We apply Rolle's theorem to the same slices and average
over their directions. This bounds the spread of the zeros of a
higher derivative of \(F\) in terms of the deficit on a fixed line.

Differentiation discards the lower coefficients of \(F\).
We recover bounds for them by interpolating on parallel real
slices with many zeros in a short interval. If the line deficit
were too small, these bounds would force all critical values
into a disk too small to contain a pair at the prescribed distance.

The resulting estimates are uniform for \(n\ge256\). In the
finite range \(4\le n\le255\), we sharpen the moment and slice
estimates according to the degree. For instance, when the critical
values are collinear, we also use Gon\c{c}alves' inequality
\cite{Gon50} in the form stated by Borwein, Mossinghoff, and
Vaaler~\cite[Eq.~(1.4)]{BMV07}. The remaining scalar inequalities
and parameter ranges are verified by rational or interval
arithmetic. In degree four the line estimate covers every collinear
configuration of positive diameter, including collisions. The two collision families
were studied previously in the research notes of Plasma AI~\cite{Pla26}.
Degree two follows from the choice of maximizer in Step~I.
For degree three we use Dahlke's comparison, stated in
Lemma~\ref{lem:cubic-comparison}.

\medskip
\noindent
\(\bullet\)\enspace\textbf{Step VI. Characterizing all maximizers.}
Step~I selects one maximizer. To characterize them all, we return
to the variation of a critical value. Eremenko and Hayman~\cite{EH99}
used subharmonicity to move a critical value to \(\mathbb T\) while
preserving maximality. They also showed that every finite critical
value of a global maximizer lies in the closed unit disk.
Kuznetsova and Tkachev~\cite[Corollary~5]{KT03} proved that at least
one lies on \(\mathbb T\). We fix the critical values already on
\(\mathbb T\) and compute the second variation. If a critical value
lies off \(\mathbb T\), the rational source velocity is nonconstant
and \eqref{eq:equality-second-variation} gives a strictly positive
Laplacian. This contradicts even local maximality.

The comparisons in Steps~IV and~V then force these critical
values to coincide. Since their common value has modulus
one, the characterization of polynomials with a single
finite critical value gives
\[
p(z)=(z-z_0)^n-e^{i\theta},
\qquad z_0\in\mathbb C,\quad \theta\in\mathbb R.
\]
Conversely, the lemniscate of any polynomial of this form
is a translate and rotation of \(\partial E_1(p_0)\),
and therefore has length \(\ell(\partial E_1(p_0))\).

\subsection{Structure of the paper}
Section~\ref{sec:reduction} collects the extremal and conformal
mapping results used in the proof. Section~\ref{sec:new-mobius}
carries out the deformation in Step~II, identifies the limiting
deficit in Step~III, and proves the two comparisons in Step~IV.
The estimates needed in Step~V are organized by degree.

\noindent\(\bullet\)\enspace
For \(n\ge28\), Proposition~\ref{prop:uniform-circle-transfer}
in Section~\ref{sec:circle256} reduces the problem to a lower bound
for the line deficit. Section~\ref{sec:angular-derivatives} relates
its angular mean to the roots of a derivative.
Section~\ref{sec:effective-line} proves the bound for \(n\ge256\).
Sections~\ref{sec:root-moment} and~\ref{sec:finite-interpolation}
give the estimates for \(28\le n\le255\).

\noindent\(\bullet\)\enspace
For \(9\le n\le27\), the critical values of the limiting polynomial
lie on the target line. Section~\ref{sec:collinear-gram} uses this
condition to obtain moment estimates and the line bound for
\(15\le n\le27\). Section~\ref{sec:conditional-low} proves the
line estimates for \(9\le n\le14\) and the circle comparison
for the full range \(9\le n\le27\).

\noindent\(\bullet\)\enspace
We prove the line estimates for \(n=8\) in Section~\ref{sec:octic}
and for \(n=6,7\) in Section~\ref{sec:last-degrees}.
Sections~\ref{sec:new-five-line} and~\ref{sec:new-five} treat \(n=5\).
In each of these four degrees, Proposition~\ref{prop:small-circle-transfer}
gives the bound for \(\ell(\partial E_1(p))\) from the corresponding
line estimate. Section~\ref{sec:quartic-global} treats \(n=4\) and
completes the proof of the length inequality for all \(n\ge1\).

Section~\ref{sec:equality} proves the strict variation used in Step~VI
and completes the classification of the maximizers.
Appendix~\ref{sec:certificates} gives the verification programs
for the finite estimates.

\subsection{Notation}
We denote a critical point of \(p\) by \(\zeta\). When a list is
needed, \(\zeta_1,\ldots,\zeta_{n-1}\) are the zeros of \(p'\)
counted with multiplicity, and \(v_j=p(\zeta_j)\) is the
corresponding value in the target plane. We specify the source
coordinate when a polynomial is translated or rescaled. In the
equality case of Theorem~\ref{thm:main} with \(n\ge2\), all the
\(\zeta_j\) coincide at the translation parameter \(z_0\).

We write \(D(w,r)=\{z\in\C:|z-w|<r\}\) for the open disk of
center \(w\) and radius \(r>0\), and \(\widehat\C=\C\cup\{\infty\}\).
The symbols \(ds\) and \(dA\) denote arclength and planar area
measure, respectively. For a polynomial \(F\), we write
\([z^j]F\) for its coefficient of \(z^j\) and
\(\|F\|_r=\max_{|z|\le r}|F(z)|\). For an integrable function
\(f\) on \(\mathbb T\), we put
\(\langle f\rangle=(2\pi)^{-1}\int_0^{2\pi}f(e^{it})\,dt\).
For \(y\ge0\), a quantity \(O_x(y)\) has absolute value at most
\(C_x y\), where \(C_x\) may depend only on \(x\).

A monic polynomial \(p\) of degree \(n\) is called
\emph{centered} if the centroid of its zeros, counted with
multiplicity, is at the origin, or equivalently, if its
coefficient of \(z^{n-1}\) vanishes.
If \(a_{n-1}=[z^{n-1}]p\), the zeros have mean
\(z_0=-a_{n-1}/n\). The polynomial \(q(z)=p(z+z_0)\) is
centered and satisfies \(\partial E_1(q)=\partial E_1(p)-z_0\).
Thus in comparing lengths we may center \(p\) by translation.

\subsection{Acknowledgments and AI tools disclosure} The author is grateful to Alexandre Eremenko for his encouragement,
support, and valuable comments.

This work was supported by the China Scholarship Council,
by the Young Elite Scientists Sponsorship Program for PhD Students
of the China Association for Science and Technology,
and by the Fundamental Research Funds for the Central Universities
at Xi'an Jiaotong University (Grant No.~xzy022024045).

The author began studying the Erd\H{o}s--Herzog--Piranian lemniscate
conjecture in September 2026, shortly after working on an explicit
degree bound for Sendov's conjecture. Lech Mazur's computer-assisted
work on Sendov's conjecture and Terence Tao's subsequent account of
the proof on \href{https://terrytao.wordpress.com/2026/08/12/a-digestion-of-the-proof-of-sendovs-conjecture/}{his blog} both influenced the author's approach. Inspired by this work, the author succeeded in resolving two strengthened variants of Sendov's conjecture: Borcea's 2-variance conjecture and the quadratic Tang--Zhang inequality. In the
case of lemniscates, Tao~\cite[Theorem~1.1(iv)]{Tao25} had already
shown that \(z^n-1\) attains the maximum length for every sufficiently
large degree \(n\). The author's first question was how far one could
lower the threshold in this theorem. Each improvement would leave
fewer degrees to examine individually. There was also a reason to
work in the opposite direction. A direct analysis of the lowest
degrees might reveal which estimates in the large-degree argument
were losing the most and which of those losses could be avoided.
The author therefore proposed studying the two ends of the degree
range together, with the expectation that rigorous computation would
play an essential role in the degrees in between.

A geometric question gave this strategy a more precise aim. For
\(n\geq2\), the model polynomial \(p_0(z)=z^n-1\) has a single
finite critical value, namely \(-1\). Suppose now that a polynomial
has two distinct critical values. Can their separation give a
quantitative lower bound for the amount by which its lemniscate
falls short of the model length? In pursuing this question, the
author paid particular attention to how arclength changes when the
target circle is deformed towards a line. Applying the same
M\"obius transformation to the unit circle and to the values of
the polynomial leaves the inverse image in the source plane
unchanged. The change in arclength enters when one also changes
the source coordinate so that the deformed map can again be
represented by a monic polynomial of degree \(n\).

Once the limiting line is identified with the real axis, the model
gives a concrete comparison. Fix \(R>0\) and consider the target
segment \([-R^n,R^n]\). Its inverse image under \(z\mapsto z^n\)
consists of \(n\) straight segments, each of length \(2R\). The
total length is therefore exactly \(2nR\). For a general monic
polynomial \(p\) of the same degree, the author proposed keeping
the target segment fixed and asking how large the difference
\(2nR-\ell\bigl(p^{-1}([-R^n,R^n])\bigr)\) must remain as
\(R\to\infty\) when two critical values of \(p\) are a prescribed
distance apart.

The intersection counts used by Fryntov and
Nazarov~\cite[\S5]{FN09} suggested how to approach this question.
Under the power map, the inverse image of the full real axis is
the union of \(n\) lines through the origin. Almost every line in
the source plane meets this union in exactly \(n\) points. For a
general polynomial, some lines have fewer intersections. Crofton's
formula expresses the limiting length deficit in terms of these
missing intersections, integrated over the directions and offsets
of the lines. This also clarifies what the passage back to the
circle requires. The lower bound for the deficit must be large
enough to absorb the error introduced by the deformation and
still yield a strict inequality for the original lemniscate.

The same requirement guided the author's proposal for the finite
computations. Even in a single fixed degree, the coefficients vary
continuously. The estimates to be verified therefore had to cover
whole regions of parameter space, including the configurations on
their boundaries. The quartic case illustrates a distinction that
matters here. Two critical points can coalesce, as they do in the
limit leading to \(z^4-1\). Two critical values can also coincide
while the corresponding critical points remain distinct. For
example, the quartic \(p(z)=(z^2-\sqrt{2})^2-1\) has distinct
critical points at \(z=\pm2^{1/4}\), both with critical value
\(-1\), while its third critical point, at the origin, has critical
value \(1\). Thus even the condition that all critical values lie
on the unit circle allows this second possibility. Accounting for
both kinds of collision helped the author determine which estimates
had to extend to the boundary and which configurations required
a separate argument.

ChatGPT 6.1 was used for English-language editing, proofreading, grammatical corrections, the preparation of illustrative figures, and exploratory discussions of possible approaches to selected parts of this paper. Throughout this process, the authors provided mathematical guidance, critically evaluated and corrected errors, and adjusted the research directions as necessary. The mathematical arguments and proofs in the original manuscript were independently developed, checked, and written by the author.

\section{Preliminaries}\label{sec:reduction}

We first recall the existence and normalization of a maximizer.

\begin{lem}[{\cite[Lemmas~4 and~5]{EH99}}]
\label{lem:extremizers}
For each \(n\ge2\), the map \(p\mapsto\ell(\partial E_1(p))\)
is continuous on the space of monic polynomials of degree \(n\).
It attains its maximum at a polynomial all of whose finite
critical values lie on \(\mathbb T\).
\end{lem}

We use Dahlke's results in degree three.

\begin{lem}[{\cite[Theorem~1.1 and Propositions~3.1 and~5.3]{Dah26}}]\label{lem:cubic-comparison}
Let \(p\) be a monic cubic and set \(p_0(z)=z^3-1\). Then
\(\ell(\partial E_1(p))\le\ell(\partial E_1(p_0))\).
The inequality is strict if \(p\) has two distinct critical
values on \(\mathbb T\).
\end{lem}

We shall also need local maximality near \(p_0(z)=z^n-1\).

\begin{lem}[{\cite[\S6]{FN09}}]\label{lem:local-model}
For each \(n\ge2\), there exists \(\varepsilon_n>0\) such that
every centered monic polynomial \(p\) of degree \(n\) satisfying
\(\max_{0\le j\le n-2}|[z^j](p-p_0)|<\varepsilon_n\) obeys
\(\ell(\partial E_1(p))\le\ell(\partial E_1(p_0))\).
\end{lem}
\begin{rem}
Lemma~\ref{lem:local-model} remains valid, with the same
\(\varepsilon_n\), if \(p_0(z)\) is replaced throughout
by \(z^n+v\) for any \(v\in\mathbb T\).
\end{rem}
\begin{proof}
Choose \(\lambda\in\mathbb T\) with \(\lambda^n=-v\), and put
\(q(z)=\lambda^{-n}p(\lambda z)\). Then \(q\) is centered and
monic, and
\[
q(z)-p_0(z)
=\lambda^{-n}\bigl(p(\lambda z)-((\lambda z)^n+v)\bigr).
\]
Thus the coefficient of \(z^j\) in \(q-p_0\) is
\(\lambda^{j-n}\) times the coefficient of \(z^j\) in
\(p-(z^n+v)\). Since \(|\lambda|=1\), these coefficients have
the same modulus. Also \(\partial E_1(q)=\lambda^{-1}\partial E_1(p)\),
and the same rotation takes the lemniscate of \(z^n+v\) to that
of \(p_0\). Lemma~\ref{lem:local-model} now gives the result with
the same \(\varepsilon_n\).
\end{proof}

For \(A\in\C\) and an integer \(j\ge0\), the rising factorial,
also called the Pochhammer symbol, is defined by
\[
(A)_0=1,
\qquad
(A)_j=A(A+1)\cdots(A+j-1)
\quad (j\ge1).
\]
For \(A,B,C\in\C\) with \(C\notin\{0,-1,-2,\ldots\}\),
Gauss's hypergeometric function is defined on \(|u|<1\) by
\[
{}_2F_1(A,B;C;u)
=
\sum_{j=0}^{\infty}
\frac{(A)_j(B)_j}{(C)_j\,j!}\,u^j.
\]
We shall use the following expansion as \(u\) approaches
\(1\) from below.

\begin{lem}[{\cite[Eqs.~15.10.21 and~5.5.5]{Olv26}}]
\label{lem:hypergeometric-connection}
Fix \(0<x<1\), and put \(a=(1+x)/2\) and \(b=(1-x)/2\).
Then, as \(u\uparrow1\),
\[
2\pi\,{}_2F_1(a,a;1;u)
=
2^xB\!\left(\frac12,\frac{x}{2}\right)(1-u)^{-x}+
2\pi\frac{\Gamma(-x)}{\Gamma(b)^2}
+
O_x\bigl((1-u)^{1-x}\bigr).
\]
\end{lem}

Beardon, Carne, and Ng~\cite[Lemmas~2.1 and~2.3]{BCN02}
proved the corresponding statements for polynomials normalized by
\(p(0)=0\) and \(p'\) monic. We give the argument for centered
monic polynomials, the normalization used here.

\begin{lem}\label{lem:critical-compactness}
Fix an integer \(n\ge2\) and \(0\le M<\infty\).
The centered monic polynomials \(p\) of degree \(n\)
satisfying \(|p(\zeta)|\le M\) at every critical point \(\zeta\)
form a compact set in coefficient space.

A monic polynomial \(p\) of degree \(n\) with a single
finite critical value \(v\) has the form
\[
p(z)=(z-z_0)^n+v
\]
for some \(z_0\in\C\).
\end{lem}

\begin{proof}
\emph{Step 1. Polynomials with a single critical value.}
Suppose that every critical value of \(p\) equals \(v\).
Write
\(p(z)-v=\prod_{\nu=1}^r(z-\xi_\nu)^{m_\nu}\),
where the \(\xi_\nu\) are distinct,
\(m_\nu\ge1\), and \(\sum_{\nu=1}^r m_\nu=n\).
The multiplicity of \(\xi_\nu\) as a zero of \(p'\)
is exactly \(m_\nu-1\). Every zero of \(p'\) lies among
the \(\xi_\nu\), since its image under \(p\) is \(v\).
Counting the zeros of \(p'\) with multiplicity therefore gives
\(n-1=\sum_{\nu=1}^r(m_\nu-1)=n-r\).
Thus \(r=1\), proving the asserted form.

\emph{Step 2. Boundedness of the coefficients.}
Suppose, to the contrary, that the coefficient vectors
in the stated family are unbounded. We can then choose
a sequence
\(p_k(z)=z^n+\sum_{j=0}^{n-2}a_{j,k}z^j\)
such that
\[
K_k:=\max_{0\le j\le n-2}
|a_{j,k}|^{1/(n-j)}
\longrightarrow\infty.
\]
Define
\begin{equation}\label{eq:critical-value-rescaling}
\begin{aligned}
q_k(z)
&:=K_k^{-n}p_k(K_kz)
=z^n+\sum_{j=0}^{n-2}b_{j,k}z^j,\\
b_{j,k}
&:=K_k^{j-n}a_{j,k}.
\end{aligned}
\end{equation}
By construction, \(\max_{0\le j\le n-2}|b_{j,k}|=1\).
After passing to a subsequence, the polynomials \(q_k\)
converge in coefficient space to a centered monic polynomial
\(q\) of degree \(n\), and
\begin{equation}\label{eq:critical-limit-coefficient}
\max_{0\le j\le n-2}\bigl|[z^j]q\bigr|=1.
\end{equation}

Let \(\zeta\) be any critical point of \(q\).
Continuity of polynomial roots gives critical points
\(\zeta_k\) of \(q_k\) with \(\zeta_k\to \zeta\); this also holds
when \(\zeta\) is a multiple zero of \(q'\).
Differentiating \eqref{eq:critical-value-rescaling} shows
that \(K_k\zeta_k\) is a critical point of \(p_k\).
Consequently,
\[
|q_k(\zeta_k)|
=
K_k^{-n}|p_k(K_k\zeta_k)|
\le MK_k^{-n}
\longrightarrow0.
\]
Coefficient convergence and \(\zeta_k\to \zeta\) imply \(q(\zeta)=0\).
Thus every critical value of \(q\) is zero.
Step~1 gives \(q(z)=(z-z_0)^n\), and centering forces
\(z_0=0\). This contradicts
\eqref{eq:critical-limit-coefficient}.
The coefficient vectors are therefore bounded.

\emph{Step 3. Closedness.}
Let \(p_k\) belong to the stated family and converge in
coefficient space to \(p\). The limit remains centered
and monic of degree \(n\).
For any critical point \(\zeta\) of \(p\), continuity of roots
again gives critical points \(\zeta_k\) of \(p_k\) with
\(\zeta_k\to \zeta\). Hence
\[
|p(\zeta)|
=
\lim_{k\to\infty}|p_k(\zeta_k)|
\le M.
\]
Thus \(p\) belongs to the same family.
The family is closed and bounded in the finite-dimensional
coefficient space \(\C^{n-1}\), and is therefore compact.
\end{proof}

We use the exterior area theorem in the following normalization.
\begin{lem}[{\cite[\S2.2]{Dur83}}]\label{lem:area-theorem}
Let \(f(w)=w+b_0+\sum_{j\ge1}b_jw^{-j}\) be univalent on
\(|w|>1\). Then \(\sum_{j\ge1}j|b_j|^2\le1\).
\end{lem}

Every centered monic polynomial \(p\) of degree \(n\) has
a unique normalized inverse germ \(f_p\) near infinity,
characterized by \(p(f_p(w))=w^n\) and
\(f_p(w)=w+O(w^{-1})\).
We choose the branch of \(\sqrt{f_p'}\) that tends to \(1\)
at infinity and write
\(h_j(p)=[w^{-j}]\sqrt{f_p'(w)}\) for \(j\ge1\).
These coefficients are defined using only the germ near
infinity, so no restriction on the location of the critical
values is needed.

We use the following exterior coefficient formula of
Kuznetsova and Tkachev.

\begin{lem}[{\cite[Theorem~4, Eqs.~(44)--(45)]{KT03}}]
\label{lem:exterior-parseval}
Let \(p\) be a centered monic polynomial of degree \(n\ge2\)
whose finite critical values lie in the closed unit disk.
The normalized inverse germ \(f_p\) extends univalently
to \(|w|>1\). Write
\(\sqrt{f_p'(w)}=\sum_{j\ge0}h_j(p)w^{-j}\),
choosing the square root that tends to \(1\) at infinity.
Then \(h_0(p)=1\), \(h_1(p)=0\), and
\begin{equation}\label{eq:exterior-parseval}
\ell\bigl(\partial E_1(p)\bigr)
=
2\pi\sum_{j\ge0}|h_j(p)|^2.
\end{equation}
\end{lem}

\begin{rem}\label{rem:exterior-critical-level}
Formula~\eqref{eq:exterior-parseval} includes critical
levels, where the boundary need not be a Jordan curve.
Indeed, for \(r>1\), put \(p_r(z)=r^{-n}p(rz)\).
Parametrizing \(\partial E_{r^n}(p)\) by \(f_p(re^{it})\)
and applying Parseval's identity gives
\begin{equation}\label{eq:exterior-parseval-radius}
\ell\bigl(\partial E_1(p_r)\bigr)
=
\frac1r\,\ell\bigl(\partial E_{r^n}(p)\bigr)
=
2\pi\sum_{j\ge0}|h_j(p)|^2r^{-2j}.
\end{equation}
As \(r\downarrow1\), the left-hand side converges to
\(\ell(\partial E_1(p))\) by Lemma~\ref{lem:extremizers},
while the series increases to
\(2\pi\sum_{j\ge0}|h_j(p)|^2\).
Thus \eqref{eq:exterior-parseval-radius} yields
\eqref{eq:exterior-parseval}.
At a critical level, boundary integrals are taken arc
by arc; the finitely many critical points have arclength
measure zero.
\end{rem}

We shall need the following weighted version.

\begin{cor}\label{cor:weighted-exterior-parseval}
Under the hypotheses of Lemma~\ref{lem:exterior-parseval},
let \(B\) be a nonzero polynomial of degree \(d\), and write
\(B(f_p(w))\sqrt{f_p'(w)}
=\sum_{j=-\infty}^{d}b_jw^j\).
Then
\begin{equation}\label{eq:weighted-exterior-parseval}
\int_{\partial E_1(p)}|B(z)|^2\,ds
=
2\pi\sum_{j=-\infty}^{d}|b_j|^2.
\end{equation}
\end{cor}

\begin{proof}
The function \(w^{-d}B(f_p(w))\) is bounded and analytic
on \(|w|>1\). Near infinity this follows from
\(f_p(w)=w+O(w^{-1})\); near the unit circle it follows
from \(p(f_p(w))=w^n\), which keeps \(f_p(w)\) bounded.

By \eqref{eq:exterior-parseval}, the radial values of
\(\sqrt{f_p'}\) converge in \(L^2(\mathbb T)\).
The same is true after multiplication by the bounded
analytic function \(w^{-d}B(f_p(w))\).
Parseval's identity therefore applies to this product.
On the unit circle, the factor \(w^{-d}\) has modulus one,
and the parametrization \(z=f_p(e^{it})\) gives
\(ds=|f_p'(e^{it})|\,dt\) on each regular boundary arc.
To see that each arc is traversed once, take a small open arc
\(I\subset\mathbb T\) whose closure contains no critical value of
\(p\). Choose a simply connected neighborhood of \(I\) avoiding
zero and all critical values. The \(n\) inverse branches of \(p\)
on this neighborhood agree on the exterior side of \(I\) with the
\(n\) maps \(w\mapsto f_p(u)\), where \(u\) ranges over the
branches of \(u^n=w\). Indeed, these maps are distinct by the
univalence of \(f_p\). Their restrictions to \(I\) therefore
parametrize the \(n\) arcs of \(p^{-1}(I)\), each once.
The fibres over the finitely many critical values on \(\mathbb T\)
have arclength zero. This proves \eqref{eq:weighted-exterior-parseval}.

The normalization at infinity also gives
\(b_d=[z^d]B(z)\). Thus
\eqref{eq:weighted-exterior-parseval} implies
\(\int_{\partial E_1(p)}|B(z)|^2\,ds\ge2\pi|b_d|^2\).
\end{proof}

Lengths of inverse images can be computed by summing
the contributions of the inverse branches.

\begin{lem}\label{lem:inverse-length}
Let \(p\) be a polynomial of degree \(n\ge1\), and let
\(\gamma\) be a rectifiable curve in the target plane,
with arclength counted without multiplicity. Then
\begin{equation}\label{eq:inverse-length}
\ell\bigl(p^{-1}(\gamma)\bigr)
=
\int_\gamma
\left(\sum_{p(z)=w}\frac1{|p'(z)|}\right)|dw|.
\end{equation}
For each regular value \(w\), the sum runs over its
\(n\) distinct preimages. At critical values, the integral
is understood improperly; both sides may be infinite.
\end{lem}

\begin{proof}
\emph{Step 1. Regular values.}
Let \(\mathcal V\) be the finite set of critical values
of \(p\). Every point of \(\C\setminus\mathcal V\)
has a disk neighborhood \(U\) on which all \(n\)
inverse branches \(\varphi_1,\ldots,\varphi_n\) are defined.
Differentiating \(p(\varphi_j(w))=w\) gives
\(\varphi_j'(w)=1/p'(\varphi_j(w))\).
Hence the chain rule yields
\begin{equation}\label{eq:inverse-branch-length}
\ell\bigl(\varphi_j(\gamma\cap U)\bigr)
=
\int_{\gamma\cap U}|\varphi_j'(w)|\,|dw|
=
\int_{\gamma\cap U}
\frac{|dw|}{|p'(\varphi_j(w))|}.
\end{equation}
The images of the inverse branches are disjoint.
Using a countable cover by such disks and a disjoint
measurable partition subordinate to that cover, we can
sum \eqref{eq:inverse-branch-length} without double counting.

\emph{Step 2. Passage across critical values.}
For \(\varepsilon>0\), set
\(\gamma_\varepsilon
=\gamma\setminus\bigcup_{v\in\mathcal V}D(v,\varepsilon)\).
Step~1 proves \eqref{eq:inverse-length} with
\(\gamma_\varepsilon\) in place of \(\gamma\).
As \(\varepsilon\downarrow0\), these sets increase to
\(\gamma\setminus\mathcal V\), and their inverse images
increase to \(p^{-1}(\gamma\setminus\mathcal V)\).
Monotone convergence applies to both sides.
Since \(p^{-1}(\mathcal V)\) is finite and therefore
has arclength measure zero, the limit gives
\eqref{eq:inverse-length}.
\end{proof}

\begin{rem}\label{rem:inverse-speed}
The singularities in \eqref{eq:inverse-length} are locally
integrable along lines and circles.
Indeed, suppose \(p(\zeta)=v\) and \(p-v\) vanishes to order
\(m\ge2\) at \(\zeta\). Locally,
\(p(z)-v=(z-\zeta)^m g(z)\), where \(g\) is analytic and
\(g(\zeta)\ne0\). Every inverse branch \(\varphi\) approaching
\(\zeta\) therefore satisfies
\begin{equation}\label{eq:inverse-speed}
|\varphi'(w)|
=
O\bigl(|w-v|^{-1+1/m}\bigr)
\qquad (w\to v).
\end{equation}
On a line or circle through \(v\), the distance
\(|w-v|\) is comparable to arclength distance from \(v\).
Since \(-1+1/m>-1\), the bound
\eqref{eq:inverse-speed} is integrable near \(v\).
In particular, the integral in \eqref{eq:inverse-length}
is finite for every line segment or circle.
\end{rem}

We use the Cauchy--Crofton formula with unoriented lines.

\begin{lem}[{\cite[\S4, Eq.~(2)]{Dol63}}]\label{lem:crofton-curve}
Let \(\gamma\) be a rectifiable, piecewise smooth planar
curve, with arclength counted without multiplicity.
For \(0\le\theta<\pi\) and \(h\in\R\), let
\(N_\gamma(\theta,h)\) denote the number of distinct
intersections of \(\gamma\) with
\(e^{i\theta}(\R+ih)\). Then
\begin{equation}\label{eq:crofton-curve}
\ell(\gamma)
=
\frac12\int_0^\pi\int_\R
N_\gamma(\theta,h)\,dh\,d\theta.
\end{equation}
The intersection count is understood for almost every
\((\theta,h)\). Endpoints and isolated singular points
do not affect \eqref{eq:crofton-curve}.
\end{lem}

For a centered monic polynomial \(p\), let \(f_p\) be its
normalized inverse germ at infinity, and set
\begin{equation}\label{eq:inverse-remainder-definition}
R_n(p):=[t^n]\left(\sqrt{f_p'(1/t)}-1-\tfrac12\log f_p'(1/t)\right).
\end{equation}
Here \(f_p'(1/t)\) has value \(1\) at \(t=0\). Choose the square
root and logarithm with values \(1\) and \(0\) there.

\begin{lem}\label{lem:inverse-log-identity}
For a centered monic polynomial \(p\) of degree \(n\ge3\),
\[
h_n(p)=\frac1{2n}\sum_{p'(\zeta)=0}p(\zeta)+R_n(p),
\]
where the zeros of \(p'\) are counted with multiplicity.
\end{lem}

\begin{proof}
Take \(r\) large enough that the positively oriented curve
\(C=f_p(\{|w|=r\})\) encloses all zeros of \(p'\).
Cauchy's formula and integration by parts give
\[
\begin{aligned}
[t^n]\tfrac12\log f_p'(1/t)
 &=-\frac1{4\pi i n}\int_{|w|=r}w^n\frac{f_p''(w)}{f_p'(w)}\,dw\\
 &=\frac1{4\pi i n}\int_C p(z)\frac{p''(z)}{p'(z)}\,dz .
\end{aligned}
\]
For the second equality, logarithmically differentiate
\(p'(f_p(w))f_p'(w)=nw^{n-1}\) and integrate. The term in
\(w^{n-1}\) has integral zero. At a zero \(\zeta\) of \(p'\) of
multiplicity \(m\), the residue of the last integrand is
\(mp(\zeta)\). The residue theorem and
\eqref{eq:inverse-remainder-definition} give the result.
\end{proof}

We write \(\operatorname{cap}K\) for logarithmic capacity.
For the normalization in the second assertion below, take the
inverse of the exterior map in Saff~\cite[discussion following
Eq.~(3.3)]{Saf10}.

\begin{lem}[{\cite[Theorems~1.9 and~3.5]{Saf10}}]
\label{lem:capacity-facts}
For every nonempty compact set \(K\subset\C\),
\begin{equation}\label{eq:capacity-transfinite-diameter}
\operatorname{cap}K
=
\lim_{m\to\infty}
\max_{z_1,\ldots,z_m\in K}
\left(
\prod_{1\le i<j\le m}|z_i-z_j|
\right)^{2/(m(m-1))}.
\end{equation}
If \(K\) is a nondegenerate continuum with simply connected
exterior, its exterior conformal map, normalized to fix
infinity and have positive leading coefficient, has expansion
\begin{equation}\label{eq:capacity-exterior-map}
\Psi(w)=dw+c+O(w^{-1}),
\qquad d=\operatorname{cap}K>0.
\end{equation}
\end{lem}

\begin{rem}\label{rem:capacity-similarity}
Formula~\eqref{eq:capacity-transfinite-diameter} immediately
gives
\begin{equation}\label{eq:capacity-similarity}
\operatorname{cap}(\lambda K+c)
=
|\lambda|\operatorname{cap}K
\qquad
(\lambda\in\C\setminus\{0\},\ c\in\C).
\end{equation}
Indeed, each pairwise distance is multiplied by
\(|\lambda|\).
\end{rem}

The capacity of a set containing the critical values
controls the first inverse coefficients.

\begin{lem}\label{lem:last-arc}
Let \(K\) be a nondegenerate compact continuum with simply
connected exterior, and let \(p\) be a centered monic
polynomial of degree \(n\ge3\) whose finite critical values
belong to \(K\). Write
\(f_p(w)=w+\sum_{j\ge1}a_jw^{-j}\) near infinity.
For \(d=\operatorname{cap}K\), one has
\begin{equation}\label{eq:capacity-inverse-coefficients}
\sum_{j=1}^{n-2}
j|a_j|^2d^{-2(j+1)/n}
\le1.
\end{equation}
\end{lem}

\begin{proof}
\emph{Step 1. Construction of the lift.}
Let \(\Psi\) be the exterior map in
\eqref{eq:capacity-exterior-map}.
Near infinity, choose the lift \(H\) characterized by
\begin{equation}\label{eq:capacity-lift}
p(H(w))=\Psi(w^n),
\qquad
\frac{H(w)}w\longrightarrow d^{1/n}.
\end{equation}
For the continuation argument, remove infinity from
both exterior domains.
Since the critical values of \(p\) lie in \(K\), the
inverse branches of \(p\) can be continued along every
path in the exterior of \(K\).
Properness of \(p\) prevents a lift from escaping to infinity
along a compact path.

A generator of the fundamental group of \(|w|>1\) is
mapped by \(w\mapsto\Psi(w^n)\) to \(n\) turns around
infinity. The monodromy of \(p\) around infinity is an
\(n\)-cycle, so continuation around this generator returns
each branch to itself. Thus \(H\) extends to a single-valued
holomorphic function on \(|w|>1\).

\emph{Step 2. Univalence.}
Put \(\omega=e^{2\pi i/n}\).
If \(H(w_1)=H(w_2)\), then
\eqref{eq:capacity-lift} and the injectivity of \(\Psi\)
give \(w_2=\omega^k w_1\) for some \(k\).
For \(1\le k\le n-1\), the functions \(H(w)\) and
\(H(\omega^k w)\) are distinct near infinity.
They cannot coincide at any other point: since \(p'\)
does not vanish over the exterior of \(K\), the inverse
function theorem would make them agree locally, and
hence everywhere.
Therefore \(k=0\), proving that \(H\) is univalent.

\emph{Step 3. Comparison of coefficients.}
Set \(G(w)=d^{-1/n}H(w)\).
Near infinity, the root \(\Psi(w^n)^{1/n}\) is chosen
to be asymptotic to \(d^{1/n}w\), and
\(H(w)=f_p(\Psi(w^n)^{1/n})\).
By \eqref{eq:capacity-exterior-map},
\(\Psi(w^n)^{1/n}=d^{1/n}w+O(w^{1-n})\).
Consequently,
\begin{equation}\label{eq:capacity-lift-expansion}
G(w)
=
w+\sum_{j=1}^{n-2}
a_jd^{-(j+1)/n}w^{-j}
+O(w^{1-n}).
\end{equation}
Apply Lemma~\ref{lem:area-theorem} to the univalent
function \(G\). Its first \(n-2\) Laurent coefficients,
identified in \eqref{eq:capacity-lift-expansion}, give
\eqref{eq:capacity-inverse-coefficients}.
\end{proof}

Since \(f_p'(1/t)=1-\sum_{j\ge1}ja_jt^{j+1}\) and
\(\sqrt{1-s}-1-\tfrac12\log(1-s)=O(s^2)\) at zero,
\eqref{eq:inverse-remainder-definition} shows that every monomial in \(R_n(p)\) contains at least two
inverse coefficients and has weighted degree \(n\),
where \(a_j\) has weight \(j+1\).
In particular, only \(a_1,\ldots,a_{n-3}\) occur.

Regard \(R_n\) as a polynomial in these coefficients,
and choose \(k_n\ge0\) that bounds its modulus on
\(\sum_{j=1}^{n-3}j|a_j|^2\le1\).
Such a bound exists because this set is compact.
For \(n=3\), the remainder vanishes, and we take \(k_3=0\).

\begin{cor}\label{cor:nonlinear-capacity}
Under the hypotheses of Lemma~\ref{lem:last-arc},
\begin{equation}\label{eq:nonlinear-capacity}
|R_n(p)|\le k_n\operatorname{cap}K.
\end{equation}
\end{cor}

\begin{proof}
Put \(d=\operatorname{cap}K\) and
\(q(z)=d^{-1}p(d^{1/n}z)\).
Its normalized inverse germ is
\(f_q(w)=d^{-1/n}f_p(d^{1/n}w)\), so its inverse
coefficients are \(d^{-(j+1)/n}a_j\).
By \eqref{eq:capacity-inverse-coefficients}, these
coefficients lie in the ball defining \(k_n\).
Weighted homogeneity gives \(R_n(p)=dR_n(q)\), proving
\eqref{eq:nonlinear-capacity}.
\end{proof}

We shall also need the capacity of a circular arc.

\begin{lem}[{\cite[Table~5.1]{Ran95}}]
\label{lem:arc-capacity}
For \(0<\alpha<\pi\), let
\(K_\alpha=\{e^{it}:|t|\le\alpha\}\).
Then \(\operatorname{cap}K_\alpha=\sin(\alpha/2)\).
\end{lem}

\begin{cor}\label{cor:moved-segment-capacity}
For \(\tau>0\), let
\(M_{i\tau}(w)=w/(1-i\tau w)\).
Then
\begin{equation}\label{eq:moved-segment-capacity}
\operatorname{cap}M_{i\tau}([-1,1])
=
\frac1{2\sqrt{1+\tau^2}}
\le\frac12.
\end{equation}
\end{cor}

\begin{proof}
For \(s\in[-1,1]\), direct calculation gives
\[
M_{i\tau}(s)
=
\frac{i}{2\tau}
-\frac{i}{2\tau}
\exp\bigl(2i\arctan(\tau s)\bigr).
\]
Hence
\(M_{i\tau}([-1,1])
=\frac{i}{2\tau}-\frac{i}{2\tau}K_\alpha\),
where \(\alpha=2\arctan\tau\in(0,\pi)\).
By \eqref{eq:capacity-similarity} and
Lemma~\ref{lem:arc-capacity}, its capacity is
\((2\tau)^{-1}\sin(\arctan\tau)\), which gives
\eqref{eq:moved-segment-capacity}.
\end{proof}

\section{Deforming the target circle}\label{sec:new-mobius}

Throughout this section, \(n\ge2\), and we write
\(x=1/n\), \(a=(1+x)/2\), and \(b=(1-x)/2\).

\subsection{Inverse images of lines}

We use the integral-geometric argument of Fryntov and
Nazarov~\cite[\S5]{FN09} to identify the limiting length deficit.
For a monic polynomial \(F\) of degree \(n\) and an affine
line \(\mathcal L=w_0+e^{i\psi}\R\), define
\begin{equation}\label{eq:line-deficit-definition}
\Delta(F,\mathcal L)
:=
\lim_{R\to\infty}
\left[
2nR-
\ell\left(
F^{-1}\bigl(w_0+e^{i\psi}[-R^n,R^n]\bigr)
\right)
\right].
\end{equation}
The next lemma proves that this limit exists.

\begin{lem}\label{lem:new-crofton}
For \(0\le\theta<\pi\) and \(h\in\R\), let
\(N(\theta,h)\) be the number of distinct real zeros,
as a function of \(X\), of
\[
\Im\left(
e^{-i\psi}\bigl(F(e^{i\theta}(X+ih))-w_0\bigr)
\right).
\]
Then the limit in \eqref{eq:line-deficit-definition}
is finite and satisfies
\begin{equation}\label{eq:new-crofton}
\Delta(F,\mathcal L)
=
\frac12\int_0^\pi\int_\R
\bigl(n-N(\theta,h)\bigr)\,dh\,d\theta.
\end{equation}
\end{lem}

\begin{proof}
A translation of the source centers \(F\) without changing
the truncated inverse lengths. Rotating the source and
target simultaneously, and translating the target, then
reduces the proof to \(\mathcal L=\R\), with \(F\) centered
and monic.

Write \(F(f(w))=w^n\), where
\(f(w)=w+O(w^{-1})\) near infinity.
For sufficiently large \(R_0\), the part of
\(F^{-1}(\R)\) outside
\(F^{-1}([-R_0^n,R_0^n])\) consists of the \(2n\) ends
\[
r\longmapsto f(re^{ij\pi/n}),
\qquad r\ge R_0,\quad 0\le j<2n.
\]
On these ends,
\(f'(re^{ij\pi/n})=1+O_F(r^{-2})\).
The bounded part has finite length by
Lemma~\ref{lem:inverse-length} and
Remark~\ref{rem:inverse-speed}.

For a compact curve \(\Gamma\), put
\(V_\theta(\Gamma)=\int_\R N_\Gamma(\theta,h)\,dh\).
This is the total variation of its projection onto the
normal direction \(ie^{i\theta}\).
Let
\(\Gamma_R=F^{-1}([-R^n,R^n])\) and
\(\Gamma_R^0=\{z:z^n\in[-R^n,R^n]\}\).
The difference between the projection variations of a
curved end and its corresponding straight ray is bounded by
\[
\int_{R_0}^{\infty}
\bigl|f'(re^{ij\pi/n})-1\bigr|\,dr.
\]
It follows that
\(V_\theta(\Gamma_R^0)-V_\theta(\Gamma_R)\)
is bounded uniformly in \(R\ge R_0\) and \(\theta\).

Except for the finitely many directions parallel to an
asymptotic ray, the projections of all ends are eventually
strictly monotone. For each remaining direction, the full
intersection count equals \(n\) when \(|h|\) is sufficiently
large. Moreover, the projected endpoints of the curved
and straight ends differ by \(O_F(R^{-1})\).
Consequently,
\[
\lim_{R\to\infty}
\bigl(V_\theta(\Gamma_R^0)-V_\theta(\Gamma_R)\bigr)
=
\int_\R\bigl(n-N(\theta,h)\bigr)\,dh
\]
for almost every \(\theta\).
Apply \eqref{eq:crofton-curve} to the two truncated curves
and use the uniform bound above to pass to the limit
under the angular integral. Since
\(\ell(\Gamma_R^0)=2nR\), this proves
\eqref{eq:new-crofton}.

For almost every source line, the defining real polynomial
is nonzero and has degree at most \(n\), so
\(N(\theta,h)\le n\). Thus
\eqref{eq:new-crofton} also gives \(\Delta(F,\mathcal L)\ge0\).
\end{proof}

\begin{rem}\label{rem:line-deficit-invariance}
The choice of \(w_0\in\mathcal L\) does not affect
\eqref{eq:line-deficit-definition}.
Changing \(w_0\) moves the two target endpoints by a bounded
distance. By \eqref{eq:inverse-length}, the corresponding
change in inverse length is \(O_F(R^{1-n})\), which tends
to zero.

Source translations and simultaneous translations of the
polynomial and target line preserve the deficit.
For \(s>0\), let \(F_s(z)=s^nF(z/s)\). Directly from
\eqref{eq:line-deficit-definition},
\[
\Delta(F_s,s^n\mathcal L)=s\,\Delta(F,\mathcal L).
\]
\end{rem}

\subsection{The polynomial deformation}

For \(|\beta|<1\), set
\(M_\beta(w)=w/(1-\beta w)\).
For each monic polynomial \(F\) of degree \(n\), polynomial
division gives
\[
F^2=F'Q_F+R_F,
\qquad \deg R_F\le n-2.
\]
The coefficients of \(Q_F\) and \(R_F\) are polynomial
functions of the coefficients of \(F\), since the leading
coefficient of the divisor \(F'\) is the fixed number \(n\).

\begin{lem}\label{lem:mobius-deformation}
Let \(p\) be centered and monic of degree \(n\), with all
finite critical values in the closed unit disk.
The equations
\begin{equation}\label{eq:new-flow}
\begin{aligned}
\partial_\beta F_\beta&=R_{F_\beta},
& F_0&=p,\\
\partial_\beta\Phi_\beta&=Q_{F_\beta}\circ\Phi_\beta,
& \Phi_0(z)&=z
\end{aligned}
\end{equation}
define a holomorphic family of centered monic polynomials
\(F_\beta\) on \(|\beta|<1\).

For every \(r<1\), the maps \(\Phi_\beta\), \(|\beta|\le r\),
are holomorphic and univalent on a common neighborhood
of \(\partial E_1(p)\), and
\begin{equation}\label{eq:mobius-deformation-identity}
F_\beta(\Phi_\beta(z))=M_\beta(p(z)).
\end{equation}
Critical points retain their multiplicities, and a critical
value \(v\) is carried to \(M_\beta(v)\).
The polynomial family extends across every
\(\xi\in\mathbb T\) for which \(1-\xi v\ne0\) at all
initial critical values.
\end{lem}

\begin{proof}
Local existence follows from the holomorphic differential
equations \eqref{eq:new-flow}.
The chain rule gives
\(\partial_\beta(F_\beta\circ\Phi_\beta)
=(F_\beta\circ\Phi_\beta)^2\).
The same equation and initial condition are satisfied by
\(M_\beta\circ p\), proving
\eqref{eq:mobius-deformation-identity} locally.
The degree, leading coefficient, and centering are preserved
because \(\deg R_{F_\beta}\le n-2\).

Locally, \(\Phi_\beta\) and \(M_\beta\) are biholomorphic.
Thus \eqref{eq:mobius-deformation-identity} transports
each critical point with its multiplicity, and its critical
value becomes \(M_\beta(v)\).
These transported points account for all \(n-1\) zeros
of \(F_\beta'\), counted with multiplicity.

List the critical values of \(p\) as \(v_1,\ldots,v_{n-1}\),
with critical multiplicities. The local transport just proved gives
\[
 \operatorname{Res}_z\!\left(F_\beta'(z)/n,\,w-F_\beta(z)\right)
 =\prod_{j=1}^{n-1}\bigl(w-M_\beta(v_j)\bigr).
\]
The first polynomial in the resultant is monic. Both sides are
polynomials in \(w\), so the identity of their coefficients persists
along every continuation of \(F_\beta\) that avoids the poles of
the right-hand side. On each compact parameter subdisk the critical
values are therefore bounded. Lemma~\ref{lem:critical-compactness}
bounds the coefficients, and the polynomial differential equation
continues along every path in the unit disk. The monodromy theorem
gives a single-valued family. The same argument extends it across
any boundary point at which all \(M_\beta(v_j)\) are finite.

Fix \(0<r<r'<1\), choose \(1<s<1/r'\), and put
\(U=\{z:|p(z)|<s\}\). This open set contains both
\(\partial E_1(p)\) and every critical point of \(p\).
For \(z\in U\) and \(|\beta|\le r'\) we have
\(|M_\beta(p(z))|\le s/(1-r's)\).
The coefficients of \(F_\beta\) are bounded on this parameter disk.
Since \(F_\beta\) is monic, \eqref{eq:mobius-deformation-identity}
therefore bounds \(\Phi_\beta(z)\) along every continuation.
The source equation in \eqref{eq:new-flow} continues along every
path in \(|\beta|<r'\); uniqueness and the monodromy theorem give
a single-valued holomorphic flow there for \(z\in U\).
Uniqueness along the reversed path gives injectivity in \(z\).
At each critical point, \eqref{eq:mobius-deformation-identity}
now preserves the local degree throughout \(|\beta|\le r\), since
both changes of coordinate are locally biholomorphic. These
multiplicities exhaust the \(n-1\) zeros of \(F_\beta'\).
\end{proof}

We now fix a centered monic polynomial \(p\) whose critical
values \(v_1,\ldots,v_{n-1}\), counted with multiplicity,
all belong to \(\mathbb T\).
Put \(t=1-|\beta|^2\) and
\(\mathcal C_\beta=M_\beta(\mathbb T)\).
The moving circle has equation
\begin{equation}\label{eq:moving-target-circle}
\mathcal C_\beta
=
\{w\in\C:t|w|^2-2\Re(\beta w)=1\}.
\end{equation}
Its center is \(\bar\beta/t\) and its radius is \(1/t\).

By \eqref{eq:mobius-deformation-identity},
\(\Phi_\beta(\partial E_1(p))\subset
F_\beta^{-1}(\mathcal C_\beta)\).
Over each regular target value, injectivity of
\(\Phi_\beta\) and equality of the polynomial degrees show
that all \(n\) preimages occur. Continuity then gives
\begin{equation}\label{eq:moving-inverse-curve}
\Phi_\beta\bigl(\partial E_1(p)\bigr)
=
F_\beta^{-1}(\mathcal C_\beta).
\end{equation}

To bring the moving target circle back to \(\mathbb T\), set
\begin{equation}\label{eq:normalized-deformation}
P_\beta(z):=tF_\beta(t^{-x}z)-\bar\beta.
\end{equation}
This polynomial is centered and monic. Its critical values are
\begin{equation}\label{eq:normalized-critical-values}
\frac{v_j-\bar\beta}{1-\beta v_j}\in\mathbb T.
\end{equation}
Equations~\eqref{eq:moving-inverse-curve}
and~\eqref{eq:normalized-deformation} give
\begin{equation}\label{eq:normalized-deformation-length}
\ell\bigl(\partial E_1(P_\beta)\bigr)
=
t^x\int_{\partial E_1(p)}|\Phi_\beta'(z)|\,ds.
\end{equation}
Taylor expansions in \(\beta\) will be taken before this
normalization, using the holomorphic families \(F_\beta\)
and \(\Phi_\beta\).

\subsection{The limit as the circle becomes a line}

Let
\(\mathcal P=\{\xi\in\mathbb T:\xi v_j=1
\text{ for some }j\}\).
For \(\xi\in\mathbb T\setminus\mathcal P\), write
\(F_\xi\) for the boundary value supplied by
Lemma~\ref{lem:mobius-deformation}, and set
\(\mathcal L_\xi=\{w\in\C:\Re(\xi w)=-1/2\}\).

\begin{lem}\label{lem:circle-line-limit}
For every \(\xi\in\mathbb T\setminus\mathcal P\),
\begin{equation}\label{eq:new-limit}
\lim_{r\uparrow1}(1-r^2)^{-x}
\left[
\ell\bigl(\partial E_1(p_0)\bigr)
-
\ell\bigl(\partial E_1(P_{r\xi})\bigr)
\right]
=
\Delta(F_\xi,\mathcal L_\xi).
\end{equation}
\end{lem}

\begin{proof}
\emph{Step 1. Comparison on bounded target sets.}
Put \(t=1-r^2\), and translate the polynomial and its
target circle by the moving critical value
\(M_{r\xi}(v_1)\). Thus set
\[
G_r=F_{r\xi}-M_{r\xi}(v_1),
\qquad
C_r=\mathcal C_{r\xi}-M_{r\xi}(v_1).
\]
The circle \(C_r\) passes through zero and has radius \(1/t\).
As \(r\uparrow1\), \(G_r\) converges in coefficient space to
\(G_1=F_\xi-M_\xi(v_1)\), and \(C_r\) converges locally to
\(\mathcal L_\xi-M_\xi(v_1)\), a line through zero.

For the map \(z\mapsto z^n\), the inverse image of \(C_r\)
has length \(t^{-x}\ell(\partial E_1(p_0))\).
By \eqref{eq:normalized-deformation-length}, the inverse
image under \(G_r\) has length
\(t^{-x}\ell(\partial E_1(P_{r\xi}))\).
Their difference is therefore the expression whose limit
appears in \eqref{eq:new-limit}.

Fix a sufficiently large target disk \(D(0,R)\).
The clipped inverse images lie in a common bounded source
region, since the polynomials \(G_r\) are monic and their
coefficients remain bounded.
Using \eqref{eq:moving-target-circle}, the real defining
polynomials of the moving inverse curves converge, together
with their first derivatives, on this region.
For almost every source line, all intersections with the
limiting curve are transverse and none lies on the clipping
boundary. The clipped intersection counts consequently
converge. They are bounded by \(2n\), so
\eqref{eq:crofton-curve} and dominated convergence prove
convergence of the clipped inverse lengths.

\emph{Step 2. Uniform control of the tails.}
For a monic polynomial \(G\), let
\(S_G(w)=\sum_{G(z)=w}|G'(z)|^{-1}\) at regular values.
The normalized inverse expansions give
\begin{equation}\label{eq:inverse-speed-tail}
S_{G_r}(w)
=
|w|^{x-1}+O_{p,\xi}(|w|^{-1-x}),
\end{equation}
uniformly for \(r\) sufficiently close to \(1\).
Indeed, each inverse branch has the form
\(f_{G_r}(\omega^j w^x)\), where
\(f_{G_r}'(s)=1+O_{p,\xi}(|s|^{-2})\);
summing the \(n\) inverse speeds gives
\eqref{eq:inverse-speed-tail}.

On either arm of a circle through zero of radius \(T\),
the arclength element in the radial variable \(\rho=|w|\)
is
\(d\rho/\sqrt{1-(\rho/(2T))^2}\).
For \(2T\ge R\),
\begin{equation}\label{eq:circle-tail-integral}
\int_R^{2T}
\frac{\rho^{-1-x}}{\sqrt{1-(\rho/(2T))^2}}\,d\rho
\le \frac{R^{-x}}x.
\end{equation}
To verify this, differentiate
\(-\rho^{-x}\sqrt{1-(\rho/(2T))^2}\);
its derivative is at least
\(x\rho^{-1-x}/\sqrt{1-(\rho/(2T))^2}\).

By \eqref{eq:inverse-length},
\eqref{eq:inverse-speed-tail}, and
\eqref{eq:circle-tail-integral}, the difference between
the two inverse lengths outside \(D(0,R)\) is
\(O_{p,\xi}(R^{-x})\), uniformly in \(r\).
The corresponding straight tails satisfy the same bound.
We may therefore first let \(r\uparrow1\) on the clipped
curves and then let \(R\to\infty\).
The limiting line difference is the deficit
\(\Delta(G_1,\mathcal L_\xi-M_\xi(v_1))\).
By Remark~\ref{rem:line-deficit-invariance}, this equals
\(\Delta(F_\xi,\mathcal L_\xi)\), proving \eqref{eq:new-limit}.
\end{proof}

\subsection{Averaging the deformed lengths}

Choose the square root of \(\Phi_\beta'\) holomorphically
in \(\beta\), with value \(1\) at \(\beta=0\), and write
\[
\sqrt{\Phi_\beta'(z)}
=
\sum_{k\ge0}A_k(z)\beta^k.
\]
For \(k\ge0\), define
\[
\mathcal E_k
:=
\int_{\partial E_1(p)}|A_k(z)|^2\,ds
-
2\pi \left(\frac{(a)_k}{k!}\right)^2,
\]
and put
\begin{equation}\label{eq:universal-constant}
U_n
:=
2\pi-
2\pi\frac{\Gamma(-x)}{\Gamma(b)^2}.
\end{equation}

\begin{prop}\label{prop:new-mass}
The quantities \(\mathcal E_k\) are nonnegative, and
\begin{equation}\label{eq:new-mass}
\ell\bigl(\partial E_1(p)\bigr)
+
\sum_{k\ge1}\mathcal E_k
+
\bigl\langle\Delta(F_\xi,\mathcal L_\xi)\bigr\rangle
\le U_n.
\end{equation}
The integrand may be assigned arbitrary values on the
finite exceptional set \(\mathcal P\).
\end{prop}

\begin{proof}
\emph{Step 1. The coefficient identity.}
The coefficient of \(\beta^k\) in \(\Phi_\beta\) is a
polynomial of degree at most \(nk+1\).
This follows inductively from \eqref{eq:new-flow}.
At the highest degree, only the leading term
\(z^{n+1}/n\) of \(Q_F\) contributes.
The resulting leading terms agree with those of
\(z(1-\beta z^n)^{-x}\).
After differentiation and taking the square root, they
agree with \((1-\beta z^n)^{-a}\).
Thus \(A_k\) has degree \(nk\) and leading coefficient
\((a)_k/k!\), with \(A_0=1\).

The weighted identity
\eqref{eq:weighted-exterior-parseval}, applied to \(A_k\),
now gives \(\mathcal E_k\ge0\).
Applying Parseval also in the parameter \(\beta\), and
using \eqref{eq:normalized-deformation-length}, gives,
for \(u=r^2<1\),
\begin{equation}\label{eq:new-parseval}
\sum_{k\ge0}\mathcal E_k u^k
=
(1-u)^{-x}
\left\langle
\ell\bigl(\partial E_1(P_{r\xi})\bigr)
\right\rangle
-
2\pi\,{}_2F_1(a,a;1;u).
\end{equation}
Here the exchanges of the coefficient sum and the contour
integral are justified by nonnegativity.

\emph{Step 2. Control near the exceptional directions.}
There exists \(d_*>0\), depending only on \(n\), such that
a centered monic polynomial whose critical values lie on
\(\mathbb T\) and have diameter at most \(d_*\) satisfies
the model length bound.
Otherwise, Lemma~\ref{lem:critical-compactness} would give
a coefficient limit with a single critical value
\(v\in\mathbb T\), hence a limit of the form \(z^n+v\).
Lemma~\ref{lem:local-model}, after rotation, would then
contradict the failure of the bound along the sequence.

For \(\beta=r\xi\), equation
\eqref{eq:normalized-critical-values} gives
\[
\left|
\frac{v_j-\bar\beta}{1-\beta v_j}+\bar\beta
\right|
=
\frac{1-r^2}{|1-r\xi v_j|}.
\]
Consequently, if all denominators are at least
\(2(1-r^2)/d_*\), the normalized critical values have
diameter at most \(d_*\).
The complementary directions form a union of arcs of
normalized measure \(O_n(1-r^2)\).

The polynomials \(P_\beta\) all have their critical values
on \(\mathbb T\). Lemmas~\ref{lem:critical-compactness}
and~\ref{lem:extremizers} therefore give a uniform
degree-dependent bound for their lengths.
Set
\[
D_r(\xi)
:=
(1-r^2)^{-x}
\left[
\ell\bigl(\partial E_1(p_0)\bigr)
-
\ell\bigl(\partial E_1(P_{r\xi})\bigr)
\right],
\]
and let \(D_r^-=\max\{-D_r,0\}\).
The preceding exceptional-arc estimate gives
\begin{equation}\label{eq:negative-part-average}
\langle D_r^-\rangle
=
O_n\bigl((1-r^2)^{1-x}\bigr)
\longrightarrow0.
\end{equation}

\emph{Step 3. Passage to the boundary average.}
Lemma~\ref{lem:hypergeometric-connection} and
\eqref{eq:new-parseval} imply
\begin{equation}\label{eq:mass-before-limit}
\sum_{k\ge0}\mathcal E_k r^{2k}
+
\langle D_r\rangle
=
U_n-2\pi+o(1).
\end{equation}
By Lemma~\ref{lem:circle-line-limit},
\(D_r(\xi)\to\Delta(F_\xi,\mathcal L_\xi)\) outside
\(\mathcal P\).
Apply Fatou's lemma to the positive parts, use
\eqref{eq:negative-part-average} for the negative parts,
and use monotone convergence for the nonnegative
coefficient sum in \eqref{eq:mass-before-limit}.
This gives
\[
\sum_{k\ge0}\mathcal E_k
+
\bigl\langle\Delta(F_\xi,\mathcal L_\xi)\bigr\rangle
\le U_n-2\pi.
\]
Since
\(\mathcal E_0=\ell(\partial E_1(p))-2\pi\),
this is \eqref{eq:new-mass}.
\end{proof}

\begin{cor}\label{cor:universal-length-bound}
Every monic polynomial \(q\) of degree \(n\ge2\) satisfies
\begin{equation}\label{eq:universal-length-bound}
\ell\bigl(\partial E_1(q)\bigr)
\le U_n
<
\ell\bigl(\partial E_1(p_0)\bigr)
+
2\pi-8\log2.
\end{equation}
\end{cor}

\begin{proof}
Choose a maximizer as in Lemma~\ref{lem:extremizers}
and center it by translation.
The nonnegativity of the coefficient terms and line
deficits in \eqref{eq:new-mass} gives the first bound
in \eqref{eq:universal-length-bound}.

For the strict scalar bound, put
\(J(s)=2^sB(1/2,s/2)\), using its gamma-function
continuation when \(s<0\), and let
\(k(y)=(1-e^{-2y})^{-1/2}-1\).
Substitution in the beta integral gives
\begin{equation}\label{eq:regularized-beta-integral}
J(s)
=
2^{s+1}
\left(
\frac1s+\int_0^\infty e^{-sy}k(y)\,dy
\right),
\qquad -1<s<1,\quad s\ne0.
\end{equation}
The identity first follows for \(s>0\) and then extends
to the stated range by analytic continuation.
The substitution \(q=\sqrt{1-e^{-2y}}\) gives
\(\int_0^\infty k(y)\,dy=\int_0^1(1+q)^{-1}\,dq=\log2\).

Taking \(s=\pm x\) in \eqref{eq:regularized-beta-integral},
we obtain
\begin{equation}\label{eq:beta-gap}
J(x)+J(-x)
=
\frac{4\sinh(x\log2)}x
+
4\int_0^\infty
\cosh\bigl(x(\log2-y)\bigr)k(y)\,dy
>
8\log2.
\end{equation}
The strict inequality follows from
\(\sinh(x\log2)>x\log2\), \(k(y)\ge0\), and
\(\cosh\bigl(x(\log2-y)\bigr)\ge1\).
Finally,
\(J(x)=\ell(\partial E_1(p_0))\), while the gamma
duplication formula gives
\(J(-x)=2\pi\Gamma(-x)/\Gamma(b)^2\).
Thus \eqref{eq:beta-gap} and
\eqref{eq:universal-constant} yield the second bound
in \eqref{eq:universal-length-bound}.
\end{proof}

\subsection{An exterior coefficient estimate}

The weighted Parseval formula allows us to retain a
specific coefficient below the leading term of each \(A_k\).

\begin{lem}\label{lem:exterior-coefficients}
In the notation of Proposition~\ref{prop:new-mass},
for every \(k\ge1\),
\begin{equation}\label{eq:exterior-coefficients}
\mathcal E_k
\ge
2\pi\binom{b}{k}^{\!2}|h_n(p)|^2.
\end{equation}
\end{lem}

\begin{proof}
Fix \(w\) of sufficiently large modulus. For initial points near
\(f_p(w)\), the source flow exists for \(\beta\) near zero, with
the parameter neighborhood allowed to depend on \(w\).
The chain rule gives
\(F_\beta(\Phi_\beta(f_p(w)))=M_\beta(w^n)\).
Since \(p'(f_p(w))\ne0\), uniqueness of the local inverse gives
\[
\Phi_\beta(f_p(w))
  =f_{F_\beta}\bigl(w(1-\beta w^n)^{-1/n}\bigr),
\]
where the power equals one at \(\beta=0\).
Differentiate with respect to \(w\) at fixed \(\beta\), and then
put \(\beta=u/w^n\). At \(u=0\), choose both square roots to be
\(\sqrt{f_p'(w)}\) with its normalization at infinity. As germs
in \(u\) the two sides then satisfy
\begin{equation}\label{eq:exterior-deformation-identity}
\sqrt{
\Phi_{u/w^n}'(f_p(w))f_p'(w)
}
=
(1-u)^{-a}
\sqrt{
f_{F_{u/w^n}}'\bigl(w(1-u)^{-x}\bigr)
}.
\end{equation}

To justify coefficient comparison at infinity, put \(t=1/w\).
The function \(y=t f_{F_\beta}(1/t)\) solves \(t^nF_\beta(y/t)=1\).
At \(t=0\) this equation reduces to \(y^n=1\), with simple root
\(y=1\). By the implicit function theorem, \(y\) and
\(f_{F_\beta}'(1/t)=y-t\partial_t y\) are holomorphic in
\((\beta,t)\) near \((0,0)\). Thus the right-hand side of
\eqref{eq:exterior-deformation-identity} is holomorphic in
\((u,t)\) near zero, and we may take Taylor coefficients in \(u\)
before taking Laurent coefficients at infinity. After the substitution
\(\beta=u/w^n\), a term of parameter degree \(s\) in
\(h_j(F_\beta)\) contributes at order \(w^{-j-ns}\). To contribute to
\(w^{-n}\) we must have \(j+ns=n\). Since \(h_0\equiv1\), the
only contribution comes from \((j,s)=(n,0)\).
The coefficient of \(w^{-n}\) on the right of
\eqref{eq:exterior-deformation-identity} is consequently
\((1-u)^b h_n(p)\).

Comparing the coefficients of \(u^k w^{-n}\) gives
\[
 [w^{nk-n}]\bigl(A_k(f_p(w))\sqrt{f_p'(w)}\bigr)
     =(-1)^k\binom{b}{k}h_n(p).
\]
The leading coefficient, at \(w^{nk}\), is \((a)_k/k!\).
Keeping these two terms in \eqref{eq:weighted-exterior-parseval}
and subtracting the leading contribution gives
\eqref{eq:exterior-coefficients}.
\end{proof}

Put
\[
I_0:=\frac{\Gamma(1+2b)}{\Gamma(1+b)^2},
\qquad
H:=\left(\frac{U_n}{2\pi}-1\right)^{1/2}.
\]

\begin{cor}\label{cor:full-exterior-coefficient}
Under the hypotheses of Proposition~\ref{prop:new-mass},
\begin{equation}\label{eq:full-exterior-coefficient}
\sum_{k\ge1}\mathcal E_k
\ge
2\pi(I_0-1)|h_n(p)|^2.
\end{equation}
\end{cor}

\begin{proof}
Parseval's identity for \((1-u)^b\), followed by the
beta integral, gives
\begin{equation}\label{eq:binomial-square-sum}
\sum_{k\ge0}\binom{b}{k}^{\!2}
=
\frac1{2\pi}\int_0^{2\pi}|1-e^{it}|^{2b}\,dt
=
\frac{\Gamma(1+2b)}{\Gamma(1+b)^2}
=
I_0.
\end{equation}
Summing \eqref{eq:exterior-coefficients} over \(k\ge1\)
and using \eqref{eq:binomial-square-sum} proves
\eqref{eq:full-exterior-coefficient}.
\end{proof}

\begin{rem}\label{rem:exterior-constant-normalizations}
The gamma recurrence gives the equivalent normalizations
\begin{equation}\label{eq:exterior-constant-normalizations}
H^2
=
-\frac{\Gamma(-x)}{\Gamma(b)^2}
=
\frac{n-1}{4}I_0
=
\frac{n}{B(b,b)}
\le\frac n\pi.
\end{equation}
For the last inequality, \(0<b\le1/2\) implies
\(B(b,b)=\int_0^1[t(1-t)]^{b-1}\,dt
\ge B(1/2,1/2)=\pi\).
\end{rem}

\subsection{Straightening a short arc}

Suppose that the critical values of \(p\) lie on
\(\mathbb T\) and have diameter \(0<d<\sqrt3\).
Choose a farthest pair.
After replacing \(p(z)\) by
\(e^{-i\vartheta}p(e^{i\vartheta/n}z)\) for a suitable
\(\vartheta\), we may take this pair to be
\(e^{\pm i\gamma}\), where
\(0<\gamma<\pi/3\) and \(d=2\sin\gamma\).
This replacement preserves monicity, centering, and
lemniscate length. We continue to write \(p\) for the
normalized polynomial.

Every other critical value is at distance at most \(d\)
from both endpoints.
Since \(\gamma<\pi/3\), the intersection of the two
corresponding angular neighborhoods is precisely
\(\{e^{it}:|t|\le\gamma\}\).
All critical values therefore lie on this short arc.
The arc is contained in the closed disk with center
\(c=\cos\gamma\) and radius
\(\delta=\sin\gamma=d/2\).

Put \(\tau=\tan(\gamma/2)\), and fix an \(n\)-th root
\((2i\tau)^x\) of \(2i\tau\).

\begin{lem}\label{lem:short-arc-normalization}
For the normalized polynomial \(p\) above, there is
a centered monic polynomial \(F\) whose critical values
belong to \([-1,1]\) and include both endpoints, such that
\begin{equation}\label{eq:short-arc-normalization}
p(z)
=
1+2i\tau F_{i\tau}\bigl((2i\tau)^{-x}z\bigr),
\qquad
d=\frac{4\tau}{1+\tau^2}.
\end{equation}
Here \(F_\beta\) denotes the deformation
\eqref{eq:new-flow} with initial polynomial \(F\).
\end{lem}

\begin{proof}
For \(r\in[-1,1]\),
\[
1+2i\tau M_{i\tau}(r)
=
\frac{1+i\tau r}{1-i\tau r}.
\]
This map sends \([-1,1]\) bijectively onto the short
arc from \(e^{-i\gamma}\) to \(e^{i\gamma}\).

Figure~\ref{fig:ehp-arc-normalization} shows the arc in the
coordinate \(u=(w-1)/(2i\tau)\), before \(M_{-i\tau}\) takes it
onto \([-1,1]\).
\par\nopagebreak[4]
\begin{figure}[H]
\centering
\includegraphics[width=\linewidth]{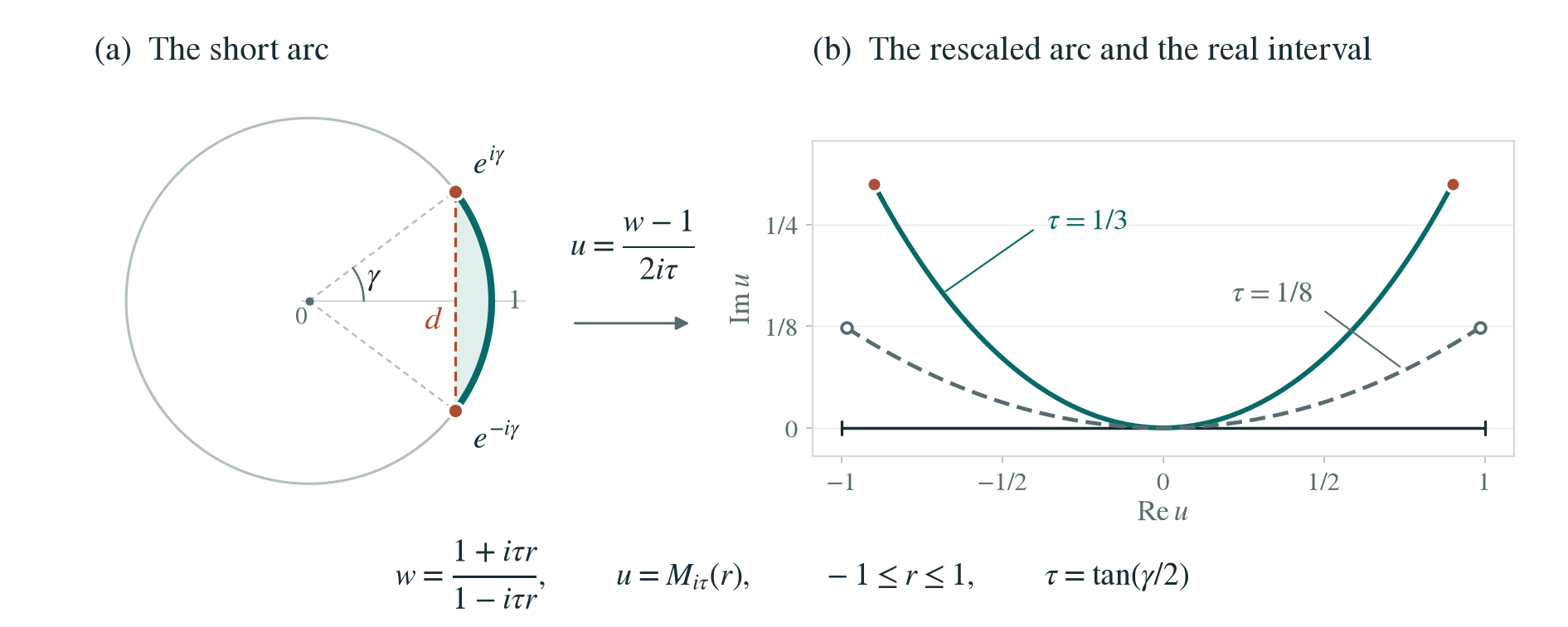}
\caption{Straightening the arc with endpoints \(e^{\pm i\gamma}\),
where \(\tau=\tan(\gamma/2)\). At left, \(\tau=1/3\) and the chord
has length \(d=2\sin\gamma\). The affine coordinate \(u=(w-1)/(2i\tau)\)
takes the arc to \(M_{i\tau}([-1,1])\). The right panel shows this image
for two values of \(\tau\), with an enlarged vertical scale.
The inverse \(M_{-i\tau}\) takes each curved arc exactly onto the black
segment \([-1,1]\).}
\label{fig:ehp-arc-normalization}
\end{figure}

Let \(\lambda=(2i\tau)^x\) and set
\(q(z)=(p(\lambda z)-1)/(2i\tau)\).
Then \(q\) is centered and monic, and its critical values
have the form \(M_{i\tau}(r_j)\), where
\(r_j\in[-1,1]\) and both endpoints occur.
Since
\[
M_{-is\tau}\bigl(M_{i\tau}(r_j)\bigr)
=
M_{i(1-s)\tau}(r_j),
\qquad 0\le s\le1,
\]
the reverse deformation encounters no pole.
In fact \(\tau<1/\sqrt3\), and the initial critical
values of \(q\) lie in the closed unit disk, so
Lemma~\ref{lem:mobius-deformation} applies throughout
this parameter segment.

Set \(F=q_{-i\tau}\).
Its critical values are the \(r_j\), with multiplicities
preserved. Uniqueness for the autonomous polynomial
equation in \eqref{eq:new-flow} gives \(F_{i\tau}=q\).
Substitution yields the polynomial identity in
\eqref{eq:short-arc-normalization}.
The relation for \(d\) follows from
\(d=2\sin\gamma\) and \(\tau=\tan(\gamma/2)\).
\end{proof}

\begin{rem}
Different choices of \((2i\tau)^x\) differ by an
\(n\)-th root of unity and only rotate the source
coordinate. The resulting lengths are unchanged.
The polynomial \(F\) may have complex coefficients,
although all its critical values are real.
\end{rem}

Whenever \(F\) has critical values in \([-1,1]\), we
use \(F_\beta\) for the deformation starting from \(F\),
and \(f_\beta\) for its normalized inverse germ,
so that \(F_\beta(f_\beta(w))=w^n\) near infinity.

\begin{cor}\label{cor:signed-short-coefficient}
Let \(n\ge3\), and suppose that \(F\) is centered and
monic with critical values in \([-1,1]\).
Assume that \(|R_n|\le\kappa\) on the coefficient ball
\(\sum_{j=1}^{n-3}j|a_j|^2\le1\).
For \(0<\tau<1\),
\begin{equation}\label{eq:signed-short-coefficient}
\max\{0,\Im h_n(F_{i\tau})\}
\le
\frac{\kappa}{2\sqrt{1+\tau^2}}
+
\frac{n-1}{2n}\tau
\le
\frac\kappa2+\frac{n-1}{2n}\tau.
\end{equation}
\end{cor}

\begin{proof}
Write the critical values of \(F\), with multiplicity,
as \(r_1,\ldots,r_{n-1}\in[-1,1]\).
Their images satisfy
\[
\Im M_{i\tau}(r_j)
=
\frac{\tau r_j^2}{1+\tau^2r_j^2}
\le\tau.
\]
Apply Lemma~\ref{lem:inverse-log-identity} to \(F_{i\tau}\).
The imaginary part of its linear critical-value term
is at most \((n-1)\tau/(2n)\).
Its critical values all lie on
\(M_{i\tau}([-1,1])\), so
\eqref{eq:nonlinear-capacity} and
\eqref{eq:moved-segment-capacity} give
\[
|R_n(F_{i\tau})|
\le
\frac{\kappa}{2\sqrt{1+\tau^2}}.
\]
These two bounds prove \eqref{eq:signed-short-coefficient}.
\end{proof}
For \(r,\tau>0\), put
\(\rho=r/\sqrt{1+\tau^2r^2}\) and
\(\vartheta(\rho)=\arcsin(\tau\rho)=\arctan(\tau r)\).
The images of the two real rays under \(M_{i\tau}\) satisfy
\begin{equation}\label{eq:curvature-arm-parametrization}
\begin{gathered}
M_{i\tau}(r)=\rho e^{i\vartheta(\rho)},
\qquad
M_{i\tau}(-r)=\rho e^{i(\pi-\vartheta(\rho))},\\
|dM_{i\tau}(\pm r)|
=
\frac{dr}{1+\tau^2r^2}
=
\frac{d\rho}{\sqrt{1-\tau^2\rho^2}}.
\end{gathered}
\end{equation}
For \(v>1\) and \(m\ge0\), write
\[
K(v,m)
:=
1+\frac1v+\frac{1+\pi m/2}{\min\{v-1,1\}}.
\]

\begin{lem}\label{lem:curvature-one-arm}
Let \(R,\tau>0\) satisfy \(\tau R\le1/2\), and put
\(s=R/\sqrt{1+\tau^2R^2}\).
For every \(v>1\) and \(\lambda\in\R\),
\begin{equation}\label{eq:curvature-one-arm}
\left|
\int_s^{1/\tau}
\frac{e^{i\lambda\vartheta(\rho)}\rho^{-1-v}}
     {\sqrt{1-\tau^2\rho^2}}\,d\rho
-
\int_R^\infty \rho^{-1-v}\,d\rho
\right|
\le
\tau s^{1-v}K(v,|\lambda|).
\end{equation}
\end{lem}

\begin{proof}
Write \(g(\rho)=\sqrt{1-\tau^2\rho^2}\).
Since \(\vartheta(\rho)\le\pi\tau\rho/2\) and
\(1-g(\rho)\le\tau\rho\), we have
\[
\left|\frac{e^{i\lambda\vartheta(\rho)}}{g(\rho)}-1\right|
\le
\left(1+\frac{\pi|\lambda|}{2}\right)
\frac{\tau\rho}{g(\rho)}.
\]
The derivative of \(-\rho^{1-v}g(\rho)\) is
\[
\frac{\rho^{-v}}{g(\rho)}
\bigl((v-1)g(\rho)^2+\tau^2\rho^2\bigr)
\ge
\min\{v-1,1\}\frac{\rho^{-v}}{g(\rho)}.
\]
Integrating from \(s\) to \(1/\tau\) bounds the error from
replacing the curved integrand by \(\rho^{-1-v}\) on this
interval by
\(\tau s^{1-v}(1+\pi|\lambda|/2)/\min\{v-1,1\}\).

It remains to account for the endpoints. Since
\(\log(R/s)=\tfrac12\log(1+\tau^2R^2)\le\tau s\)
when \(\tau R\le1/2\), and since \(\tau s<1\), we obtain
\[
\int_s^R \rho^{-1-v}\,d\rho
\le s^{-v}\log(R/s)\le\tau s^{1-v},
\qquad
\int_{1/\tau}^{\infty}\rho^{-1-v}\,d\rho
=\frac{\tau^v}{v}\le\frac{\tau s^{1-v}}{v}.
\]
Adding the two endpoint contributions gives
\eqref{eq:curvature-one-arm}.
\end{proof}

\begin{lem}\label{lem:curvature-general}
Let \(F\) be a centered monic polynomial of degree \(n\ge3\)
whose critical values lie in \([-1,1]\), and let \(F_\beta\)
denote the deformation \eqref{eq:new-flow} starting from \(F\).
Suppose that
\(R>1\), \(0<\tau<b_0<1-1/R\), and \(\tau R\le1/2\).
Put \(x=1/n\), fix an \(n\)-th root \(\sigma\) of \(2i\tau\),
and set
$
p(z)=1+2i\tau F_{i\tau}(z/\sigma).
$
Write
\[
H=\left(\frac{U_n}{2\pi}-1\right)^{1/2},
\qquad
c=\frac1{1-b_0},
\qquad
s=\frac R{\sqrt{1+\tau^2R^2}},
\qquad
h=\max\{0,\Im h_n(F_{i\tau})\}.
\]
Define
\[
\begin{aligned}
e&=\frac{\tau R}{\sqrt{1-\tau^2R^2}},
&
M&=\frac{U_n}{\pi}\left(\frac R2\right)^x
 B\!\left(\frac{1-x}{2},\frac12\right),\\
A&=\frac{HR^x(c/R)^2}
        {(2-x)\sqrt{1-(c/R)^4}},
&
B_{\mathrm{tail}}
&=\frac{2H^2c^{4x}R^{-3x}}
        {3x\bigl(1-(c/R)^4\bigr)},
\end{aligned}
\]
and put
\[
\begin{aligned}
\mathcal A
&=
\left(\sum_{q\ge2}K(q-x,q)^2s^{-2q}\right)^{1/2},\\
\mathcal T
&=
\sup_{j\ge2}
\sum_{\substack{k\ge2,\ k\equiv j\pmod n\\ j+k>n+1}}
K\bigl((j+k-1)x,|j-k|x\bigr)s^{-(j+k)x}.
\end{aligned}
\]
Then
\begin{equation}\label{eq:curvature-general}
\begin{aligned}
\Delta(F,\R)
-(2\tau)^{-x}
\left[
\ell\bigl(\partial E_1(p_0)\bigr)
-
\ell\bigl(\partial E_1(p)\bigr)
\right]
&\quad\le
M(e+e^2)
+\frac{(1+x)\tau^2R^{2+x}}{2+x}
+(4A+2B_{\mathrm{tail}})\frac{\tau}{b_0-\tau}\\
&\qquad\quad
+2hB\!\left(\frac x2,\frac12\right)\tau^{1-x}
+2\tau s^{1+x}\bigl(H^2\mathcal T+2H\mathcal A\bigr).
\end{aligned}
\end{equation}
\end{lem}

\begin{proof}
Write \(G=F_{i\tau}\), \(I=[-R,R]\), and \(J=\R\setminus I\).
The assumptions give \(c<R\). They also give
\(\tau R<R-1\), so \(1+\tau^2R^2<R^2\) and hence \(s>1\).

For a monic polynomial \(Q\) of degree \(n\), let
\(S_Q(w)=\sum_{Q(z)=w}|Q'(z)|^{-1}\) at regular values,
and put
\[
\mathcal I_Q(\gamma)
:=
\int_\gamma\bigl(S_Q(w)-|w|^{x-1}\bigr)\,|dw|
\]
for the target curves considered below.
Here \(|w|^{x-1}\) is the inverse speed of \(z\mapsto z^n\).
The local estimate \eqref{eq:inverse-speed} gives integrability
at finite critical values on lines and circles. At infinity,
the normalized inverse expansion gives
\(S_Q(w)-|w|^{x-1}=O_Q(|w|^{-1-x})\).

By \eqref{eq:line-deficit-definition} and
\eqref{eq:inverse-length}, we have
\(\Delta(F,\R)=-\mathcal I_F(\R)\).
The curve \(M_{i\tau}(\R)\) traces the circle
\(|1+2i\tau w|=1\) apart from the point \(i/\tau\).
Its inverse image under \(G\) has length
\((2\tau)^{-x}\ell(\partial E_1(p))\).
Under the power map, the corresponding length is
\((2\tau)^{-x}\ell(\partial E_1(p_0))\), by rotation
and dilation. Splitting at \(|r|=R\) therefore gives
\begin{equation}\label{eq:curvature-decomposition}
\begin{aligned}
\Delta(F,\R)
-(2\tau)^{-x}
\bigl[\ell(\partial E_1(p_0))-\ell(\partial E_1(p))\bigr]=
\bigl[\mathcal I_G(M_{i\tau}(I))-\mathcal I_F(I)\bigr]+
\bigl[\mathcal I_G(J)-\mathcal I_F(J)\bigr]+
\bigl[\mathcal I_G(M_{i\tau}(J))-\mathcal I_G(J)\bigr].
\end{aligned}
\end{equation}

\emph{Step 1. The bounded interval.}
Set \(X=F^{-1}(I)\) and \(\mathcal H=L^2(X,ds)\).
On every compact subdisk of \(|\beta|<1/R\),
\(M_\beta\circ F\) is uniformly bounded in a neighborhood
of \(X\). The continuation argument in
Lemma~\ref{lem:mobius-deformation} therefore applies on
this neighborhood. The source map \(\Phi_\beta\) is univalent
there and carries \(X\) onto \(F_\beta^{-1}(M_\beta(I))\).
Indeed, it carries each regular fibre onto the full
corresponding fibre, and continuity gives the same conclusion
at critical values.

Choose \(u(\beta)=\sqrt{\Phi_\beta'}\) with \(u(0)=1\).
This is a holomorphic \(\mathcal H\)-valued function, and
\(\|u(\beta)\|_{\mathcal H}^2
=\ell(F_\beta^{-1}(M_\beta(I)))\).
For \(\Im\beta\ne0\), the image of the real line lies on
a circle of radius \(1/(2|\Im\beta|)\), as follows from
\(\Im(1/M_\beta(r))=-\Im\beta\) for real \(r\ne0\).
For any target circle of radius \(T\) and center \(w_0\),
applying \eqref{eq:universal-length-bound} to the monic
polynomial \(T^{-1}(F_\beta(T^xz)-w_0)\) bounds the inverse
length by \(U_nT^x\).

Write \(u(\beta)=\sum_{k\ge0}u_k\beta^k\).
For \(0<\eta<1/R\), Parseval's identity and the circle
bound just obtained give
\[
\begin{aligned}
\sum_{k\ge0}\|u_k\|_{\mathcal H}^2\eta^{2k}
&=
\frac1{2\pi}\int_0^{2\pi}
\|u(\eta e^{i\phi})\|_{\mathcal H}^2\,d\phi\\
&\le
\frac{U_n(2\eta)^{-x}}{\pi}
B\!\left(\frac{1-x}{2},\frac12\right).
\end{aligned}
\]
Letting \(\eta\uparrow1/R\) gives
\(\sum_{k\ge0}\|u_k\|_{\mathcal H}^2R^{-2k}\le M\).
Put \(m=\|u_0\|_{\mathcal H}^2\).
Cauchy--Schwarz applied to the nonconstant coefficients yields
\[
\|u(i\tau)-u_0\|_{\mathcal H}
\le
\sqrt{M-m}
\left(\sum_{k\ge1}(\tau R)^{2k}\right)^{1/2}
=e\sqrt{M-m}.
\]
It follows that
\(\|u(i\tau)\|_{\mathcal H}^2-m
\le 2e\sqrt{m(M-m)}+e^2(M-m)\le M(e+e^2)\).

For the power map, the curved inverse length has density
\(|r|^{x-1}(1+\tau^2r^2)^{-(1+x)/2}\) in the original
real parameter \(r\).
Using \(1-(1+y)^{-(1+x)/2}\le(1+x)y/2\) for \(y\ge0\)
and integrating over \(I\), we obtain
\begin{equation}\label{eq:curvature-bounded-interval}
\mathcal I_G(M_{i\tau}(I))-\mathcal I_F(I)
\le
M(e+e^2)+\frac{(1+x)\tau^2R^{2+x}}{2+x}.
\end{equation}

\emph{Step 2. Varying the polynomial on the real rays.}
Write \(f_\beta=f_{F_\beta}\) and
\(h_j(\beta)=h_j(F_\beta)\), with
\(h_0(\beta)=1\) and \(h_1(\beta)=0\).
For \(|\beta|\le b_0\), every critical value of \(F_\beta\)
has modulus at most \(c\).
The inverse coefficients of \(c^{-1}F_\beta(c^xz)\) are
\(c^{-jx}h_j(\beta)\). Thus
\eqref{eq:exterior-parseval} and
\eqref{eq:universal-length-bound} give
\begin{equation}\label{eq:curvature-exterior-energy}
\sum_{j\ge2}|h_j(\beta)|^2c^{-2jx}\le H^2.
\end{equation}

Whenever \(\rho\) exceeds the moduli of all critical values
of \(F_\beta\), differentiating
\(F_\beta(f_\beta(w))=w^n\) gives
\begin{equation}\label{eq:curvature-fibre-average}
S_{F_\beta}(\rho e^{i\theta})
=
\frac{\rho^{x-1}}n
\sum_{\nu=0}^{n-1}
\left|
\sum_{j\ge0}h_j(\beta)\rho^{-jx}
e^{-ij(\theta+2\pi\nu)/n}
\right|^2.
\end{equation}
Taking \(\theta=0,\pi\) together amounts to averaging over
the \(2n\)-th roots of unity.
For \(r\ge R\) and \(0\le\nu<2n\), put
\[
V_\beta(r,\nu)
=
\sqrt{f_\beta'(r^xe^{i\pi\nu/n})}-1,
\qquad
g(\beta)
=
\sum_{q\ge1}
\frac{h_{2nq}(\beta)R^{x-2q}}{2q-x}.
\]
Give the radial variable the measure \(2r^{x-1}\,dr\)
and the index \(\nu\) normalized counting measure.
Denote the resulting \(L^2\) space by \(\mathcal K\).
Expanding the square in \eqref{eq:curvature-fibre-average}
and integrating over the two rays gives
\begin{equation}\label{eq:curvature-straight-excess}
\mathcal I_{F_\beta}(J)
=
4\Re g(\beta)+\|V_\beta\|_{\mathcal K}^2.
\end{equation}
The indices \(2nq\) in the linear term are exactly those
that survive the average over the \(2n\) directions.

By \eqref{eq:curvature-exterior-energy} and Cauchy--Schwarz,
\[
|g(\beta)|
\le
HR^x
\left(
\sum_{q\ge1}\frac{(c/R)^{4q}}{(2q-x)^2}
\right)^{1/2}
\le A.
\]
For the quadratic term, group the indices \(j\ge2\) by
their residue classes modulo \(2n\).
Within each class, the sum of the weights \((c/r)^{2jx}\)
is at most \((c/r)^{4x}/(1-(c/r)^4)\).
Applying Cauchy--Schwarz within each class and then using
\eqref{eq:curvature-exterior-energy}, we obtain
\[
\|V_\beta\|_{\mathcal K}^2
\le
2H^2\int_R^\infty
r^{x-1}\frac{(c/r)^{4x}}{1-(c/r)^4}\,dr
\le B_{\mathrm{tail}}.
\]

The coefficient series for \(g(\beta)\) and \(V_\beta\)
converge locally uniformly in their respective norms.
Indeed, \(b_0\) can be increased slightly while keeping
\((1-b_0)^{-1}<R\). Both functions are therefore holomorphic
on a neighborhood of \(|\beta|\le b_0\).
Cauchy's coefficient estimates, summed over the positive
powers of \(\tau/b_0\), give
\[
|g(i\tau)-g(0)|
\le A\frac{\tau}{b_0-\tau},
\qquad
\|V_{i\tau}-V_0\|_{\mathcal K}
\le\sqrt{B_{\mathrm{tail}}}\frac{\tau}{b_0-\tau}.
\]
Since both vectors have norm at most
\(\sqrt{B_{\mathrm{tail}}}\), the difference of their
squared norms is at most
\(2\sqrt{B_{\mathrm{tail}}}\|V_{i\tau}-V_0\|_{\mathcal K}\)
in absolute value. Hence \eqref{eq:curvature-straight-excess}
gives
\begin{equation}\label{eq:curvature-straight-bound}
\mathcal I_G(J)-\mathcal I_F(J)
\le
(4A+2B_{\mathrm{tail}})\frac{\tau}{b_0-\tau}.
\end{equation}

\emph{Step 3. Varying the target rays with \(G\) fixed.}
Write \(h_j=h_j(G)\).
Every critical value of \(G\) is \(M_{i\tau}(r_j)\)
for some \(r_j\in[-1,1]\), and therefore has modulus
at most one. Applying \eqref{eq:exterior-parseval}
directly to \(G\) gives \(\sum_{j\ge2}|h_j|^2\le H^2\)
by \eqref{eq:universal-length-bound}.

Expanding \eqref{eq:curvature-fibre-average} at \(\beta=i\tau\)
leaves only pairs \(j\equiv k\pmod n\).
For each such pair, put \(v=(j+k-1)x\) and \(m=(j-k)x\).
Its contribution to \(S_G(\rho e^{i\theta})\) is
\(h_j\overline{h_k}e^{-im\theta}\rho^{-1-v}\).
The pair \(j=k=0\) cancels in the definition of
\(\mathcal I_G\).

For a diagonal pair \(j=k\), the change on each ray,
in the original parameter \(r\ge R\), has density
\[
|h_j|^2r^{-1-v}
\left((1+\tau^2r^2)^{(v-1)/2}-1\right).
\]
This is nonpositive when \(v\le1\).
Among the off-diagonal pairs with \(v\le1\), only
\((0,n)\) and \((n,0)\) can contribute, since \(h_1=0\)
and distinct indices \(j,k\ge2\) congruent modulo \(n\)
satisfy \(j+k\ge n+4\).

The contributions of \((0,n)\) and \((n,0)\) cancel
on the two straight rays. On the curved rays their angular
factors sum to
\[
2\Re(h_ne^{-i\vartheta})
+
2\Re(h_ne^{-i(\pi-\vartheta)})
=
4\Im(h_n)\sin\vartheta.
\]
Since \(\sin\vartheta=\tau\rho\),
\eqref{eq:curvature-arm-parametrization} bounds their
combined contribution from above by
\begin{equation}\label{eq:curve-first256}
4h\tau\int_s^{1/\tau}
\frac{\rho^{x-1}}{\sqrt{1-\tau^2\rho^2}}\,d\rho
\le
2hB\!\left(\frac x2,\frac12\right)\tau^{1-x}.
\end{equation}
The last inequality follows by extending the integral
to \(0\) and substituting \(y=\tau^2\rho^2\).

For \(v>1\), factor out the phase at the straight ray.
Lemma~\ref{lem:curvature-one-arm} bounds the change
of the \((j,k)\) term on each arm by
\[
\tau s^{1+x}|h_jh_k|\,
K(v,|m|)\,s^{-(j+k)x}.
\]
For \(j,k\ge2\), the weights
\(K(v,|m|)s^{-(j+k)x}\) form a symmetric nonnegative
matrix on the admissible pairs, with row sums at most
\(\mathcal T\). Using
\(2|h_jh_k|\le |h_j|^2+|h_k|^2\), its quadratic form
at \((|h_j|)_{j\ge2}\) is at most \(H^2\mathcal T\).

If one index is zero, the other has the form \(nq\)
with \(q\ge2\). The two orders of each pair contribute
a total weight bounded by
\[
2\sum_{q\ge2}|h_{nq}|K(q-x,q)s^{-q}
\le 2H\mathcal A.
\]
Thus the contribution of all pairs with \(v>1\) on
both arms is at most
\(2\tau s^{1+x}(H^2\mathcal T+2H\mathcal A)\).
These sums converge absolutely since \(s>1\) and
\(v-1\ge1/n\) on the admissible pairs with \(v>1\).
Termwise integration follows from geometric convergence
on \(\rho\ge s\) and the integrable endpoint factor
in \eqref{eq:curvature-arm-parametrization}.

Inserting \eqref{eq:curvature-bounded-interval},
\eqref{eq:curvature-straight-bound},
\eqref{eq:curve-first256}, and the bound for \(v>1\)
into \eqref{eq:curvature-decomposition} proves
\eqref{eq:curvature-general}.
\end{proof}
\subsection{Nonlinear inverse coefficients in degrees five through eight}

We first bound the remainder in Lemma~\ref{lem:inverse-log-identity}
under the sole assumption that the Laurent coefficients satisfy the
area inequality.

\begin{lem}\label{lem:small-nonlinear}
Suppose that \(\sum_{j\ge1}j|a_j|^2\le1\). Then
\(|R_n|\le\kappa_n\) for \(5\le n\le8\), where
\[
 \kappa_5=\frac{1}{4\sqrt2},\qquad
 \kappa_6=\frac{27}{100},\qquad
 \kappa_7=\frac9{25},\qquad
 \kappa_8=\frac{53}{100}.
\]
\end{lem}

\begin{proof}
Put \(r_j=\sqrt j\,|a_j|\), so that \(\sum_{j\ge1}r_j^2\le1\).
Since \(Z(t)=\sum_{j\ge1}ja_jt^{j+1}\) vanishes to order two at
the origin, only the first three terms of
\[
 \sqrt{1-Z}-1-\tfrac12\log(1-Z)
 =\tfrac18Z^2+\tfrac5{48}Z^3+\tfrac{11}{128}Z^4+O(Z^5)
\]
can contribute in degrees at most eight. Extracting coefficients gives
\begin{equation}\label{eq:small-remainder-expansions}
\begin{aligned}
 R_5&=\tfrac12a_1a_2,\\
 R_6&=\tfrac34a_1a_3+\tfrac12a_2^2+\tfrac5{48}a_1^3,\\
 R_7&=a_1a_4+\tfrac32a_2a_3+\tfrac58a_1^2a_2,\\
 R_8&=\tfrac54a_1a_5+2a_2a_4+\tfrac98a_3^2
       +\tfrac{15}{16}a_1^2a_3+\tfrac54a_1a_2^2
       +\tfrac{11}{128}a_1^4.
\end{aligned}
\end{equation}
For \(n=5\), the first identity gives
\(\lvert R_5\rvert=r_1r_2/(2\sqrt2)\le1/(4\sqrt2)\).

\emph{Degree six.}
Set \(t=r_1\) and \(m=\tfrac14+\tfrac1{16}t^2\).
By \eqref{eq:small-remainder-expansions},
\[
 |R_6|
 \le \tfrac{\sqrt3}{4}t r_3+\tfrac14r_2^2+\tfrac5{48}t^3
 \le m(1-t^2)+\frac{3t^2}{64m}+\tfrac5{48}t^3.
\]
Here we completed the square in \(r_3\) and used
\(r_2^2+r_3^2\le1-t^2\). Convexity of \(s\mapsto(1/4+s/16)^{-1}\)
on \([0,1]\) gives \(1/m\le4-4t^2/5\). It follows that
\[
 |R_6|\le\tfrac14+\tfrac5{48}t^3-\tfrac1{10}t^4.
\]
The polynomial on the right attains its maximum at \(t=25/32\).
Its value there is
\(\tfrac14+\tfrac5{192}(25/32)^3<27/100\).

\emph{Degree seven.}
This time put \(t=r_2\), \(d=5\sqrt2/16\), and choose
\(m=(dt+\sqrt{d^2t^2+1/4})/2\). The identity
\(m(m-dt)=1/16\) gives
\(\tfrac12r_1r_4+dt r_1^2\le m(r_1^2+r_4^2)\).
The remaining term in the bound for \(|R_7|\) is
\((\sqrt6/4)t r_3\), which is at most \(mr_3^2+3t^2/(32m)\).
Thus
\[
 |R_7|\le m(1-t^2)+\frac{3t^2}{32m}
 =\frac{1+t^2/2}{4}\sqrt{1+\frac{25t^2}{32}}
   +dt\left(\frac12-\frac54t^2\right).
\]
Using \(\sqrt{1+y}\le1+y/2\) and \(7/5<\sqrt2<10/7\), we obtain
\[
 |R_7|\le
 \tfrac14+\tfrac{25}{112}t+\tfrac{57}{256}t^2
          -\tfrac{35}{64}t^3+\tfrac{25}{512}t^4
 \le \tfrac14+\tfrac{25}{112}(t+t^2)-\tfrac{125}{256}t^3,
\]
where the last inequality uses \(0\le t\le1\). The last cubic
attains its maximum at \(t=4/7\). Its value is
\(493/1372<9/25\).

\emph{Degree eight.}
Put \(t=r_1\) and \(m=3/8+5t/12\). In terms of the \(r_j\),
the last line of \eqref{eq:small-remainder-expansions} gives
\[
 |R_8|\le
 \tfrac{\sqrt5}{4}t r_5+\tfrac1{\sqrt2}r_2r_4
 +\tfrac38r_3^2+\tfrac{5\sqrt3}{16}t^2r_3
 +\tfrac58t r_2^2+\tfrac{11}{128}t^4.
\]
The quadratic \(m(m-5t/8)\) is concave on \([0,1]\), with endpoint
values \(9/64\) and \(19/144\), both greater than \(1/8\). Hence
\(\tfrac58t r_2^2+r_2r_4/\sqrt2\le m(r_2^2+r_4^2)\).
Completing the squares in \(r_3\) and \(r_5\), and using
\(m-3/8=5t/12\), now gives
\[
 |R_8|\le m(1-t^2)+\frac{5t^2}{64m}
             +\tfrac{45}{256}t^3+\tfrac{11}{128}t^4.
\]
The bound at \(t=0\) follows by continuity. Convexity gives
\(1/m\le8/3-80t/57\). Substituting this bound and using
\(t^4\le t^3\), we find
\[
 |R_8|\le \tfrac38+\tfrac5{12}t-\tfrac16t^2
                    -\tfrac{1287}{4864}t^3
 \le \tfrac38+\tfrac5{12}t-\tfrac16t^2-\tfrac14t^3.
\]
The last cubic attains its maximum at \(t=5/9\), where its value
is \(2987/5832<53/100\).
\end{proof}

\subsection{Two exterior coefficients at a maximizer}

We retain the notation of Proposition~\ref{prop:new-mass}. In
particular, \(x=1/n\), \(b=(1-x)/2\), and
\[
 I_0=\frac{\Gamma(1+2b)}{\Gamma(1+b)^2},
 \qquad H=\left(\frac{U_n}{2\pi}-1\right)^{1/2}.
\]
Write \(G(w)=\sqrt{f_p'(w)}\), with \(G(\infty)=1\), and let
\(\gamma_2(p)=[w^0]\bigl(A_2(f_p(w))G(w)\bigr)\). At a local
maximum, the first variation of length vanishes. We use this
condition to strengthen \eqref{eq:full-exterior-coefficient}.

\begin{lem}\label{lem:stationary-coefficients}
Let \(p\) be a centered monic polynomial of degree \(n\ge2\)
whose finite critical values lie on \(\mathbb T\). If \(p\) is
a local maximizer of lemniscate length among centered monic
polynomials of degree \(n\), then
\begin{equation}\label{eq:stationary-coefficients}
 \sum_{k\ge1}\mathcal E_k
 \ge 2\pi\left\{
 \left(I_0-1+\frac{b^2}{H^2}\right)|h_n(p)|^2
 +|\gamma_2(p)|^2\right\}.
\end{equation}
\end{lem}

\begin{proof}
For the normalized deformation \(P_\beta\) of
Proposition~\ref{prop:new-mass}, the length is
\((1-|\beta|^2)^x\int_{\partial E_1(p)}|\sqrt{\Phi_\beta'}|^2\,ds\).
Since \(\sqrt{\Phi_\beta'}=1+\beta A_1+O_p(|\beta|^2)\) uniformly
on \(\partial E_1(p)\), its first variation at \(\beta=0\) is
\(2\Re\bigl(\beta\int_{\partial E_1(p)}A_1\,ds\bigr)\).
Local maximality in the real and imaginary parameter directions gives
\begin{equation}\label{eq:stationary-first-moment}
 \int_{\partial E_1(p)}A_1\,ds=0.
\end{equation}

Put \(B=A_1(f_p)G\), and denote the strictly negative Laurent
parts of \(B\) and \(G\) by \(B_-\) and \(G_-\), respectively.
Use the normalized circle inner product
\(\langle B,G\rangle_2=\langle B\overline G\rangle\).
By the exterior parametrization,
\(\langle B,G\rangle_2=(2\pi)^{-1}\int_{\partial E_1(p)}A_1\,ds=0\)
by \eqref{eq:stationary-first-moment}. The coefficient identity
in the proof of Lemma~\ref{lem:exterior-coefficients} is
\begin{equation}\label{eq:stationary-selected-mode}
 [w^{n(k-1)}]\bigl(A_k(f_p(w))G(w)\bigr)
 =(-1)^k\binom bk h_n(p),\qquad k\ge1.
\end{equation}
For \(k=1\), this gives \([w^0]B=-bh_n(p)\). Since
\(G=1+G_-\), orthogonality gives
\(\langle B_-,G_-\rangle_2=bh_n(p)\). Moreover,
Lemma~\ref{lem:exterior-parseval} and Proposition~\ref{prop:new-mass}
give
\[
 \|G_-\|_2^2=\frac{\ell(\partial E_1(p))}{2\pi}-1\le H^2.
\]
It follows from Cauchy--Schwarz that
\(H\|B_-\|_2\ge b|h_n(p)|\).

By \eqref{eq:weighted-exterior-parseval}, \(\mathcal E_k/(2\pi)\)
is the sum of the squared moduli of the Laurent coefficients of
\(A_k(f_p)G\) below degree \(nk\). The coefficients in
\eqref{eq:stationary-selected-mode} contribute
\((I_0-1)|h_n(p)|^2\) after summation in \(k\). In degree \(k=1\)
we may also retain all strictly negative powers, whose contribution
is \(\|B_-\|_2^2\). In degree \(k=2\) we may retain the constant
coefficient \(\gamma_2(p)\), since the power already used there is
\(w^n\). These three contributions are disjoint and give
\eqref{eq:stationary-coefficients}.
\end{proof}

To describe how \(\gamma_2\) changes under a translation and a
dilation of the target, we need its definition for any centered
monic polynomial \(P\) of degree \(n\). Let \(F_\beta^P\) and
\(\Phi_\beta^P\) be the local solutions of \eqref{eq:new-flow}
with initial data \(P\) and the identity. For each fixed \(z\)
the solutions exist for \(\beta\) near zero. Set
\(A_2^P(z)=[\beta^2]\sqrt{(\Phi_\beta^P)'(z)}\), with the square
root equal to one at \(\beta=0\). Differentiating the flow twice
at zero shows that \(A_2^P\) is a polynomial in \(z\). We set
\(\gamma_2(P)=[w^0]\bigl(A_2^P(f_P(w))\sqrt{f_P'(w)}\bigr)\),
using the Laurent expansion at infinity. This agrees with the
preceding definition when the critical values of \(P\) lie on
\(\mathbb T\).

\begin{lem}\label{lem:affine-second-coefficient}
Let \(Q\) be centered and monic of degree \(n\ge2\), and let
\(p(z)=c+\delta Q(\delta^{-x}z)\), where \(c\in\C\) and \(\delta>0\).
Put
\[
 \eta=[\beta]h_n(F_\beta^Q),\qquad
 s_b=\frac{b(b+1)}2,\qquad
 t_b=\frac{b^2(b-1)(b+3)}4.
\]
Then
\begin{equation}\label{eq:affine-inverse}
 \sqrt{f_p'(w)}
 =(1-cw^{-n})^{-b}
   \sqrt{f_Q'\!\left(\delta^{-x}w(1-cw^{-n})^x\right)},
 \qquad h_n(p)=bc+\delta h_n(Q),
\end{equation}
and
\begin{equation}\label{eq:second-constant}
 \gamma_2(p)
 =t_b\left(c^2+\frac{2\delta c}{b}h_n(Q)\right)
   -b\delta^2\eta+s_b\delta^2h_{2n}(Q).
\end{equation}
\end{lem}

\begin{proof}
All expansions in this proof are taken near infinity, with the
branches normalized there. The inverse equation \(p(f_p(w))=w^n\)
gives
\(f_p(w)=\delta^x f_Q(\delta^{-x}w(1-cw^{-n})^x)\).
Differentiating and taking the normalized square root gives the
first identity in \eqref{eq:affine-inverse}. Its coefficient of
\(w^{-n}\) gives the second.

For \(|w|\) sufficiently large, the inverse branches starting at
\(f_p(w)\) agree for \(\beta\) near zero since \(p'(f_p(w))\ne0\).
Differentiating in \(w\) at fixed \(\beta\) and then setting
\(\beta=u/w^n\) gives \eqref{eq:exterior-deformation-identity}
for the local flow starting at \(p\). The square roots in that
identity have common value \(\sqrt{f_p'(w)}\) at \(u=0\).
The coefficient of \(u^2w^{-2n}\) on the left is \(\gamma_2(p)\).
On the right, the coefficient of \(w^{-2n}\) is
\((1-u)^{1+b}h_{2n}(p)
 +u(1-u)^b[\beta]h_n(F_\beta^p)\).
Taking its coefficient of \(u^2\) yields
\begin{equation}\label{eq:second-coefficient-identity}
 \gamma_2(p)=s_bh_{2n}(p)-b[\beta]h_n(F_\beta^p).
\end{equation}

To compute the last term, put
\(c_\beta=c/(1-\beta c)\),
\(d_\beta=\delta/(1-\beta c)^2\), and
\(\nu_\beta=\beta\delta/(1-\beta c)\). For \(\beta\) near zero,
choose the holomorphic branch of \(d_\beta^x\) with \(d_0^x=\delta^x\).
With this choice, the inverse-germ calculation in \eqref{eq:affine-inverse}
also applies to the complex dilation \(d_\beta\). We claim that
\begin{equation}\label{eq:affine-polynomial-deformation}
 F_\beta^p(z)
 =c_\beta+d_\beta F_{\nu_\beta}^Q(d_\beta^{-x}z).
\end{equation}
Indeed, let \(Y_\beta\) denote the right-hand side. The identities
\(c_\beta'=c_\beta^2\), \(d_\beta'=2c_\beta d_\beta\), and
\(\nu_\beta'=d_\beta\), together with polynomial division in
\eqref{eq:new-flow}, give
\(\partial_\beta Y_\beta\equiv Y_\beta^2\pmod{\partial_zY_\beta}\).
The derivative of the source dilation contributes a multiple of
\(z\partial_zY_\beta\) and hence vanishes in this congruence.
Since \(Y_\beta\) is centered and monic, \(\partial_\beta Y_\beta\)
has degree at most \(n-2\). It is therefore the remainder prescribed
by \eqref{eq:new-flow}. The initial condition \(Y_0=p\) proves
\eqref{eq:affine-polynomial-deformation} by uniqueness.

Reading the coefficient of \(w^{-2n}\) in
\eqref{eq:affine-inverse}, and applying its \(w^{-n}\) identity
to \eqref{eq:affine-polynomial-deformation} and differentiating
at zero, gives, respectively,
\begin{equation}\label{eq:affine-second-coefficients}
\begin{aligned}
 h_{2n}(p)
 &=s_bc^2+(b+1)\delta c\,h_n(Q)+\delta^2h_{2n}(Q),\\
 [\beta]h_n(F_\beta^p)
 &=bc^2+2\delta c\,h_n(Q)+\delta^2\eta.
\end{aligned}
\end{equation}
Substituting \eqref{eq:affine-second-coefficients} into
\eqref{eq:second-coefficient-identity} proves
\eqref{eq:second-constant}, since
\(s_b^2-b^2=t_b\) and \(s_b(b+1)-2b=2t_b/b\).
\end{proof}

\subsection{From the line deficit to the circle in degrees five through eight}

The next proposition specifies how large a line deficit suffices
in each of the four degrees. The constants \(L_n\) in its hypothesis
are listed in Table~\ref{tab:small-circle}.
\par\nopagebreak[4]
\begin{table}[H]
\centering
\begin{tabular}{ccccrr}
\toprule
\(n\) & \(L_n\) & \(\tau_0^{-1}\) & Error below
      & Bisection & Trisection \\
\midrule
5 & \(17/16\)   & 120 & \(9/10\)     &  4 &  4 \\
6 & \(73/100\)  & 200 & \(681/1000\) & 10 & 12 \\
7 & \(367/500\) & 300 & \(131/250\)  & 10 & 12 \\
8 & \(39/50\)   & 300 & \(7/10\)     &  6 &  6 \\
\bottomrule
\end{tabular}
\par\smallskip
\caption{Line deficits sufficient for the circle comparison.
The fourth column bounds the error in \eqref{eq:curvature-general}
for \(0<\tau\le\tau_0\). The last two columns count the closed
intervals used to verify \eqref{eq:small-circle-test}.}
\label{tab:small-circle}
\end{table}

\begin{prop}\label{prop:small-circle-transfer}
Fix \(n\in\{5,6,7,8\}\). Suppose that
\(\Delta(F,\mathcal L)>L_n\) for every centered monic polynomial \(F\)
of degree \(n\) whose finite critical values lie on a line \(\mathcal L\)
and have diameter two. Then every maximizing polynomial of degree
\(n\) whose finite critical values lie on \(\mathbb T\) has a
single finite critical value.
\end{prop}

\begin{proof}
Let \(p\) be a maximizer as in the statement. A source translation
allows us to assume that \(p\) is centered. Write \(2\delta\) for
the diameter of its critical values and suppose that \(\delta>0\).
We shall show that
\(\ell(\partial E_1(p))<\ell(\partial E_1(p_0))\).

\emph{Step 1. Critical values on a short arc.}
Set \(\delta_0=2\tau_0/(1+\tau_0^2)\), with \(\tau_0\) as in
Table~\ref{tab:small-circle}. If \(0<\delta\le\delta_0\), the
normalization \eqref{eq:short-arc-normalization} gives
\(0<\tau\le\tau_0\) and a centered monic polynomial \(F\) whose
critical values lie in \([-1,1]\) and include both endpoints.
The hypothesis therefore gives \(\Delta(F,\R)>L_n\).

Apply Lemma~\ref{lem:curvature-general} with \(R=3\) and \(b_0=2/5\).
Its restrictions hold since
\(\tau\le1/120<b_0<1-1/R\) and \(\tau R\le1/40\).
By Lemma~\ref{lem:small-nonlinear} and
\eqref{eq:signed-short-coefficient}, the signed coefficient \(h\)
in that lemma satisfies \(h\le\kappa_n/2+(n-1)\tau/(2n)\).

After this substitution the right-hand side of
\eqref{eq:curvature-general} is nondecreasing in \(\tau\). To see
this for the two series terms, write
\(s=R/\sqrt{1+\tau^2R^2}\). After multiplication by \(s^{1+x}\),
the powers in the row sums defining \(\mathcal T\) are \(s^{1-v}\)
with \(v=(j+k-1)/n>1\). Those inside the norm defining
\(s^{1+x}\mathcal A\) have the same form with \(v=q-x>1\).
Each increases as \(s\) decreases. The remaining terms are
nondecreasing directly. It thus suffices to evaluate the upper
bound at \(\tau_0\). The rational enclosures in
\href{https://github.com/zhangteng2000/ehp-lemniscate-verification-20261108/blob/ca8c451bfbc7312b9c57acc5e672c0b69cd601ae/verification/stationary_circle_transfer.py}{\nolinkurl{stationary_circle_transfer.py}} and
\href{https://github.com/zhangteng2000/ehp-lemniscate-verification-20261108/blob/ca8c451bfbc7312b9c57acc5e672c0b69cd601ae/verification/independent_stationary_circle_transfer.py}{\nolinkurl{independent_stationary_circle_transfer.py}} give the
strict bounds in the fourth column of Table~\ref{tab:small-circle}.
Each is less than \(L_n\), so \eqref{eq:curvature-general} gives
the desired strict length inequality.

\emph{Step 2. Lower bounds for the two coefficients.}
Put \(\delta_*=173/200\). Suppose first that
\(0<\delta\le\delta_*<\sqrt3/2\). After rotation, a pair of
critical values at distance \(2\delta\) has the form
\(e^{\pm i\alpha}=c\pm i\delta\), where
\(c=\sqrt{1-\delta^2}\) and \(0<\alpha<\pi/3\). Every critical
value is within chordal distance \(2\delta\) of both endpoints.
Since \(\alpha<\pi/3\), this places all of them on the shorter arc
between the endpoints.

Write \(p(z)=c+\delta Q(\delta^{-x}z)\). The critical values
\(q_i\) of \(Q\), counted with critical multiplicity, satisfy
\(\Re q_i\ge0\) and \(|q_i|\le1\). By Lemma~\ref{lem:arc-capacity},
the image of the shorter arc under \(w\mapsto(w-c)/\delta\)
has capacity
\[
 g=\frac{\sin(\alpha/2)}{\delta}
   =\frac1{\sqrt{2(1+c)}}.
\]
Lemma~\ref{lem:inverse-log-identity} expresses \(h_n(Q)\) as
\((2n)^{-1}\sum_iq_i+R_n(Q)\). Lemma~\ref{lem:small-nonlinear}
and \eqref{eq:nonlinear-capacity} hence give
\(\Re h_n(Q)\ge-\kappa_n g\). Lemma~\ref{lem:exterior-parseval}
and Corollary~\ref{cor:universal-length-bound} also give
\(|h_{2n}(Q)|\le H\).

We also need a lower bound for
\(\Re\eta=\Re[\beta]h_n(F_\beta^Q)\).
For \(|\beta|=1/2\), the critical values \(M_\beta(q_i)\) lie in
the image of the closed unit disk under \(M_\beta\), a disk of
radius \(4/3\). Thus \eqref{eq:nonlinear-capacity} bounds
\(|R_n(F_\beta^Q)|\) by \(4\kappa_n/3\) on this parameter circle.
Cauchy's estimate gives
\(\bigl|[\beta]R_n(F_\beta^Q)\bigr|\le8\kappa_n/3\).
Since \([\beta]M_\beta(q_i)=q_i^2\) and \(|q_i|\le1\), we obtain
\(\Re\eta\ge-K_n\), where \(K_n=(n-1)/(2n)+8\kappa_n/3\).

For real \(y\), write \([y]_+=\max(y,0)\).
Taking real parts in \eqref{eq:affine-inverse} and negative real
parts in \eqref{eq:second-constant} now gives
\begin{equation}\label{eq:small-circle-coefficients}
\begin{aligned}
 |h_n(p)|
 &\ge [bc-\kappa_ng\delta]_+,\\
 |\gamma_2(p)|
 &\ge\left[
   -t_b\left(c^2-\frac{2\kappa_ng\delta c}{b}\right)
   -\delta^2(bK_n+s_bH)\right]_+.
\end{aligned}
\end{equation}
The sign in the second estimate uses \(t_b<0\).

\emph{Step 3. The diameter intervals.}
Choose critical values \(v,v'\) with \(|v-v'|=2\delta\). For every
\(\xi\in\mathbb T\) outside the finite pole set, their images under
\(M_\xi\) are separated by
\(2\delta/(|1-\xi v|\,|1-\xi v'|)\). Rescaling the critical-value
diameter to two and applying the line hypothesis gives
\[
 \Delta(F_\xi,\mathcal L_\xi)
 \ge L_n\delta^x|1-\xi v|^{-x}|1-\xi v'|^{-x}.
\]
The circle mean of the product on the right is smallest when
\(v\) and \(v'\) are antipodal. Indeed, the superlevel sets of
each factor are arcs. For two arcs of given lengths their overlap
is smallest when their centers are antipodal. Integrating the
overlap over the two levels proves the assertion. We conclude that
\begin{equation}\label{eq:small-mean-line}
 \big\langle\Delta(F_\xi,\mathcal L_\xi)\big\rangle
 \ge L_n\delta^xJ_x,\qquad
 J_x=\big\langle|1-\xi^2|^{-x}\big\rangle
     =\frac{2^{-x}}{\pi}B\!\left(\frac{1-x}{2},\frac12\right)
     =\sum_{j\ge0}\left(\frac{(x/2)_j}{j!}\right)^2.
\end{equation}
The two expressions for \(J_x\) follow from the beta integral and
from Parseval's formula applied to \((1-u^2)^{-x/2}\), respectively.

On a closed interval \([l,u]\subset[\delta_0,\delta_*]\), put
\(c_-=\sqrt{1-u^2}\), \(c_+=\sqrt{1-l^2}\), and
\(g_+=(2(1+c_-))^{-1/2}\). Uniform lower bounds for the two
quantities in \eqref{eq:small-circle-coefficients} are
\[
 d_1=[bc_--\kappa_ng_+u]_+,\qquad
 d_2=\left[-t_b\left(c_-^2-\frac{2\kappa_n g_+u c_+}{b}\right)
                  -u^2(bK_n+s_bH)\right]_+.
\]
By \eqref{eq:stationary-coefficients}, throughout this interval
\(\sum_{k\ge1}\mathcal E_k\) is at least
\(e(l,u)=2\pi\{(I_0-1+b^2/H^2)d_1^2+d_2^2\}\).
On intervals contained in \([\delta_*,1]\) we use \(e(l,u)=0\).

The two programs in Step~1 verify, by rational interval arithmetic,
that every interval in their respective partitions satisfies
\begin{equation}\label{eq:small-circle-test}
 U_n-\ell(\partial E_1(p_0))
 < e(l,u)+L_nJ_x l^x.
\end{equation}
The first uses bisection and encloses the beta integral in
\eqref{eq:small-mean-line}. The second uses trisection and bounds
\(J_x\) below by the terms \(0\le j\le100\) of its positive series.
Each program checks that the first and last endpoints are
\(\delta_0\) and \(1\), and that successive closed intervals share
an endpoint. Thus for every \(\delta\in[\delta_0,1]\),
\eqref{eq:new-mass}, \eqref{eq:small-mean-line}, and
\eqref{eq:small-circle-test} give
\(\ell(\partial E_1(p))<\ell(\partial E_1(p_0))\).
Together with Step~1 this contradicts maximality whenever
\(\delta>0\).
\end{proof}

\Needspace{8\baselineskip}
\subsection{Real slices and their cubic derivatives}

We now turn to the line deficit itself. A slice with fewer than
\(n\) real intersections contributes to the Crofton integral
\eqref{eq:new-crofton}. One way to detect such a slice is to
differentiate its defining polynomial until only a cubic remains.
The discriminant of this cubic depends on just two coefficients
of \(F\).

\begin{lem}\label{lem:cubic-slice}
Let \(n\ge4\), let
\(F(z)=z^n+Az^{n-2}+Bz^{n-3}+O(z^{n-4})\) at infinity, and let
\(\mathcal L=w_0+e^{i\psi}\R\). For \(0\le\theta<\pi\), write
\[
\begin{aligned}
 E&=X_0+iY=e^{i(n\theta-\psi)},&
 J&=\frac{3A}{\binom n2}e^{i((n-2)\theta-\psi)},\\
 b_1&=\Im J,&
 c&=\Im\!\left(\frac{B}{\binom n3}e^{i((n-3)\theta-\psi)}\right),
 \qquad Z=E\overline J .
\end{aligned}
\]
For \(h\in\R\), set
\begin{equation}\label{eq:new-D}
 D(\theta,h)
 =(h\Im Z+Yc)^2+\frac4{27}Yb_1^3
   +h^2\left(4hX_0c-\frac43b_1^2
                  +4h^2(\Re Z-h^2)\right).
\end{equation}
Then \(\Delta(F,\mathcal L)\) is at least the area, with respect to
\(d\theta\,dh\), of
\(\{(\theta,h)\in[0,\pi)\times\R:D(\theta,h)>0\}\).
\end{lem}

\begin{proof}
The intersections counted in \eqref{eq:new-crofton} are the
distinct real zeros of
\(S(X)=\Im(e^{-i\psi}(F(e^{i\theta}(X+ih))-w_0))\).
Differentiating \(n-3\) times removes all terms below degree
\(n-3\) and gives
\[
 \frac6{n!}S^{(n-3)}(X)
 =YX^3+3hX_0X^2+(b_1-3h^2Y)X+c+h\Re J-h^3X_0.
\]
Substitution into the cubic discriminant, using \(X_0^2+Y^2=1\),
gives \(-27D(\theta,h)\). Apart from the finitely many directions
with \(Y=0\), this is a cubic of degree three. If \(D(\theta,h)>0\),
its discriminant is negative, so it has exactly one real zero.
Repeated application of Rolle's theorem then gives at most
\(1+(n-3)=n-2\) distinct real zeros of \(S\). On this set the
integrand \((n-N(\theta,h))/2\) in \eqref{eq:new-crofton} is
at least one. It is nonnegative for almost every other line.
Integration proves the claim.
\end{proof}

We finish by giving a test for \(D>0\) throughout a rectangular
cell. Put \(f=3|A|/\binom n2\) and \(g=|B|/\binom n3\). Write
\(A=|A|e^{i\phi_A}\) and \(B=|B|e^{i\phi_B}\), choosing an
arbitrary phase for a zero coefficient. After replacing the angular
coordinate by \(\theta-\phi_A/2\) and calling it \(\theta\) again,
put
\[
 \alpha=\tfrac n2\phi_A-\psi,\qquad
 \beta=\phi_B-\tfrac32\phi_A,\qquad
 \phi_0=n\theta+\alpha,\quad
 \phi_1=(n-2)\theta+\alpha,\quad
 \phi_2=(n-3)\theta+\alpha+\beta.
\]
Then \(E=e^{i\phi_0}\), \(J=fe^{i\phi_1}\),
\(c=g\sin\phi_2\), and \(Z=fe^{2i\theta}\).
The translation of the angular interval has no effect on the area
in Lemma~\ref{lem:cubic-slice}, since
\(D(\theta+\pi,-h)=D(\theta,h)\). In these coordinates,
\eqref{eq:new-D} becomes
\begin{equation}\label{eq:cubic-phase-expansion}
\begin{aligned}
 D={}&g^2\sin^2\phi_0\sin^2\phi_2
       +\tfrac4{27}f^3\sin\phi_0\sin^3\phi_1\\
    &+2hfg\sin(2\theta)\sin\phi_0\sin\phi_2
       +h^2f^2\left(\sin^2(2\theta)-\tfrac43\sin^2\phi_1\right)\\
    &+4h^3g\cos\phi_0\sin\phi_2
       +4h^4f\cos(2\theta)-4h^6.
\end{aligned}
\end{equation}

On a closed rectangular cell in \((\theta,\alpha,\beta,f,g,h)\),
let \(f_*,g_*\) be the upper endpoints of the \(f,g\) intervals,
and put \(r=\max|h|\), \(v=\min|h|\), and
\(W=g_*^2+4f_*^3/27\). The following constants bound the
corresponding second derivatives from above throughout the cell:
\begin{equation}\label{eq:cubic-hessian}
\begin{aligned}
 U_\theta
 &=(4n-6)^2W+2(2n-1)^2f_*g_*r
       +\left(8+\tfrac83(n-2)^2\right)f_*^2r^2\\
 &\qquad +(4n-6)^2g_*r^3+16f_*r^4,\\
 U_\alpha&=16W+8f_*g_*r+\tfrac{16}{3}f_*^2r^2+16g_*r^3,\\
 U_\beta&=4g_*^2+2f_*g_*r+4g_*r^3,\\
 U_f&=2r^2+\tfrac89f_*,
 \qquad U_g=2,\\
 U_h&=\max\{0,48f_*r^2+24g_*r+2f_*^2-120v^4\}.
\end{aligned}
\end{equation}
Indeed, the second derivative of a product of sines and cosines
is bounded in absolute value by the square of the sum of their
absolute frequencies. Applying this term by term in
\eqref{eq:cubic-phase-expansion} gives the angular bounds. In the
coefficient of \(h^2\), differentiating the two squares directly
gives \(8+\tfrac83(n-2)^2\) for the \(\theta\)-derivative.
Differentiation in \(f\) and \(g\) gives \(U_f\) and \(U_g\).
For the \(h\)-derivative the term \(-4h^6\) contributes
\(-120h^4\le-120v^4\), which gives the last bound in
\eqref{eq:cubic-hessian}.

Let \(w_j\) be the width of the cell in coordinate \(j\), and let
\(\mathcal I D\) be the multilinear interpolant of the vertex
values of \(D\). In one variable, the difference between the
linear interpolant and the function is bounded above by \(Uw^2/8\)
when the second derivative is bounded above by \(U\) on an
interval of width \(w\). Applying this estimate one coordinate
at a time gives
\[
 \mathcal I D-D\le\frac18\sum_j U_jw_j^2,
\]
since each interpolation operator preserves order and constants.
The interpolant is a convex combination of the vertex values.
Thus the strict test
\[
 \min_{\text{vertices}}D>\frac18\sum_jU_jw_j^2
\]
certifies \(D>0\) throughout the cell. For fixed
\(\alpha,\beta,f,g\), the areas of certified \((\theta,h)\)-rectangles
with disjoint interiors may therefore be added in
Lemma~\ref{lem:cubic-slice}.

\section{Angular averages for polynomial inverse images of lines}
\label{sec:angular-derivatives}

\subsection{The mean deficit and the roots}

Let \(G\) be a monic polynomial of degree \(m\), with roots
\(\xi_1,\ldots,\xi_m\) counted with multiplicity. We write
\(D_G=\max_{i,j}|\xi_i-\xi_j|\) and denote the mean deficit over
target lines through the origin by
\[
 U(G)=\frac1\pi\int_0^\pi\Delta(G,e^{i\phi}\R)\,d\phi.
\]
Both \(U(G)\) and \(D_G\) are invariant under translation of the
source. We first express \(U(G)\) directly in terms of the roots.

For a finite positive measure \(\mu=\sum_i w_i\delta_{\xi_i}\)
of finite support, put
\[
 \mathcal U(\mu)
 =\frac1\pi\int_{\C}
 \left(\sum_i\frac{w_i}{|z-\xi_i|}
       -\left|\sum_i\frac{w_i}{z-\xi_i}\right|\right)dA(z).
\]
The integrand is nonnegative and locally integrable. At infinity
the arguments of the nonzero summands \(w_i/(z-\xi_i)\) differ
from \(-\arg z\) by \(O_\mu(|z|^{-1})\). Taking their projections
in that direction bounds the integrand by \(O_\mu(|z|^{-3})\).
In particular, \(\mathcal U(\mu)\) is finite.

Tao~\cite[Lemmas~7.1--7.2]{Tao25} used the integrated
defect in the triangle inequality to estimate the dispersion of
critical points. For the roots of \(G\), this defect has the exact
interpretation
\begin{equation}\label{eq:angular-root-energy}
 U(G)=\mathcal U\!\left(\sum_{j=1}^m\delta_{\xi_j}\right).
\end{equation}
Indeed, fix a source line \(z=e^{i\theta}(X+ih)\) that misses the
roots of \(G\). As the target angle \(\phi\) runs from \(0\) to
\(\pi\), the integral of the intersection count is the total variation
of \(\arg G(z)\) along the source line. This follows from the
one-dimensional area formula applied to \(\arg G(z)\) modulo \(\pi\).
The density of this variation is
\(\bigl|\Im(e^{i\theta}G'(z)/G(z))\bigr|\). Each factor \(z-\xi_j\)
has total argument variation \(\pi\) on the same line. Thus, on
averaging \eqref{eq:new-crofton} in \(\phi\), the difference inside
the integral over \(X\) is
\[
 \sum_{j=1}^m
 \left|\Im\frac{e^{i\theta}}{z-\xi_j}\right|
 -
 \left|\Im\left(e^{i\theta}\sum_{j=1}^m
                         \frac1{z-\xi_j}\right)\right|.
\]
This difference is nonnegative. Tonelli's theorem, the change of
variables \(dX\,dh=dA(z)\), and the identity
\(\int_0^\pi|\Im(e^{i\theta}v)|\,d\theta=2|v|\) give
\eqref{eq:angular-root-energy}.

The constant for two roots of equal weight will occur several times.
Set
\[
 c_*=\frac1{2\pi}\left(\int_0^\pi\sqrt{\sin t}\,dt\right)^2
     =\frac{B(3/4,1/2)^2}{2\pi}.
\]
We shall use \(c_*>913/1000\). For an explicit check, take \(t=3/2\)
in \eqref{eq:split-beta-series} and retain the terms \(0\le j\le7\)
in both positive series. Using \(2^{-3/4}>22/37\) and
\(2^{-1/2}>70/99\) gives \(B(3/4,1/2)>599/250\). The stated
bound follows from \(\pi<22/7\).

\begin{lem}\label{lem:weighted-root-mean}
Let \(G\) be centered and monic of degree \(m\ge3\), and put
\(t=(m-2)/2\). Then
\begin{equation}\label{eq:root-mean256}
 U(G)\ge\frac{c_*m}{2(m-1)}\sum_{j=1}^m|\xi_j|,
 \qquad
 U(G)\ge\frac{2t}{1+t+\sqrt t}\,D_G.
\end{equation}
\end{lem}

\begin{proof}
\emph{Step 1. Two roots of equal weight.}
The triangle inequality makes \(\mathcal U\) superadditive in
\(\mu\), and its definition gives positive homogeneity in the
weights. A similarity of the support with ratio \(d\) multiplies
\(\mathcal U\) by \(d\). We claim that
\begin{equation}\label{eq:angular-equal-pair}
 \mathcal U(\delta_a+\delta_b)=c_*|a-b|.
\end{equation}

By translation it suffices to use the two roots of \(z^2+A\).
For the target line \(e^{i\phi}\R\), the discriminant of the real
slice on \(z=e^{i\theta}(X+ih)\), divided by four, is
\(h^2-\sin(2\theta-\phi)\Im(e^{-i\phi}A)\).
For almost every \(\theta\), the slice loses two real intersections
precisely where this expression is negative. Write
\(I=\int_0^\pi\sqrt{\sin t}\,dt\) and \([y]_+=\max(y,0)\).
Integration first in \(h\) and then in \(\theta\) gives
\[
 \Delta(z^2+A,e^{i\phi}\R)
 =2\int_0^\pi
   \sqrt{[\sin(2\theta-\phi)\Im(e^{-i\phi}A)]_+}\,d\theta
 =I\sqrt{|\Im(e^{-i\phi}A)|}.
\]
Averaging in \(\phi\) gives \(U(z^2+A)=I^2\sqrt{|A|}/\pi\).
The roots are \(2\sqrt{|A|}\) apart, so
\eqref{eq:angular-root-energy} proves \eqref{eq:angular-equal-pair}.

Each root occurs in \(m-1\) unordered pairs. By superadditivity,
homogeneity, and \eqref{eq:angular-equal-pair},
\[
 U(G)\ge\frac{c_*}{m-1}\sum_{i<j}|\xi_i-\xi_j|.
\]
Since \(\sum_j\xi_j=0\), the triangle inequality gives
\(\sum_j|\xi_i-\xi_j|\ge|\sum_j(\xi_i-\xi_j)|=m|\xi_i|\).
Summing over \(i\) proves the first bound in \eqref{eq:root-mean256}.

\emph{Step 2. Two roots at maximal distance.}
Relabel the roots so that \(|\xi_1-\xi_2|=D_G\), and write
\(\mathcal C(t)=\mathcal U(t\delta_0+\delta_1)\).
The full root measure is the average, over \(3\le j\le m\), of
\[
 (\delta_{\xi_1}+t\delta_{\xi_j})
 +(\delta_{\xi_2}+t\delta_{\xi_j}).
\]
Superadditivity and similarity therefore give
\[
 U(G)\ge\frac{\mathcal C(t)}{m-2}
 \sum_{j=3}^m\bigl(|\xi_1-\xi_j|+|\xi_2-\xi_j|\bigr)
 \ge\mathcal C(t)D_G.
\]
It remains to estimate \(\mathcal C(t)\).

Put \(r=|z|\), \(q=|z-1|\), \(u=r+q\), \(v=r-q\), and \(A=r+tq\).
After multiplication by \(rq\), the integrand defining
\(\pi\mathcal C(t)\) becomes
\(A-\sqrt{A^2-t(1-v^2)}\). The two half-planes together have
area element
\(2rq\,du\,dv/\sqrt{(u^2-1)(1-v^2)}\), with \(u>1\) and
\(-1<v<1\). Rationalizing the difference of square roots bounds it
below by \(t(1-v^2)/(2A)\). Since
\(2A=(1+t)(u-\eta v)\), where \(\eta=(t-1)/(t+1)\), integration
in \(v\) gives
\begin{equation}\label{eq:angular-unequal-pair}
 \mathcal C(t)\ge\frac{2t}{1+t}
 \int_1^\infty
 \frac{du}{(u+\sqrt{u^2-\eta^2})\sqrt{u^2-1}}.
\end{equation}
Here we used
\(\int_{-1}^1\sqrt{1-v^2}\,(a+bv)^{-1}\,dv
=\pi/(a+\sqrt{a^2-b^2})\) for \(a>|b|\), obtained by putting
\(v=\cos\theta\) and then \(s=\tan(\theta/2)\).

In \eqref{eq:angular-unequal-pair}, put \(u=\cosh y\) and
\(\sigma=\sqrt{1-\eta^2}=2\sqrt t/(1+t)\). The inequality
\(\sqrt{\sinh^2y+\sigma^2}\le\sinh y+\sigma\) bounds the integral
below by
\[
 \int_0^\infty\frac{dy}{e^y+\sigma}
 =\int_0^1\frac{ds}{1+\sigma s}
 \ge\frac1{1+\sigma/2},
\]
where the last inequality is Jensen's inequality. This gives
\(\mathcal C(t)\ge2t/(1+t+\sqrt t)\) and proves the second bound
in \eqref{eq:root-mean256}.
\end{proof}

\subsection{Sampling the target angle after differentiation}

Differentiating a real slice \(k\) times changes the target angle
by \(k\theta\), where \(\theta\) is the direction of the source
line. The next estimate compares the resulting deficit with the
full angular mean for the derivative.

\begin{lem}\label{lem:sampled-derivative}
Let \(F\) be centered and monic of degree \(n\), let
\(1\le k\le n-2\), and put \(m=n-k\) and \(G=m!F^{(k)}/n!\).
Let \(R\) be the largest modulus of a root of \(G\). With
\[
 \alpha_k=\frac1k\csc\frac{\pi}{2k}-\frac2\pi,
 \qquad
 \beta_k=\frac2\pi-\frac1k\cot\frac{\pi}{2k},
\]
every affine target line \(\mathcal L\) satisfies
\begin{equation}\label{eq:sampled256}
 \Delta(F,\mathcal L)\ge U(G)
 -2\pi mR\left(\alpha_k+
       \frac{4\beta_k}{2n-1}\left(\frac{95}{117}\right)^n\right).
\end{equation}
\end{lem}

\begin{proof}
\emph{Step 1. An integral formula for the sampled angles.}
Write \(\mathcal L=w_0+e^{i\psi}\R\). On the source line
\(z=e^{i\theta}(X+ih)\), the \(k\)-th derivative in \(X\) of
\(\Im(e^{-i\psi}(F(z)-w_0))\) is
\((n!/m!)\Im(e^{ik\theta-i\psi}G(z))\). Let
\(N_k(\theta,h;\psi)\) count the distinct real zeros of this last
slice, and put
\[
 \Lambda_k(\psi)=\frac12\int_0^\pi\int_\R
                   \bigl(m-N_k(\theta,h;\psi)\bigr)\,dh\,d\theta.
\]
Repeated application of Rolle's theorem and \eqref{eq:new-crofton}
give \(\Delta(F,\mathcal L)\ge\Lambda_k(\psi)\). Averaging in \(\psi\)
modulo \(\pi\) gives \(U(G)\), since for fixed \(\theta\) the
target angle is \(\psi-k\theta\).

Write \(V=G'/G\) and introduce the \(\pi\)-periodic function
\(P_k(v)=k^{-1}\sum_{j=0}^{k-1}|\sin((v+j\pi)/k)|\). We shall use
\begin{equation}\label{eq:sampling-integral}
 \Lambda_k(\psi)=\lim_{T\to\infty}
 \left\{2mT-\frac12\int_{|z|\le T}|V(z)|
 P_k\bigl(\psi-\arg G(z)+k\arg V(z)\bigr)\,dA(z)\right\}.
\end{equation}
At a zero of \(V\) the integrand is taken to be zero. The roots of
\(G\) may be omitted from the area integral.

To prove \eqref{eq:sampling-integral}, first count only intersections
in \(|z|\le T\). For \(G(z)\ne0\), the equation
\(k\theta=\psi-\arg G(z)\pmod\pi\) gives \(k\) source directions
\(\theta_j\) modulo \(\pi\). In each local choice of direction, the
map from \(z\) to its line parameters is
\[
 z\longmapsto
 \bigl(\theta_j(z),\Im(e^{-i\theta_j(z)}z)\bigr).
\]
Its absolute Jacobian is
\(|\Im(e^{i\theta_j}V(z))|/k\). Indeed, in the Jacobian the term
in the second differential proportional to \(d\theta_j\) drops out,
leaving the directional derivative of \(\theta_j\) along
\(e^{i\theta_j}\). The area formula, summed over \(j\), gives the
integral in \eqref{eq:sampling-integral}. Lines through a root of
\(G\) form a null set in the line parameters. For the model \(z^m\)
the phase is \(\psi-n\arg z\), so the corresponding integral,
including the factor \(1/2\), is \(2mT\).

We justify passage to the full intersection counts. Since \(G\) is
centered, each end of \(G^{-1}(e^{i(\psi-k\theta)}\R)\) has the
form \(z(s)=f_G(s\omega)\), where \(|\omega|=1\),
\(z(s)=s\omega+O_G(s^{-1})\), and
\(z'(s)=\omega+O_G(s^{-2})\), uniformly in \(\omega\).
For \(T\) sufficiently large this end exits the disk at a unique
\(s_T=T+O_G(T^{-1})\). Its projected height
\(H(s)=\Im(e^{-i\theta}z(s))\) satisfies
\(H(s)=as+O_G(s^{-1})\) and \(H'(s)=a+O_G(s^{-2})\), where
\(a=\Im(e^{-i\theta}\omega)\).

An end has \(a=0\) only if \(n\theta=\psi\pmod\pi\). Omit these
finitely many directions for each fixed \(\psi\). On every other end,
\(H(s)\) tends to \(+\infty\) or \(-\infty\) according to the sign
of \(a\). At almost every height the signed crossing count is thus
the same as for the half-line starting at \(H(s_T)\) in the
direction of \(a\). Each crossing in the opposite direction adds
two to the unsigned count. By the one-dimensional area formula,
their total contribution, integrated over heights, is at most
\(2\int_{s_T}^\infty|H'(s)-a|\,ds\).
The corresponding model half-line starts at \(aT\). Thus the
integrated discrepancy between the two omitted ends is at most
\[
 |H(s_T)-aT|+2\int_{s_T}^\infty|H'(s)-a|\,ds=O_G(T^{-1}).
\]
The implied constant is independent of \(\theta\) and \(\psi\).
Summing over the \(2m\) ends and integrating in \(\theta\) leaves
an error \(O_G(T^{-1})\), uniformly in \(\psi\). This proves
\eqref{eq:sampling-integral} and permits averaging it in \(\psi\).

\emph{Step 2. The contribution from \(|z|\le2R\).}
The Fourier series of \(|\sin v|\), averaged over the \(k\)
directions, gives
\begin{equation}\label{eq:angular-sampling-series}
 P_k(v)=\frac2\pi-\frac4\pi\sum_{q\ge1}
                         \frac{\cos(2qv)}{4k^2q^2-1}.
\end{equation}
On \(0\le v\le\pi\), summing the sines directly gives
\(P_k(v)=\cos((v-\pi/2)/k)/(k\sin(\pi/(2k)))\).
Hence \(\max P_k=2/\pi+\alpha_k\) and
\(\min P_k=P_k(0)=2/\pi-\beta_k\). In particular, the sum of
the absolute values of the cosine coefficients in
\eqref{eq:angular-sampling-series} is \(\beta_k\).

Subtracting \eqref{eq:sampling-integral} from its angular mean
expresses \(U(G)-\Lambda_k(\psi)\) as the limit of the integral
of \(\tfrac12|V|(P_k-2/\pi)\) over \(|z|\le T\).
On \(|z|\le2R\), use \(|V(z)|\le\sum_j|z-\xi_j|^{-1}\).
The integral of each kernel over this disk is at most \(4\pi R\),
as is seen by centering the disk at the pole and applying circular
rearrangement. The contribution from \(|z|\le2R\) is therefore
at most \(2\pi mR\alpha_k\).

\emph{Step 3. The exterior contribution.}
If \(R=0\), then \(G=z^m\) and \(U(G)=0\), so the assertion
follows from \(\Delta(F,\mathcal L)\ge0\). Otherwise scale the source
so that \(R=1\). Put
\[
 B(t)=\prod_{j=1}^m(1-\xi_jt),\qquad
 A(t)=\frac1m\sum_{j=1}^m\frac1{1-\xi_jt}.
\]
Then \(G(z)=z^m B(1/z)\) and \(V(z)=mA(1/z)/z\). On
\(|t|\le p<1\) we have
\begin{equation}\label{eq:derivative-disk}
 (1-p)^m\le|B(t)|\le(1+p)^m,\qquad
 \frac1{1+p}\le|A(t)|\le\frac1{1-p}.
\end{equation}
For the last two bounds, observe that \(v\mapsto(1-v)^{-1}\)
maps \(|v|\le p\) onto the closed disk with center \(1/(1-p^2)\)
and radius \(p/(1-p^2)\). The average defining \(A\) lies in this
disk as well. We take the analytic square root with \(A(0)^{1/2}=1\).

For \(q\ge1\), let \(C_q(r)\) be the circular mean on \(|z|=r\)
of \(|V(z)|\exp(2iq(k\arg V(z)-\arg G(z)))\). Define
\(H_+=A^{kq+1/2}B^{-q}\), \(H_-=A^{-kq+1/2}B^q\), and
\(H_-^\#(t)=\overline{H_-(\bar t)}\). Writing \(z=r\omega\) gives
\[
 C_q(r)=\frac mr\,\frac1{2\pi i}
 \int_{|\omega|=1}\omega^{-2nq}
       H_+\!\left(\frac1{r\omega}\right)
       H_-^\#\!\left(\frac{\omega}{r}\right)
       \frac{d\omega}{\omega}.
\]
The integrand is analytic in \(r^{-1}<|\omega|<r\). For \(r\ge2\),
move the contour to \(|\omega|=cr\), with \(c=9/10\). The arguments
of \(H_+\) and \(H_-^\#\) now have moduli \(1/(cr^2)\) and \(c\),
respectively. By \eqref{eq:derivative-disk}, the two square-root
factors contribute at most
\(\bigl((1-c)(1-1/(4c))\bigr)^{-1/2}<4\). The remaining powers
give
\begin{equation}\label{eq:angular-exterior-mode}
 |C_q(r)|\le\frac{4m}{r}
 \left(\frac{1+c}{c^2r^2-c}\right)^{nq}
 \le\frac{4m}{r}
 \left(\frac{95}{117}\left(\frac2r\right)^2\right)^{nq}.
\end{equation}
Integrating \eqref{eq:angular-exterior-mode} in \(r\), with the
area factor and the factor \(1/2\) from \eqref{eq:sampling-integral},
bounds the contribution of the \(q\)-th mode before its Fourier
coefficient by
\(8\pi m(95/117)^{nq}/(2nq-1)\). Summing against the absolute
cosine coefficients in \eqref{eq:angular-sampling-series} gives
at most \(8\pi m\beta_k(95/117)^n/(2n-1)\). Restoring \(R\),
adding the bound from Step~2, and using
\(\Delta(F,\mathcal L)\ge\Lambda_k(\psi)\) proves \eqref{eq:sampled256}.
\end{proof}

\subsection{A product estimate for inverse beta moments}

The reconstruction of a polynomial from a higher derivative will
require a bound for a product involving a single beta-distributed
random variable. The binomial coefficients in its negative moments
allow us to estimate the entire product at once.

\begin{lem}\label{lem:inverse-beta-moments}
Let \(n\ge256\) and \(2\le k\le33n/256\) be integers, and put
\(m=n-k\). Let \(\lambda\) have density
\(\lambda^m(1-\lambda)^{k-2}/B(m+1,k-1)\) on \((0,1)\), and set
\[
 Z=\frac{m}{(n-1)\lambda}-1.
\]
If \(w_1,\ldots,w_m\in\C\) satisfy \(\sum_i|w_i|\le46/5\), then
\[
 \left|\mathbb E\prod_{i=1}^m(1-w_iZ)-1\right|<\frac1{25}.
\]
\end{lem}

\begin{proof}
\emph{Step 1. Positivity of the centered moments.}
Put \(N=n-1\) and \(p=(k-1)/N\), so \(m=N(1-p)\). The beta
integral gives
\begin{equation}\label{eq:angular-beta-raw-moments}
 \mathbb E\lambda^{-j}
 =\frac{B(m+1-j,k-1)}{B(m+1,k-1)}
 =\frac{\binom Nj}{\binom mj},\qquad 0\le j\le m.
\end{equation}
In particular, for \(M_j=\mathbb E Z^j\) we have \(M_0=1\) and
\(M_1=0\). We claim that
\begin{equation}\label{eq:centered-beta-recurrence}
 (m-j+1)M_j
 =(j-1)\bigl(pM_{j-2}+(1+p)M_{j-1}\bigr),
 \qquad 2\le j\le m.
\end{equation}
Indeed, \(Y=\lambda^{-1}-1\) has density
\(\rho(y)\) proportional to \(y^{k-2}(1+y)^{-n}\) on \(y>0\).
Differentiation gives
\((y(1+y)\rho(y))'=(k-1-my)\rho(y)\).
Integrating by parts and using \(Z=(m/N)(1+Y)-1\), we obtain
\[
 m\,\mathbb E[Zg(Z)]
 =\mathbb E[(p+Z)(1+Z)g'(Z)]
\]
for every polynomial \(g\) of degree at most \(m-1\). The boundary
term vanishes at zero since \(k\ge2\), and at infinity it is
\(O(y^{-1})\) even in degree \(m-1\). Taking \(g(Z)=Z^{j-1}\)
gives \eqref{eq:centered-beta-recurrence}. Induction now shows that
\(M_j\ge0\) for \(0\le j\le m\).

\emph{Step 2. Summing the product.}
Let \(W=46/5\) and \(s=W/m<1\). For the \(j\)-th elementary
symmetric polynomial, Maclaurin's inequality gives
\[
 |e_j(w_1,\ldots,w_m)|
 \le e_j(|w_1|,\ldots,|w_m|)
 \le\binom mj s^j.
\]
Indeed, with the sum of two nonnegative entries fixed, averaging
them can only increase \(e_j\), since the coefficient of their
product is nonnegative. At a maximum of least variance all entries
must therefore agree. Expanding the product, using \(M_1=0\)
and \(M_j\ge0\), gives
\begin{equation}\label{eq:angular-beta-product}
 \left|\mathbb E\prod_{i=1}^m(1-w_iZ)-1\right|
 \le \sum_{j=2}^m\binom mj s^jM_j
 =\mathbb E(1+sZ)^m-1.
\end{equation}

Now \(1+sZ=1-s+sm/(N\lambda)\). Expanding in \(\lambda^{-1}\)
and applying \eqref{eq:angular-beta-raw-moments} cancels the
factor \(\binom mj\) exactly. Since the summands are nonnegative,
\[
 \mathbb E(1+sZ)^m
 =(1-s)^m\sum_{j=0}^m\binom Nj
       \left(\frac{sm}{N(1-s)}\right)^j
 \le \frac{(1-ps)^N}{(1-s)^{Np}}.
\]
Taking logarithms and expanding at \(s=0\), we find
\[
 \log\mathbb E(1+sZ)^m
 \le N\sum_{j\ge2}\frac{p-p^j}{j}s^j
 \le \frac{Np\,s^2}{2(1-s)}
 =\frac{(k-1)W^2}{2m(m-W)}.
\]
The degree assumptions give \(m\ge223\) and
\((k-1)/m\le33/223\). Thus the last expression is at most
\(33W^2/(2\cdot223(223-W))<1/32\). By
\eqref{eq:angular-beta-product}, the required error is less than
\(e^{1/32}-1<1/31<1/25\).
\end{proof}

\section{An effective lower bound for the inverse-line deficit}\label{sec:effective-line}

\begin{thm}\label{thm:line256}
Let \(F\) be a monic polynomial of degree \(n\ge256\), and suppose
that its finite critical values have diameter \(d>0\). Then, for
every affine line \(\mathcal L\subset\C\),
\[
 \Delta(F,\mathcal L)>\frac45\left(\frac d2\right)^{1/n}.
\]
\end{thm}

The line \(\mathcal L\) is arbitrary. The proof uses the smallness of
\(\Delta(F,\mathcal L)\) twice. First, it bounds the sum of the moduli
of the roots of a derivative of order about \(\sqrt n\). This
controls the terms of \(F\) that survive differentiation. To
recover the remaining terms, we find a family of parallel lines
on which the real slices of \(F\) have many zeros in a short
interval. Two lines from this family suffice to bound the derivative
of the part that was lost.

\begin{proof}
\emph{Step 1. The roots of a higher derivative.}
Replacing \(F(z)\) by
\(\sigma^{-n}(F(\sigma z+z_*)-w_0)\) divides source lengths by
\(|\sigma|\) and critical-value distances by \(|\sigma|^n\).
Choose \(w_0\in\mathcal L\), take \(\sigma^n\) parallel to \(\mathcal L\)
with \(|\sigma|^n=d/2\), and choose \(z_*\) to center the
polynomial. We may thus assume that \(d=2\), that \(F\) is
centered, and that \(\mathcal L=\R\). Suppose, for a contradiction,
that \(\Delta(F,\R)\le L_0=4/5\).

Put \(k=\lceil2\sqrt n\rceil\), \(m=n-k\), and \(u=33/256\).
Then \(k/n\le u\), \(m\ge223\), and \(k^2\ge4n\). Let
\(G=m!F^{(k)}/n!\), and write \(\xi_1,\ldots,\xi_m\) for its
roots, counted with multiplicity. We claim that
\begin{equation}\label{eq:rootcaps256}
 R:=\max_i|\xi_i|<\frac{577}{1000},
 \qquad
 S:=\sum_{i=1}^m|\xi_i|<\frac{229}{100}.
\end{equation}

To see this, put \(v=\pi/(2k)\). The inequalities
\(\sin v\ge v-v^3/6\), \(\sin v\le v\), and
\(\cos v\ge1-v^2/2\) give
\[
 k^2\alpha_k\le\frac{\pi/12}{1-\pi^2/(24k^2)},
 \qquad k^2\beta_k\le\frac{\pi}{4}.
\]
Since \(k^2\ge4n\) and \(n\ge256\), the error term in
\eqref{eq:sampled256} is at most \(e_0R\), where
\[
 e_0=\frac{\pi}{2}
 \left(\frac{\pi/12}{1-\pi^2/(24\cdot1024)}
                +\frac{\pi}{511}\right)<\frac{423}{1000}.
\]
The polynomial \(G\) is centered, so
\(|\xi_i|=m^{-1}|\sum_j(\xi_i-\xi_j)|\le D_G\) for every \(i\).
Moreover, \(t=(m-2)/2\ge221/2\), and
\(2t/(1+t+\sqrt t)>181/100\) on this range. Applying
\eqref{eq:root-mean256} and \eqref{eq:sampled256}, we obtain
\[
 D_G\le\frac{L_0}{181/100-e_0},
 \qquad
 S\le\frac{2(L_0+e_0D_G)}{c_*}.
\]
Using \(R\le D_G\) and \(c_*>913/1000\) proves
\eqref{eq:rootcaps256}.

\emph{Step 2. The terms determined by \(F^{(k)}\).}
Write \(G(z)=\sum_{j=0}^m g_jz^{m-j}\), with \(g_0=1\), and
decompose \(F=T+Q\), where \(\deg Q<k\) and
\begin{equation}\label{eq:top256}
 T(z)=\sum_{j=0}^m\frac{\binom nj}{\binom mj}\,g_jz^{n-j}.
\end{equation}
Indeed, the \(k\)-th derivative of the right-hand side is
\((n!/m!)G\). Since \((-1)^jg_j\) is the \(j\)-th elementary
symmetric polynomial in the \(\xi_i\), Maclaurin's inequality
applied to their moduli gives
\(\lvert g_j\rvert\le\binom mj(S/m)^j\). Thus the denominator
in \eqref{eq:top256} cancels, and for every \(\rho>0\) we have
\begin{equation}\label{eq:top-growth-all-radii256}
 \|T\|_\rho
 \le \sum_{j=0}^m\binom nj(S/m)^j\rho^{n-j}
 \le (\rho+S/m)^n
 \le \rho^n\exp\!\left(\frac{nS}{m\rho}\right).
\end{equation}
Take \(\rho_0=23/48\) and \(r=19/20\). Since
\((229/100)/((1-u)\rho_0)<6\) and \(e^6<500\),
\eqref{eq:rootcaps256} and \eqref{eq:top-growth-all-radii256} give
\begin{equation}\label{eq:topnorm256}
 \|T\|_\rho\le500\rho^n,\qquad \rho_0\le\rho\le r.
\end{equation}

We also need a lower bound for \(T'\) on \(|z|=r\). Let
\(\lambda\) have the beta distribution with parameters \(m+1\)
and \(k-1\). Its negative moments are
\(\mathbb E\lambda^{-j}=\binom{n-1}{j}/\binom mj\) for
\(0\le j\le m\). Differentiating \eqref{eq:top256} and expanding
the product over the roots of \(G\) gives
\begin{equation}\label{eq:inverse-beta256}
 \frac{T'(z)}{nz^{n-1}}
 =\sum_{j=0}^m\frac{\binom{n-1}{j}}{\binom mj}\,g_jz^{-j}
 =\mathbb E\prod_{i=1}^m\left(1-\frac{\xi_i}{\lambda z}\right),
 \qquad z\ne0.
\end{equation}
Put \(a=(n-1)/m\), \(Z=\lambda^{-1}/a-1\), and
\(w_i=a\xi_i/(z-a\xi_i)\). For \(|z|=r\),
\eqref{eq:rootcaps256} gives \(aR<r\) and
\[
 \sum_i|w_i|
 \le\frac{S}{r/a-R}
 \le\frac{S}{r(1-u)-R}<W,\qquad W=\frac{46}{5}.
\]
Factoring \(1-a\xi_i/z\) from the \(i\)-th factor in
\eqref{eq:inverse-beta256} leaves \(1-w_iZ\). Hence, for the
comparison polynomial
\(B(z)=nz^{k-1}\prod_{i=1}^m(z-a\xi_i)\),
Lemma~\ref{lem:inverse-beta-moments} gives
\(\lvert T'(z)/B(z)-1\rvert<1/25\) on \(|z|=r\).

All zeros of \(B\) lie in \(|z|<r\), so Rouch\'e's theorem places
all zeros of \(T'\) there as well. Moreover,
\(\log(1-y)\ge-y/(1-y)\) for \(0\le y<1\), and therefore
\[
 \prod_i|1-a\xi_i/z|
 \ge \exp\!\left(-\sum_i\frac{a|\xi_i|}{r-a|\xi_i|}\right)
 \ge e^{-W}.
\]
We have proved
\begin{equation}\label{eq:toplower256}
 \min_{|z|=r}|T'(z)|
 \ge \frac{24}{25}\,nr^{n-1}e^{-W}.
\end{equation}

\emph{Step 3. Real zeros on parallel lines.}
We give the count for arbitrary \(s,H>0\) and \(0<\gamma<1\),
since we shall use it again with a different choice of scales.
Parametrize the source lines by \(z=e^{i\theta}(X+ih)\), with
\(0\le\theta<2\pi\) and \(0<h<H\). Let \(N\) be the number of
distinct real zeros of the slice
\(X\mapsto\Im F(e^{i\theta}(X+ih))\), and let \(J\) count the
roots of \(G\), with multiplicity, for which
\(\Im(e^{-i\theta}\xi_i)>h\). We omit the null set of lines
through a root of \(G\) and the finitely many directions with
\(\sin(n\theta)=0\).

A root at the origin contributes \(2\arctan(s/h)\) to the
decrease of the argument along \(-s\le X\le s\). For a root
\(\xi_i\) below the line, both endpoint arguments can be taken
in the upper half-plane after rotating the line to the horizontal.
If \(|\xi_i|<\sqrt{s^2+h^2}\), moving the root from zero to
\(\xi_i\) changes each endpoint argument by at most
\(\arcsin(|\xi_i|/\sqrt{s^2+h^2})\). The loss in the argument
decrease is therefore at most \(\pi|\xi_i|/s\), since
\(\arcsin y\le\pi y/2\) on \([0,1]\). For a larger root below
the line the same bound follows because both argument decreases
belong to \([0,\pi]\). A root above the line costs at most \(2\pi\).
Thus the total decrease of \(\arg G\), divided by \(\pi\), is
at least \(mA-S/s-2J\), where \(A=(2/\pi)\arctan(s/H)\).

Put \(H_{\theta,h}(X)=\Im F(e^{i\theta}(X+ih))\) for real \(X\).
Since \(G=m!F^{(k)}/n!\), we have
\((m!/n!)H_{\theta,h}^{(k)}(X)
=\Im(e^{ik\theta}G(e^{i\theta}(X+ih)))\).
Take a continuous argument of \(e^{ik\theta}G(e^{i\theta}(X+ih))\)
on \([-s,s]\). Each multiple of \(\pi\) strictly between its endpoint
values gives a distinct zero of \(H_{\theta,h}^{(k)}\) in \((-s,s)\).
The factor \(e^{ik\theta}\) does not affect the change in argument,
so there are at least \(mA-S/s-2J-1\) such zeros.
On each of the two complementary open rays, Rolle's
theorem bounds the number of zeros of the original slice by
the number of zeros of its \(k\)-th derivative plus \(k\).
Since that derivative has degree \(m\), the number of zeros of
the original slice in \([-s,s]\) is at least
\[
 mA-k-\frac{S}{s}-1-\bigl(2J+n-N\bigr).
\]

It remains to choose one direction for which this estimate holds
with a useful bound on \(2J+n-N\) at many heights. Its mean,
with respect to \(d\theta\,dh/(2\pi H)\), is at most
\((2S+\Delta(F,\R))/(\pi H)\). Indeed, each root contributes at
most \(2|\xi_i|\) to the integral of \(J\), since
\(\int_0^{2\pi}[\Im(e^{-i\theta}\xi_i)]_+\,d\theta=2|\xi_i|\).
For \(n-N\), use \eqref{eq:new-crofton}, observing that a full
turn with positive heights counts the same unoriented lines as
a half turn with heights of both signs.

Markov's inequality and Fubini's theorem now give a direction
\(\theta\) and a measurable set \(E\subset(0,H)\) with
\(|E|\ge\gamma H\) such that, for every \(h\in E\), the number
\(M_h\) of distinct real zeros in \([-s,s]\) satisfies
\begin{equation}\label{eq:parallel-count256}
 M_h\ge mA-k-\frac{S}{s}-1
       -\frac{2S+\Delta(F,\R)}{\pi H(1-\gamma)}.
\end{equation}

For the present theorem take \(s=1/6\) and \(H=\gamma=1/16\).
The inequalities \(\arctan(3/8)<3/8\) and \(\pi>25/8\) give
\(A=1-(2/\pi)\arctan(3/8)>19/25\). By \eqref{eq:rootcaps256},
we also have \(S/s<55/4\) and
\((2S+L_0)/(\pi H(1-\gamma))<59/2\). Consequently,
\eqref{eq:parallel-count256} gives
\[
 \frac{M_h}{n}>
 \frac{223}{256}\frac{19}{25}-\frac{33}{256}
 -\frac{55/4+1+59/2}{256}
 =\frac{9223}{25600}>\frac9{25},
 \qquad h\in E.
\]

\emph{Step 4. Recovery of the lower coefficients.}
Fix the direction obtained in Step~3. For a polynomial \(C\),
write \(C_\theta(X)=C(e^{i\theta}X)\) and
\(C_\theta^\#(X)=\overline{C_\theta(\bar X)}\). The polynomial
that agrees with \(\Im C(e^{i\theta}(X+ih))\) on the real axis is
\begin{equation}\label{eq:complex-slice256}
 \mathcal S_h C(X)
 =\frac{C_\theta(X+ih)-C_\theta^\#(X-ih)}{2i}.
\end{equation}
Put \(P_h=\mathcal S_hQ\). For \(b>s\), we have
\(\|\mathcal S_hT\|_b\le\|T\|_{b+H}\) whenever \(0<h<H\).

If a polynomial \(f\) has at least \(M\) zeros in \(|X|\le s\),
then \(\|f\|_s\le q^M\|f\|_b\), where
\(q=2sb/(s^2+b^2)\). Indeed, divide \(f\) by the Blaschke factors
of its zeros in the disk of radius \(b\). Each factor has modulus
one on \(|X|=b\) and modulus at most \(q\) on \(|X|\le s\) when
its zero lies in \(|X|\le s\). Apply this to
\(\mathcal S_hF=P_h+\mathcal S_hT\). Since \(\deg P_h<k\),
the maximum principle applied to the reciprocal polynomial gives
\(\|P_h\|_b\le(b/s)^k\|P_h\|_s\). We obtain
\begin{equation}\label{eq:slice-absorption256}
 \|P_h\|_s
 \le\frac{1+q^M}{1-q^M(b/s)^k}\,\|T\|_{b+H}
 \qquad\text{if }M\le M_h\text{ and }q^M(b/s)^k<1.
\end{equation}

Here \(b=5/12\), so \(b+H=\rho_0\), \(b/s=5/2\), and \(q=20/29\).
For \(M=9n/25\) and \(k\le un\), we have
\(q^M(b/s)^k<1/20\). At \(n=256\) this follows from the integer
comparison \(20^{25}(20/29)^{2304}(5/2)^{825}<1\), and the
bound decreases with \(n\). The fraction in
\eqref{eq:slice-absorption256} is therefore at most \(21/19<2\).
By \eqref{eq:topnorm256},
\begin{equation}\label{eq:slice-small256}
 \|P_h\|_s\le1000\rho_0^n,\qquad h\in E.
\end{equation}

We next recover \(Q'\) from two of these slices. We give the
estimate in terms of \(M_*=\sup_{h\in E}\|P_h\|_s\).
Provided \(16kH/s\ge1\), partition \((0,H)\) into
\(\lceil16kH/s\rceil\) equal intervals. Their lengths lie between
\(s/(32k)\) and \(s/(16k)\), and one contains a proportion at
least \(\gamma\) of \(E\) in the sense that \(E\) occupies at
least \(\gamma\) times its length. We can thus choose
\(h_0<h_1\) in this interval, both in \(E\), with
\(\gamma s/(64k)\le\eta:=h_1-h_0\le s/(16k)\).

Put \(A_0(X)=Q_\theta(X+ih_0)\). Subtracting the two reflected
expressions in \eqref{eq:complex-slice256} gives
\[
 D_\eta A_0(X)
 :=\frac{A_0(X+i\eta)-A_0(X-i\eta)}{2i\eta}
 =\frac{P_{h_1}(X)-P_{h_0}(X-i\eta)}{\eta}.
\]
Hence \(\|D_\eta A_0\|_{3s/4}\le2M_*/\eta\). To pass from the
difference to the derivative, use the polynomial identity
\[
 A_0'=\frac{\eta\partial_X}{\sin(\eta\partial_X)}D_\eta A_0.
\]
It follows by Taylor expansion of the two translations. All
operator series terminate on polynomials. The Taylor coefficients
of \(w/\sin w\) have modulus at most two, by Cauchy's estimate
on \(|w|=1\) and
\(\lvert1-\sin w/w\rvert\le\sum_{j\ge1}1/(2j+1)!<1/2\).
Cauchy's derivative estimate from \(|X|\le3s/4\) to \(|X|\le s/2\)
costs \(j!(4/s)^j\) at order \(j\). Since \(j\le k\) and
\(4k\eta/s\le1/4\), the total operator norm is less than
\(2\sum_{j\ge0}4^{-j}<3\). Thus \(\|A_0'\|_{s/2}\le6M_*/\eta\).

The last disk corresponds in the original coordinate to a disk
centered at \(ie^{i\theta}h_0\). Polynomial growth from that disk
gives, for every \(\varrho>0\),
\begin{equation}\label{eq:two-height-growth256}
 \|Q'\|_\varrho
 \le\frac{384kM_*}{\gamma s}
      \max\!\left\{1,\frac{2(\varrho+H)}s\right\}^{k}.
\end{equation}
In particular, \eqref{eq:slice-small256} yields
\begin{equation}\label{eq:lowpoly256}
 \|Q'\|_r\le Pk\rho_0^nB_0^k,\qquad
 P=\frac{384000}{\gamma s},\qquad B_0=\frac{2(r+H)}s.
\end{equation}

To finish, note that \((\rho_0/r)B_0^u<3/4\). Equivalently,
\((115/171)^{256}(243/20)^{33}<1\).
Also \(4e^WPur(3/4)^{256}<1\), as follows already from
\(e^W<3^{10}\) and \(P=36\,864\,000\). Since \(k\le un\),
\eqref{eq:lowpoly256} now gives
\(\|Q'\|_r<nr^{n-1}e^{-W}/4\). By \eqref{eq:toplower256},
Rouch\'e's theorem places all zeros of \(F'=T'+Q'\) in
\(|z|<r\). On this disk,
\[
 |F(z)-Q(0)|
 \le \|T\|_r+r\|Q'\|_r
 <(500+n/4)r^n<\tfrac14.
\]
For the last inequality, both \(r^n\) and \(nr^n\) decrease for
\(n\ge256\), and \(r^{16}<1/2\). All critical values of \(F\)
therefore lie in a disk of radius \(1/4\), contradicting their
diameter \(2\). This contradiction proves the theorem.
\end{proof}

The numerical comparisons in Theorem~\ref{thm:line256} can also
be checked by rational arithmetic with
\href{https://github.com/zhangteng2000/ehp-lemniscate-verification-20261108/blob/ca8c451bfbc7312b9c57acc5e672c0b69cd601ae/verification/effective_line_constants.py}{\nolinkurl{effective_line_constants.py}} and
\href{https://github.com/zhangteng2000/ehp-lemniscate-verification-20261108/blob/ca8c451bfbc7312b9c57acc5e672c0b69cd601ae/verification/audit_effective_line_256.py}{\nolinkurl{audit_effective_line_256.py}}. Both checks use only the
weaker factor \(1/2\) in \eqref{eq:toplower256}. In the proof of
Lemma~\ref{lem:inverse-beta-moments}, we obtain
\(\log\mathbb E[(1+sZ)^m]\le34914/1191935<1/32\).
An alternative estimate obtained by splitting the moment sum, together
with checks of \eqref{eq:centered-beta-recurrence} against the raw
moments, is given in
\href{https://github.com/zhangteng2000/ehp-lemniscate-verification-20261108/blob/ca8c451bfbc7312b9c57acc5e672c0b69cd601ae/verification/inverse_beta_moment_audit.py}{\nolinkurl{inverse_beta_moment_audit.py}}.

\subsection{The limiting deficit in terms of the critical points}

The same argument also gives an asymptotically sharp estimate
in terms of the critical points. Put
\(D(F)=\operatorname{diam}\{\zeta\in\C:F'(\zeta)=0\}\).

\begin{thm}\label{thm:sharp-line-limit}
For each \(n\ge3\), take the infimum over all monic polynomials
\(F\) of degree \(n\) with \(D(F)>0\) and all affine lines
\(\mathcal L\subset\C\). Then
\[
 \lim_{n\to\infty}
 \inf_{\substack{\deg F=n,\ D(F)>0\\ \mathcal L\ \mathrm{an\ affine\ line}}}
 \frac{\Delta(F,\mathcal L)}{D(F)}=2.
\]
\end{thm}

\begin{proof}
\emph{Step 1. The coefficients lost under differentiation.}
For the lower bound, consider any sequence with degrees tending
to infinity and with bounded ratio \(\Delta(F,\mathcal L)/D(F)\).
The deficit and \(D(F)\) scale by the same factor under a source
similarity. We may therefore normalize \(D(F)=1\), center \(F\),
and take \(\mathcal L=\R\). In particular, \(\Delta(F,\R)\) is bounded
along the sequence.

Set \(k=\lceil n^{2/3}\rceil\) and \(m=n-k\), and use the
decomposition \(F=T+Q\) in \eqref{eq:top256} for this choice of
\(k\). The critical points of \(F\) have barycenter zero and
diameter one, so they belong to the closed unit disk. By the
Gauss--Lucas theorem, the roots of \(G=m!F^{(k)}/n!\) belong to
their convex hull. In particular, \(R\le1\) and \(D_G\le1\).
Since \(\alpha_k,\beta_k=O(k^{-2})\) and \(n/k^2\to0\),
\eqref{eq:sampled256} gives
\begin{equation}\label{eq:sharp-line-mean}
 \Delta(F,\R)\ge U(G)-o(1).
\end{equation}
The first inequality in \eqref{eq:root-mean256} now gives a
uniform bound \(S\le S_*\) along the sequence.

We claim that for every fixed \(r>0\),
\begin{equation}\label{eq:sharp-line-remainder}
 \lim_{n\to\infty}\|Q'\|_r^{1/n}=0.
\end{equation}
Fix \(\rho>0\), and take \(s=\rho/4\), \(H=\rho/8\), \(b=2s\),
and \(\gamma=1/16\). Since \(n/m\to1\), the general bound
\eqref{eq:top-growth-all-radii256} gives
\(\|T\|_{b+H}\le\|T\|_\rho=O_\rho(\rho^n)\).
In \eqref{eq:parallel-count256}, the terms involving \(S\) and
\(\Delta(F,\R)\) are \(O_\rho(1)\), and \(k=o(n)\).
It follows that on the selected heights \(M_h\ge cn\) for all
sufficiently large \(n\), where one may take
\(c=\pi^{-1}\arctan 2>0\). Since \(q=4/5\) and
\((4/5)^{cn}2^k\to0\), \eqref{eq:slice-absorption256} gives
\(\sup_{h\in E}\|P_h\|_s=O_\rho(\rho^n)\).
By \eqref{eq:two-height-growth256},
\(\|Q'\|_r\le C_\rho k\rho^n C_{\rho,r}^{\,k}\), with constants
independent of \(n\). Taking \(n\)-th roots gives a limsup at
most \(\rho\). Letting \(\rho\downarrow0\) proves
\eqref{eq:sharp-line-remainder}.

\emph{Step 2. Comparison of the critical points with the roots of \(G\).}
For every \(0<\varepsilon<1\), we have
\begin{equation}\label{eq:sharp-line-derivative-comparison}
 \sup_{\varepsilon\le|z|\le1}
 \left|\frac{F'(z)}{nz^{n-1}}-\frac{G(z)}{z^m}\right|
 \longrightarrow0.
\end{equation}
Indeed, the coefficient of \(z^{-j}\) in \(T'(z)/(nz^{n-1})\)
is \(\binom{n-1}{j}g_j/\binom mj\), by
\eqref{eq:inverse-beta256}. For each fixed \(j\), this binomial
ratio tends to one. Moreover, the coefficient bound used in
\eqref{eq:top-growth-all-radii256} gives, for all sufficiently
large \(n\),
\[
 \frac{\binom{n-1}{j}}{\binom mj}|g_j|
 \le\binom{n-1}{j}(S_*/m)^j
 \le\frac{(2S_*)^j}{j!}.
\]
The coefficients of \(G(z)/z^m\) satisfy
\(|g_j|\le S_*^j/j!\). These bounds are summable after
multiplication by \(\varepsilon^{-j}\), giving uniform convergence
of the difference with \(T'\) in place of \(F'\). The contribution
of \(Q'\) tends uniformly to zero by \eqref{eq:sharp-line-remainder}.

To see what \eqref{eq:sharp-line-derivative-comparison} implies for
the roots, suppose \(z\ne0\) is at distance at least
\(\varepsilon\) from every \(\xi_i\). Then
\[
 \left|\frac{z^m}{G(z)}\right|
 =\prod_{i=1}^m\left|1+\frac{\xi_i}{z-\xi_i}\right|
 \le\exp(S_*/\varepsilon).
\]
Thus on \(\varepsilon\le|z|\le1\) a zero of \(F'\) must, for
all sufficiently large \(n\), be within \(\varepsilon\) of a
root of \(G\). For \(|z|<\varepsilon\), the distance to the
nearest root of \(G\) is at most \(\varepsilon+S_*/m\), since
\(\min_i|\xi_i|\le S_*/m\). As \(\varepsilon>0\) was arbitrary,
every critical point of \(F\) is within \(o(1)\) of a root of
\(G\), uniformly along the sequence.

Choose two critical points of \(F\) at distance one. Their nearby
roots of \(G\) are at distance at least \(1-o(1)\). Since
\(D_G\le1\), we obtain \(D_G\to1\). With \(t=(m-2)/2\), the
factor \(2t/(1+t+\sqrt t)\) tends to two. The second inequality
in \eqref{eq:root-mean256} and \eqref{eq:sharp-line-mean}
therefore give \(\liminf\Delta(F,\R)\ge2\).

\emph{Step 3. Sharpness.}
Take \(F(z)=z^{n-1}(z-a_n)\), where \(a_n=n/(n-1)\).
Since \(F'(z)=nz^{n-2}(z-1)\), we have \(D(F)=1\). The roots of
\(F\) are zero with multiplicity \(n-1\) and \(a_n\) with
multiplicity one. By \eqref{eq:angular-root-energy} and scaling,
\(U(F)=a_n\mathcal C(n-1)\).

We claim that \(\mathcal C(t)\le2\) for every \(t>0\). Projecting
\(t/z+1/(z-1)\) onto the direction of \(1/z\) bounds the integrand
in \(\pi\mathcal C(t)\) from above by
\[
 \frac1{|z-1|}-\Re\frac{z}{|z|(z-1)}.
\]
For \(T>1\), the integral of the second term over \(|z|<T\) is
\(2\pi(T-1)\). Indeed, in polar coordinates the circular mean of
\(z/(z-1)\) is zero for \(|z|<1\) and one for \(|z|>1\).
The integral of the first term is at most \(2\pi T\), by circular
rearrangement. Integrating the nonnegative triangle defect over
\(|z|<T\) and letting \(T\to\infty\) proves the claim.

The infimum over target lines is at most the mean over lines
through the origin. For this choice of \(F\), it is therefore
at most \(2a_n=2n/(n-1)\). Letting \(n\to\infty\) gives the
matching upper bound.
\end{proof}

\section{From lines to circles}\label{sec:circle256}

We first specialize Lemma~\ref{lem:curvature-general} to a fixed
choice of radii. This gives an error bound in terms of the degree
and the angular parameter alone.

\begin{lem}\label{lem:curvature256}
Let \(F\) be a centered monic polynomial of degree \(n\ge6\) whose
finite critical values lie in \([-1,1]\). Put \(x=1/n\), and let
\(F_\beta\) denote the deformation \eqref{eq:new-flow} starting at
\(F\). If \(0<\tau\le1/1000\) and
\(p(z)=1+2i\tau F_{i\tau}((2i\tau)^{-x}z)\), then
\begin{equation}\label{eq:curvature256}
 \Delta(F,\R)
 -(2\tau)^{-x}
   \bigl(\ell(\partial E_1(p_0))-\ell(\partial E_1(p))\bigr)
 \le 64n^2\tau+8n^{3/2}\tau^{1-x}.
\end{equation}
\end{lem}

\begin{proof}
Take \(R=3\) and \(b_0=1/3\) in
Lemma~\ref{lem:curvature-general}, so that \(c=3/2\).
The restrictions on \(\tau\) in that lemma are satisfied.
Write \(b=(1-x)/2\). By \eqref{eq:exterior-constant-normalizations},
\[
 H^2=\frac{n}{B(b,b)}\le\frac n\pi,
 \qquad U_n=2\pi(1+H^2)\le2\pi+2n,
\]
since \(b\le1/2\) and \(B(b,b)\ge B(1/2,1/2)=\pi\).
For both beta factors in \eqref{eq:curvature-general} we use
\begin{equation}\label{eq:shortarc-beta-bound}
 B(y,1/2)\le \frac1y+2\log 2,\qquad y>0.
\end{equation}
Indeed, after subtracting \(1/y\) from the beta integral, its
remainder is at most
\(\int_0^1((1-t)^{-1/2}-1)\,dt/t=2\log 2\), as follows by
putting \(s=\sqrt{1-t}\).

\emph{Step 1. The bounded interval and the straight tails.}
Since \(b\ge5/12\), \eqref{eq:shortarc-beta-bound} gives
\(B(b,1/2)<6\), and hence \(M<10n\) in the notation of
Lemma~\ref{lem:curvature-general}. Also,
\(e=3\tau/\sqrt{1-9\tau^2}\le10\tau/3\), so \(e+e^2<4\tau\).
The term \((1+x)\tau^2 3^{2+x}/(2+x)\) is less than \(12\tau^2\).
Thus the first two terms on the right of
\eqref{eq:curvature-general} contribute at most \(41n\tau\).

The constants for the straight tails simplify to
\[
 A=\frac{H3^x}{(2-x)\sqrt{15}}<\frac H3,\qquad
 B_{\mathrm{tail}}=\frac{32}{45}\,nH^2\left(\frac3{16}\right)^x
   <\frac{n^2}{4}.
\]
Since \(b_0-\tau>1/4\), \(H\le n\), and \(n\ge6\), it follows that
\((4A+2B_{\mathrm{tail}})\tau/(b_0-\tau)<4n^2\tau\).

\emph{Step 2. The curved tails.}
Here \(s=3/\sqrt{1+9\tau^2}>5/2\). In each row defining
\(\mathcal T\), put \(v=(j+k-1)/n\) and \(l=|j-k|/n\).
The restrictions on \(j,k\) give \(v-1\ge1/n\) and
\(l\in\Z_{\ge0}\). Hence \(K(v,l)\le2n(1+l)\) and
\(s^{-(j+k)/n}\le(2/5)^l\). For fixed \(j\), there is at most
one term with \(l=0\) and at most two with each \(l\ge1\).
Similarly, \(K(q-x,q)\le5q\) for \(q\ge2\), since
\(\min(q-x-1,1)\ge5/6\). The resulting sums satisfy
\[
 \mathcal T
 \le2n\left(1+2\sum_{l\ge1}(1+l)(2/5)^l\right)
 =\frac{82}{9}n<10n,\qquad
 \mathcal A\le5\sum_{q\ge2}q(2/5)^q=\frac{32}{9}<4.
\]
Using \(s^{1+x}\le3^{7/6}<4\) and \(H^2<n/3\), we obtain
\[
 2\tau s^{1+x}(H^2\mathcal T+2H\mathcal A)
 <8\tau\left(\frac{10}{3}n^2+8n\right)<39n^2\tau.
\]

The critical values of \(F_{i\tau}\) lie in the closed unit disk,
because \(|M_{i\tau}(v)|\le1\) for \(v\in[-1,1]\).
Lemma~\ref{lem:exterior-parseval} and
Corollary~\ref{cor:universal-length-bound} therefore give
\(h\le|h_n(F_{i\tau})|\le H\).
By \eqref{eq:shortarc-beta-bound}, \(B(x/2,1/2)<4n\), so the
remaining term is at most \(8n^{3/2}\tau^{1-x}\).
Substitution in \eqref{eq:curvature-general} proves
\eqref{eq:curvature256}, since \(41n+43n^2<64n^2\).
\end{proof}

We now give a line-deficit bound sufficient for the circle
comparison in every degree \(n\ge28\).

\begin{prop}\label{prop:uniform-circle-transfer}
Fix \(n\ge28\). Suppose that \(\Delta(F,\mathcal L)>4/5\) for every
centered monic polynomial \(F\) of degree \(n\) whose finite
critical values lie on an affine line \(\mathcal L\) and have diameter
two. Then every maximizing polynomial of degree \(n\) whose
finite critical values lie on \(\mathbb T\) has a single finite
critical value.
\end{prop}

\begin{proof}
Let \(p\) be such a maximizer, centered by a source translation.
Suppose its critical values have diameter \(2\delta>0\). Put
\(x=1/n\) and \(b=(1-x)/2\). Since
\(\ell(\partial E_1(p))\ge\ell(\partial E_1(p_0))\),
Proposition~\ref{prop:new-mass} gives
\begin{equation}\label{eq:circle-gap-28}
 \sum_{k\ge1}\mathcal E_k
 +\big\langle\Delta(F_\xi,\mathcal L_\xi)\big\rangle
 \le U_n-\ell(\partial E_1(p_0))
 <2\pi-8\log2<\frac{739}{1000}.
\end{equation}
For the last numerical inequality one may use \(\pi<355/113\)
and the first four terms of the positive series
\(\log2=2\sum_{j\ge0}3^{-(2j+1)}/(2j+1)\).

\emph{Step 1. Small and large diameters.}
Set \(\tau_0=(128n^2)^{-1}\). If
\(\delta\le2\tau_0/(1+\tau_0^2)\), the normalization
\eqref{eq:short-arc-normalization} gives \(0<\tau\le\tau_0\)
and a centered monic polynomial \(F\) whose real critical values
include both endpoints of \([-1,1]\) and lie in that interval.
Thus \(\Delta(F,\R)>4/5\). On the other hand, the right-hand side
of \eqref{eq:curvature256} is at most
\[
 \frac12+8\,128^{-(1-x)}n^{-1/2+2/n}
 \le\frac12+\frac8{128^{5/6}}<\frac45.
\]
It follows that \(\ell(\partial E_1(p))<\ell(\partial E_1(p_0))\),
contrary to maximality. We may therefore assume that
\(\delta>2\tau_0/(1+\tau_0^2)>\tau_0\).

Choose critical values \(v,v'\) of \(p\) with
\(|v-v'|=2\delta\). Away from the finite pole set, their images
under \(M_\xi\) are separated by
\(2\delta/(|1-\xi v|\,|1-\xi v'|)\). Rescaling the critical-value
diameter of \(F_\xi\) to two and applying the hypothesis gives
\begin{equation}\label{eq:circle-mean-28}
 \big\langle\Delta(F_\xi,\mathcal L_\xi)\big\rangle
 \ge\frac45\,\delta^x
 \big\langle |1-\xi v|^{-x}|1-\xi v'|^{-x}\big\rangle
 \ge\frac45\,\delta^x.
\end{equation}
The last inequality follows from Jensen's inequality and
\(\langle\log|1-\xi v|\rangle=0\) for \(v\in\mathbb T\).
Put \(q=739/800\). If \(\delta\ge q^n\), then
\eqref{eq:circle-mean-28} already contradicts \eqref{eq:circle-gap-28}.

\emph{Step 2. The intermediate diameters.}
It remains to consider \(\tau_0<\delta<q^n\). Since
\(q^{28}<1/9\), we may rotate a farthest pair to
\(c\pm i\delta\), where \(c=\sqrt{1-\delta^2}\). All critical
values lie on the shorter arc between these endpoints, as in
\eqref{eq:short-arc-normalization}. This arc lies in the closed
disk of center \(c\) and radius \(\delta\). Hence, on writing
\(p(z)=c+\delta Q(\delta^{-x}z)\), we obtain a centered monic
polynomial \(Q\) whose critical values lie in the closed unit
disk. Lemma~\ref{lem:exterior-parseval} and
Corollary~\ref{cor:universal-length-bound} give \(|h_n(Q)|\le H\).
By \eqref{eq:exterior-constant-normalizations}, \(H\le\sqrt{n/\pi}\).
It follows from \eqref{eq:affine-inverse} that
\[
 |h_n(p)|\ge bc-\sqrt{n/\pi}\,\delta.
\]

It suffices to retain the \(k=1\) term of
\eqref{eq:exterior-coefficients}. Indeed, \(b\ge27/56\) and
\(\delta<1/9\) give \(bc>7/15\) and \(2\pi b^2>4/3\).
The elementary inequality \(y^2\ge y/2-1/16\), applied to
\(y=|h_n(p)|\), now gives
\begin{equation}\label{eq:circle-linear-bound-28}
 \mathcal E_1+\big\langle\Delta(F_\xi,\mathcal L_\xi)\big\rangle
 \ge g_n(\delta),\qquad
 g_n(t)=\frac45t^{1/n}
 +\frac43\left\{\frac12\left(\frac7{15}
                       -\sqrt{n/\pi}\,t\right)-\frac1{16}\right\}.
\end{equation}
The function \(g_n\) is concave, so its minimum on
\([\tau_0,q^n]\) occurs at an endpoint.

At the left endpoint, \(\sqrt{n/\pi}\,\tau_0<1/120\). Also
\(\tau_0^{1/n}>13/20\), since \((128n^2)^{-1/n}\) increases for
\(n\ge28\) and \(128\cdot28^2(13/20)^{28}<1\). Thus
\[
 g_n(\tau_0)>
 \frac45\frac{13}{20}
 +\frac43\left(\frac{11}{48}-\frac1{16}\right)
 =\frac{167}{225}>\frac{739}{1000}.
\]
At the right endpoint, \(\sqrt n\,q^n\) decreases for \(n\ge28\),
since \(29q^2<28\). The inequalities \(\sqrt{28/\pi}<3\) and
\(q^{28}<1/9\) therefore give \(\sqrt{n/\pi}\,q^n<1/3\). Hence
\[
 g_n(q^n)>
 \frac{739}{1000}
 +\frac43\left\{\frac12\left(\frac7{15}-\frac13\right)
                       -\frac1{16}\right\}
 =\frac{739}{1000}+\frac1{180}.
\]
By \eqref{eq:circle-linear-bound-28}, these two endpoint bounds
contradict \eqref{eq:circle-gap-28}. This excludes every positive
critical-value diameter.
\end{proof}

Theorem~\ref{thm:line256} supplies the hypothesis of
Proposition~\ref{prop:uniform-circle-transfer} for \(n\ge256\).
By Lemma~\ref{lem:extremizers}, we may choose a centered maximizer
with all its finite critical values on \(\mathbb T\). The proposition
and Lemma~\ref{lem:critical-compactness} show that this maximizer
has the form \(z^n+v\) with \(|v|=1\). Its lemniscate is a rotation
of \(\partial E_1(p_0)\). This proves the length comparison for
every monic polynomial of degree \(n\ge256\).
Section~\ref{sec:equality} treats equality for an arbitrary
maximizer.

\section{The first moment of the roots}\label{sec:root-moment}

We return to the functional \(\mathcal U\) of
Section~\ref{sec:angular-derivatives}. For the root counting
measure \(\mu_G\) of a monic polynomial \(G\), we have
\(U(G)=\mathcal U(\mu_G)\). Our aim is to bound the first absolute
moment of \(\mu_G\) about its barycenter. We begin with two point
masses and keep track of their relative weights. Recall that
\(c_*=\mathcal U(\delta_0+\delta_1)=B(3/4,1/2)^2/(2\pi)\).

\begin{lem}\label{lem:centered-pair}
For \(a,b>0\) and \(z_1,z_2\in\C\),
\begin{equation}\label{eq:pair-mass-sharp}
 \mathcal U(a\delta_{z_1}+b\delta_{z_2})
 \ge 2c_*\,\frac{ab}{a+b}|z_1-z_2|.
\end{equation}
\end{lem}

\begin{proof}
The assertion is immediate if \(z_1=z_2\). Otherwise, by
homogeneity in the masses and similarity in the support, it
suffices to consider masses \(t,1\) at \(0,1\), respectively.
Write \(\mathcal C(t)=\mathcal U(t\delta_0+\delta_1)\) and put
\[
 \eta=\frac{t-1}{t+1},\qquad \sigma=\sqrt{1-\eta^2},\qquad
 K_\eta=\frac{1+t}{2t}\mathcal C(t).
\]
We shall show that \(K_\eta\ge K_0\).

Use the elliptic coordinates \(u=|z|+|z-1|\) and
\(v=|z|-|z-1|\) from the proof of
Lemma~\ref{lem:weighted-root-mean}. Put \(d=u-\eta v\).
Rationalizing the triangle defect gives
\begin{equation}\label{eq:pair-elliptic-kernel}
 K_\eta=\frac2\pi\int_1^\infty\frac{du}{\sqrt{u^2-1}}
 \int_{-1}^1
 \frac{\sqrt{1-v^2}\,dv}
      {d+\sqrt{d^2-\sigma^2(1-v^2)}}.
\end{equation}
Indeed, if \(r_0=|z|\) and \(r_1=|z-1|\), multiplication of the
defect by \(r_0r_1\) gives
\(A-\sqrt{A^2-t(1-v^2)}\), where \(A=r_0+tr_1=(1+t)d/2\).
The area element for the two half-planes together is
\(2r_0r_1\,du\,dv/\sqrt{(u^2-1)(1-v^2)}\).

To compare the weights in \eqref{eq:pair-elliptic-kernel}, use
the elementary identity
\[
 \int_{-1}^1\frac{\sqrt{1-s^2}}{d-hs}\,ds
 =\frac{\pi}{d+\sqrt{d^2-h^2}},\qquad d>|h|.
\]
Here \(h=\sigma\sqrt{1-v^2}\), and \(d>|h|\) follows from
\(u>1\) and \(\eta v+\sigma\sqrt{1-v^2}\le1\).
Write \(\sqrt{1-s^2}=\int_{|s|}^1 r(1-r^2)^{-1/2}\,dr\).
After interchanging the nonnegative integrals, put
\(x=v\) and \(y=s\sqrt{1-v^2}/r\) for each \(0<r<1\).
The resulting domain in \((x,y)\) is the unit disk, and
\begin{equation}\label{eq:pair-disk-average}
 K_\eta=\frac2{\pi^2}
 \int_1^\infty\frac{du}{\sqrt{u^2-1}}
 \int_0^1\frac{r^2\,dr}{\sqrt{1-r^2}}
 \iint_{x^2+y^2<1}\frac{dx\,dy}{u-\eta x-r\sigma y}.
\end{equation}

By rotation of the disk and reflection in \(x=0\), the innermost
integral in \eqref{eq:pair-disk-average} equals
\[
 \iint_{x^2+y^2<1}\frac{u\,dx\,dy}{u^2-\rho^2x^2},
 \qquad \rho^2=\eta^2+r^2\sigma^2=r^2+(1-r^2)\eta^2.
\]
For fixed \(u,r\), this integral increases with \(\rho^2\).
It is therefore smallest at \(\eta=0\). It follows that
\(K_\eta\ge K_0=\mathcal C(1)=c_*\), proving
\eqref{eq:pair-mass-sharp}. For distinct support points, equality
holds precisely when \(a=b\).
\end{proof}

\begin{lem}\label{lem:central-absolute}
Let \(\mu\) be a nonzero positive measure of finite support, with
barycenter \(c\). Then
\begin{equation}\label{eq:central-absolute}
 \mathcal U(\mu)\ge
 \frac{\sqrt3}{2}\,c_*\int_{\C}|z-c|\,d\mu(z).
\end{equation}
\end{lem}

\begin{proof}
Translate \(c\) to zero. With the support and total mass fixed,
the nonnegative weights having barycenter zero form a compact
polytope. The functional \(\mathcal U\) is concave in the weights,
whereas the first absolute moment is linear. It therefore suffices
to prove \eqref{eq:central-absolute} at the extreme points of
this polytope. At most three weights can be positive at an extreme
point. Indeed, any four vectors \((1,\Re z_i,\Im z_i)\) are
linearly dependent over \(\R\), and a dependence allows a small
perturbation of the weights in both directions.

An atom at the origin may be discarded by superadditivity.
For two remaining masses \(a,b\), the barycenter condition gives
\(az_1=-bz_2\), and hence
\[
 \frac{2ab}{a+b}|z_1-z_2|=a|z_1|+b|z_2|.
\]
Thus \eqref{eq:pair-mass-sharp} gives the assertion in this case.

It remains to consider three positive masses \(a_i\) at nonzero
points \(z_i\). Put \(r_i=a_i|z_i|\),
\(S=r_1+r_2+r_3\), and \(t_i=1-2r_i/S\).
Since the three vectors \(a_i z_i\) sum to zero, their lengths
satisfy the triangle inequalities. Thus \(t_i\ge0\) and
\(\sum_i t_i=1\). For \(\{i,j,k\}=\{1,2,3\}\), assign to the
pair \((z_j,z_k)\) the masses \(St_i/(2|z_j|)\) and
\(St_i/(2|z_k|)\), omitting the pair if \(t_i=0\).
The assigned masses at \(z_j\) sum to
\(S(1-t_j)/(2|z_j|)=a_j\), so the three pair measures sum to
\(\mu\).

Let \(\phi_i\in[0,\pi]\) be the angle between \(z_j\) and \(z_k\).
Since
\(|z_j-z_k|/(|z_j|+|z_k|)\ge\sin(\phi_i/2)\),
superadditivity and \eqref{eq:pair-mass-sharp} give
\begin{equation}\label{eq:centered-three-angles}
 \frac{\mathcal U(\mu)}{c_*S}
 \ge\sum_{i=1}^3 t_i\sin(\phi_i/2).
\end{equation}
The cosine rule applied to \(a_iz_i=-(a_jz_j+a_kz_k)\) gives
\[
 \sin^2(\phi_i/2)
 =\frac{(r_j+r_k)^2-r_i^2}{4r_jr_k}
 =\frac{t_i}{(1-t_j)(1-t_k)}.
\]
In particular, \(\sin(\phi_i/2)>0\) when \(t_i>0\), and
\[
 \sum_{t_i>0}\frac{t_i}{\sin^2(\phi_i/2)}
 \le\sum_{i=1}^3(1-t_j)(1-t_k)
 =1+\sum_{i<j}t_it_j\le\frac43,
\]
where \(j,k\) are the two indices other than \(i\) in the middle
sum. Jensen's inequality for \(y\mapsto y^{-1/2}\), with weights
\(t_i>0\), now gives
\(\sum_i t_i\sin(\phi_i/2)\ge(4/3)^{-1/2}=\sqrt3/2\).
Together with \eqref{eq:centered-three-angles}, this proves
\eqref{eq:central-absolute}, including the collinear case.
\end{proof}

The first-moment bound controls the sum of the root moduli. To
bound their higher moments, we also need to know how much of that
sum can be concentrated among the largest roots.

\begin{lem}\label{lem:ordered-radii}
Let \(G\) be centered and monic of degree \(m\ge5\), and suppose
that \(U(G)\le2c_*A\) for some \(A>0\). Write
\(b_1\ge\cdots\ge b_m\) for the moduli of its roots, counted with
multiplicity, and put \(S=1000A/433\) and \(D=(S+A)/(m-1)\).
Then
\begin{equation}\label{eq:ordered-moments}
 b_1\le A,\qquad
 \sum_{i=1}^m b_i^q\le A^q+(S-A)D^{q-1}
 \quad\text{for every real }q\ge1.
\end{equation}
\end{lem}

\begin{proof}
Order the roots \(\xi_i\) so that \(|\xi_i|=b_i\).
By \eqref{eq:central-absolute} and \(\sqrt3>433/250\),
\(\sum_i b_i\le4A/\sqrt3\le S\).

Fix \(1\le h<m\), and let \(v\) be the mean of the unmarked roots
\(\xi_{h+1},\ldots,\xi_m\). For each \(j>h\), place mass one at
each of \(\xi_1,\ldots,\xi_h\) and mass \(m-h\) at \(\xi_j\).
The average of these \(m-h\) measures is the root counting measure
of \(G\). Split the mass at \(\xi_j\) equally among the \(h\) pairs
\((\xi_i,\xi_j)\), \(1\le i\le h\). Concavity in the averaging
step, followed by superadditivity and \eqref{eq:pair-mass-sharp}
in the splitting step, gives
\[
 U(G)\ge\frac{2c_*}{m}\sum_{i\le h<j}|\xi_i-\xi_j|
 \ge\frac{2c_*(m-h)}{m}\sum_{i=1}^h|\xi_i-v|.
\]
Put \(y_i=\xi_i-v\) for \(i\le h\). Centering gives
\(v=-m^{-1}\sum_{i\le h}y_i\). Thus
\(\xi_i=(1-1/m)y_i-m^{-1}\sum_{j\le h,\,j\ne i}y_j\).
Summing the triangle inequality and using \(U(G)\le2c_*A\), we get
\begin{equation}\label{eq:ordered-prefix}
 \sum_{i=1}^h b_i
 \le\frac{m+h-2}{m}\sum_{i=1}^h|y_i|
 \le A+\frac{2A(h-1)}{m-h}.
\end{equation}
Taking \(h=1\) gives \(b_1\le A\).

Since \(m\ge5\), we have \(0<D<A\). We claim that
\begin{equation}\label{eq:ordered-excess}
 \sum_{i=1}^m(b_i-D)_+\le A-D,
 \qquad (y)_+=\max(y,0).
\end{equation}
Let \(h\) be the number of \(b_i>D\). The claim is immediate if
\(h=0\). If \(S\le A+(h-1)D\), the total-mass bound gives
\(\sum_{i\le h}b_i\le A+(h-1)D\). Otherwise,
\(A+(h-1)D<S=(m-1)D-A\), so \(m-h>2A/D\).
In this case \eqref{eq:ordered-prefix} gives the same bound.
Subtracting \(hD\) proves \eqref{eq:ordered-excess}.

Finally, convexity of \(y^q\) on \([0,A]\) gives
\[
 y^q\le D^{q-1}y+
       \frac{A^q-AD^{q-1}}{A-D}(y-D)_+,\qquad 0\le y\le A.
\]
The right-hand side interpolates \(y^q\) linearly on \([0,D]\)
and on \([D,A]\). Applying this inequality to each \(b_i\), then
using \(\sum_i b_i\le S\) and \eqref{eq:ordered-excess}, yields
\(\sum_i b_i^q\le SD^{q-1}+A^q-AD^{q-1}\), as required.
\end{proof}
\section{Interpolation on real slices}\label{sec:finite-interpolation}

In degrees \(28\) through \(255\), a small line deficit forces the
successive derivatives of a real slice to have many zeros in a short
interval. Interpolation at these zeros bounds the coefficients lost
under differentiation. The numerical part of the proof consists of
scalar inequalities for each degree.

We use the following consequence of \eqref{eq:sampled256} throughout
both this section and Section~\ref{sec:collinear-gram}. Let \(F\) be
centered and monic of degree \(n\), let \(1\le k\le n-2\), and put
\(m=n-k\) and \(G=m!F^{(k)}/n!\). Write \(\xi_1,\ldots,\xi_m\)
for the roots of \(G\), counted with multiplicity, and put
\(R=\max_i|\xi_i|\). Then, for every affine line \(\mathcal L\),
\begin{equation}\label{eq:finite-root-radius}
 2c_*R\le U(G)\le\Delta(F,\mathcal L)+E_{n,k}R,\qquad
 E_{n,k}=2\pi m\left(\alpha_k+
            \frac{4\beta_k}{2n-1}\left(\frac{95}{117}\right)^n\right).
\end{equation}
For the first inequality, mark a root \(x\) of modulus \(R\).
Averaging \eqref{eq:pair-mass-sharp} for the pairs with mass one at
\(x\) and mass \(m-1\) at each of the other roots \(y\) gives
\(U(G)\ge (2c_*/m)\sum_y|x-y|\). Since the roots sum to zero,
\(\sum_y|x-y|\ge m|x|\). This argument applies also when \(m=2\).
We use the rational lower bound \(c_0^*=913/1000<c_*\).

\begin{prop}\label{prop:finite-line}
Let \(28\le n\le255\). If \(F\) is monic of degree \(n\) and its
finite critical values have diameter \(d>0\), then, for every
affine line \(\mathcal L\subset\C\),
\[
 \Delta(F,\mathcal L)>\frac45\left(\frac d2\right)^{1/n}.
\]
\end{prop}

\begin{proof}
\emph{Step 1. The coefficients determined by a derivative.}
As in Theorem~\ref{thm:line256}, normalize \(d=2\), take
\(\mathcal L=\R\), and center \(F\). Suppose that
\(\Delta(F,\R)\le L_0=4/5\). Choose \(k\) as in
Table~\ref{tab:finite-interpolation-parameters} below.
\par\nopagebreak[4]
\begin{table}[H]
\centering\small
\begin{tabular}{@{}lcccccc@{}}
\toprule
\(n\) & \(k\) & \(r\) & \(s\) & \(H\) & \(\gamma\) & \(\sigma\)\\
\midrule
\(28\)--\(30\) & 7 & \(1\) & \(11/50\) & \(9/50\) & \(1/25\) & \(65/64\)\\
\(31\) & 7 & \(26/25\) & \(17/100\) & \(4/25\) & \(1/25\) & \(65/64\)\\
\(32\)--\(36\) & 8 & \(24/25\) & \(3/20\) & \(4/25\) & \(3/200\) & \(65/64\)\\
\(37\) & 8 & \(24/25\) & \(3/20\) & \(4/25\) & \(3/200\) & \(1\)\\
\(38\)--\(255\) & \(\lceil31\sqrt n/20\rceil\) & \(97/100\)
 & \(3/20\) & \(7/50\) & \(1/25\) & \(1\)\\
\bottomrule
\end{tabular}
\par\smallskip
\caption{Parameters for interpolation on real slices.}
\label{tab:finite-interpolation-parameters}
\end{table}

In every
case \(m=n-k\ge5\) and \(E_{n,k}<2c_0^*\). By
\eqref{eq:finite-root-radius}, with \(A=L_0/(2c_0^*-E_{n,k})\),
we have \(R\le A\) and \(U(G)\le L_0+E_{n,k}A=2c_0^*A\).
Thus Lemma~\ref{lem:ordered-radii} applies with this \(A\).
Put \(S=1000A/433\), \(D=(S+A)/(m-1)\), and
\(p_i=A^i+(S-A)D^{i-1}\) for \(i\ge1\).

Write \(G(z)=\sum_{j=0}^m g_jz^{m-j}\), and decompose \(F=T+Q\)
as in \eqref{eq:top256}, so that \(\deg Q<k\). Set \(g_0^*=1\),
\(g_1^*=0\), and define
\begin{equation}\label{eq:finite-coefficients}
 \begin{aligned}
 g_j^*&=\min\left\{\binom mj(S/m)^j,\,
              \frac1j\sum_{i=2}^j p_i g_{j-i}^*\right\}
                                      &&(2\le j\le m),\\
 t_j&=\frac{\binom nj}{\binom mj}\,g_j^*
                                      &&(0\le j\le m).
 \end{aligned}
\end{equation}
Then \(|g_j|\le g_j^*\) and \(|[z^{n-j}]T|\le t_j\).
Indeed, Maclaurin's inequality for the root moduli gives the first
bound in the minimum. For the second, the first power sum vanishes
by centering, and the \(i\)-th power sum has modulus at most \(p_i\)
by \eqref{eq:ordered-moments}. Newton's identities give the bound
by induction. The assertion about \(T\) follows from
\eqref{eq:top256}.

\emph{Step 2. Counting zeros at several derivative orders.}
Fix \(s,H>0\). On the real axis write
\(H_{\theta,h}(X)=\Im F(e^{i\theta}(X+ih))\), and denote its
extension as a real-coefficient polynomial by the same symbol.
Choose \(\theta\) uniformly in \((0,2\pi)\) and \(h\) uniformly,
and independently, in \((0,H)\). Let \(N\) count the distinct real
zeros of \(H_{\theta,h}\). Let \(V\) be the decrease in a continuous
argument of \(G(e^{i\theta}(X+ih))\) as \(X\) runs from \(-s\)
to \(s\), divided by \(\pi\). We omit lines through a root of \(G\)
and the finitely many directions with \(\sin(n\theta)=0\). A single
root contributes at most one to \(V\).

For a root \(\xi\) of modulus \(b\), its mean contribution at almost
every fixed height is
\[
 \frac1\pi\int_{-s}^s
       \frac{h}{X^2+h^2}\,{\bf1}_{\{X^2+h^2>b^2\}}\,dX.
\]
Indeed, the circular mean of \((X+ih-e^{-i\theta}\xi)^{-1}\) is
\((X+ih)^{-1}\) when \(|X+ih|>b\), and zero when \(|X+ih|<b\).
Taking negative imaginary parts gives the formula. Fubini's
theorem applies jointly in \(X,h,\theta\), since the absolute
value of the kernel is a translate of \(1/|X+ih|\), which is
locally integrable in the plane.

Compared with a root at zero, the loss at height \(0<h<b\) is
at most \(2\arccos(h/b)/\pi\), obtained by integrating over
\(|X|\le\sqrt{b^2-h^2}\). At \(h>b\) the loss is zero. Its
integral in \(h\) is therefore at most \(2b/\pi\).
For a root at zero the contribution is \(2\arctan(s/h)/\pi\).
It follows from \(\sum|\xi_i|\le S\) that
\[
 \mathbb E(m-V)\le m\chi(H/s)+\frac{2S}{\pi H},
 \qquad
 \chi(x)=\frac2\pi\int_0^1\arctan(xv)\,dv.
\]
Direct integration gives \(\pi x\chi(x)=2x\arctan x-\log(1+x^2)\).
By \eqref{eq:new-crofton},
\(\mathbb E(n-N)\le L_0/(\pi H)\). Here a full turn with positive
heights parametrizes the same unoriented lines as a half turn
with heights of both signs.

Fix \(0<\gamma<1\), and put
\begin{equation}\label{eq:finite-count}
 c=m-k-1-\frac{m\chi(H/s)+(2S+L_0)/(\pi H)}{1-\gamma}.
\end{equation}
The random variable \(m-V+n-N\) is nonnegative. Markov's inequality
and Fubini's theorem give one direction \(\theta\) and a measurable
set \(\mathcal H\subset(0,H)\), with \(|\mathcal H|\ge\gamma H\),
on which this variable is at most the fraction subtracted in
\eqref{eq:finite-count}.

Fix \(h\in\mathcal H\). The \(k\)-th derivative of \(H_{\theta,h}\)
is \((n!/m!)\Im(e^{ik\theta}G(e^{i\theta}(X+ih)))\) on the real
axis. Each multiple of \(\pi\) strictly between the endpoint
arguments of \(e^{ik\theta}G\) gives a distinct zero in \((-s,s)\).
There are at least \(V-1\) such zeros. For \(0\le j<k\), Rolle's
theorem gives at least \(N-j\) real zeros of \(H_{\theta,h}^{(j)}\)
on the whole line. On each open ray outside \([-s,s]\), the number
of its zeros is at most the corresponding count for
\(H_{\theta,h}^{(k)}\) plus \(k-j\). As \(H_{\theta,h}^{(k)}\)
has degree \(m\), the number of distinct zeros of
\(H_{\theta,h}^{(j)}\) in \([-s,s]\) is at least
\[
 (N-j)-\{m-(V-1)+2(k-j)\}
 =m-k-1+j-(m-V+n-N)\ge c+j .
\]
The parameters below give \(c>0\). We may thus use any integer
\(1\le M\le\min(k,\lceil c\rceil)\), and retain \(M+j\) zeros
at each order \(j<k\).

\emph{Step 3. Interpolation and recovery of \(Q'\).}
Set \(v=(s^2+H^2)^{1/2}\) and
\[
 B=\min\left\{\sum_{i=0}^m g_i^*v^{m-i},\
          (v+S/m)^m,\ (v^2+p_2/m)^{m/2}\right\}.
\]
Then \(|G(z)|\le B\) for \(|z|\le v\). The first bound follows
from \eqref{eq:finite-coefficients}. The other two follow by
applying the arithmetic--geometric mean inequality to
\(|z-\xi_i|\) and \(|z-\xi_i|^2\), respectively. For the latter,
centering gives
\(\sum_i|z-\xi_i|^2=m|z|^2+\sum_i|\xi_i|^2\le m v^2+p_2\).

For each \(j<k\), start at \(i=j\). At order \(i\), take
\(l=\min(M+i,k-i)\) and pass to order \(i+l\). Denote the successive
step sizes by \(l_1,\ldots,l_e\), so \(\sum l_a=k-j\). Define
\[
 d_j=\frac{s^{l_1}(2s)^{k-j-l_1}}{\prod_{a=1}^e l_a!},
 \qquad
 K_j=\frac{n!}{m!\,j!}\,d_jB+
       \sum_{i=0}^m t_i\binom{n-i}{j}H^{n-i-j}.
\]
If a real polynomial \(f\) has distinct zeros \(x_1,\ldots,x_l\)
in \([-s,s]\), the interpolation remainder gives
\[
 |f(X)|\le\frac{\prod_{a=1}^l|X-x_a|}{l!}
                 \sup_{[-s,s]}|f^{(l)}|,\qquad -s\le X\le s.
\]
The numerator is at most \(s^l\) at \(X=0\) and at most \((2s)^l\)
throughout the interval. Applying this bound at the successive
orders, the first time at zero, gives
\(|H_{\theta,h}^{(j)}(0)|\le d_j n!B/m!\).

Write \(P_h\) for the real-coefficient polynomial that agrees with
\(\Im Q(e^{i\theta}(X+ih))\) on the real axis. For a term of
degree \(n-i\) in \(T\), the modulus of the \(j\)-th Taylor
coefficient of its slice is at most
\(t_i\binom{n-i}{j}H^{n-i-j}\). Subtracting the slice of \(T\)
from that of \(F\) therefore gives
\begin{equation}\label{eq:finite-low-coeff}
 |[X^j]P_h|\le K_j,\qquad 0\le j<k,\quad
 h\in\overline{\mathcal H}.
\end{equation}
The extension to the closure follows from continuity of each
coefficient in \(h\).

Take \(h_0=\min\overline{\mathcal H}\) and
\(h_1=\max\overline{\mathcal H}\). Then
\(\gamma H\le\eta:=h_1-h_0\le H\). Put
\(A_0(X)=Q(e^{i\theta}(X+ih_0))\) and write
\(A_0^\#(X)=\overline{A_0(\bar X)}\). The reflected terms in
\(P_{h_1}(X)\) and \(P_{h_0}(X-i\eta)\) are both
\(-A_0^\#(X-i\eta)/(2i)\). Hence
\[
 \frac{A_0(X+i\eta)-A_0(X-i\eta)}{2i\eta}
 =\frac{P_{h_1}(X)-P_{h_0}(X-i\eta)}{\eta}.
\]
Applying \(\eta\partial_X/\sin(\eta\partial_X)\) recovers
\(A_0'(X)\), as in the proof of \eqref{eq:lowpoly256}. Its Taylor
series acts by a finite sum on polynomials.

Let \(\alpha_l=[w^l](w/\sin w)\), fix \(r>0\), and put \(u=r+H\).
For \(X=e^{-i\theta}z-ih_0\) and \(|z|\le r\), both \(X\) and
\(X-i\eta=e^{-i\theta}z-ih_1\) have modulus at most \(u\).
Use \(\eta^{-1}\le1/(\gamma H)\) for the constant term of the
operator and \(\eta^{l-1}\le H^{l-1}\) for its positive orders.
Since its odd coefficients vanish, \eqref{eq:finite-low-coeff} gives
\begin{equation}\label{eq:finite-V}
 \|Q'\|_r\le V_0,\qquad
 V_0=2\sum_{j=0}^{k-1}K_j
 \left\{\frac{u^j}{\gamma H}
 +\sum_{\substack{2\le l\le j\\ l\ {\rm even}}}
 |\alpha_l|\frac{j!}{(j-l)!}H^{l-1}u^{j-l}\right\}.
\end{equation}

\emph{Step 4. Locating the critical points and their values.}
Choose \(a\ge0\) with \(w_1=aA/r<1\), and put \(w_2=aD/r\).
Since \(D<A\), also \(w_2<1\). Set
\[
 \begin{aligned}
 y&=w_1-\frac{(S-A)D(a/r)^2}{2(1-w_2)},\\
 Z&=(1-w_1)e^y-
 \sum_{j=2}^m\frac{g_j^*}{r^j}
 \left|\frac{\binom{n-1}{j}}{\binom mj}-a^j\right|.
 \end{aligned}
\]
On \(|z|=r\), the product \(\prod_i(1-a\xi_i/z)\) has modulus at
least \((1-w_1)e^y\). To see this, expand the logarithm of its
modulus. The first power sum is zero. For \(j\ge2\) use
\eqref{eq:ordered-moments}. The contribution to the upper bound
for the negative logarithm from the terms \(A^j\) sums to
\(-\log(1-w_1)-w_1\). In the remaining terms, \(1/j\le1/2\)
gives the geometric sum \((S-A)D(a/r)^2/(2(1-w_2))\).

The comparison polynomial
\(B_a(z)=nz^{k-1}\prod_i(z-a\xi_i)\) has all its zeros in
\(|z|<r\). Differentiating \eqref{eq:top256} shows that the
sum subtracted in \(Z\) bounds
\(|T'(z)-B_a(z)|/(nr^{n-1})\) on the circle. Thus, if \(Z>0\)
and \(V_0<nr^{n-1}Z\), \eqref{eq:finite-V} gives
\(|F'-B_a|<|B_a|\) there. Rouch\'e's theorem places every zero
of \(F'\) in \(|z|<r\).

At a critical point \(\zeta\), subtracting \(\zeta F'(\zeta)/n\)
from \(F(\zeta)\) removes the leading term. Since \(t_1=0\) and
\(|Q(\zeta)-Q(0)|\le rV_0\), we obtain
\begin{equation}\label{eq:finite-value}
 |F(\zeta)-Q(0)|\le Y_0,\qquad
 Y_0=\sum_{j=2}^m\frac jn t_jr^{n-j}+r(1+1/n)V_0 .
\end{equation}
In particular \(Y_0<1\) contradicts the critical-value diameter
\(d=2\).

\emph{Step 5. Verification of the scalar inequalities.}
Table~\ref{tab:finite-interpolation-parameters} specifies all the
parameters, with \(a=\sigma(n-1)/m\).

For each of the \(228\) degrees, two independent computations
verify \(E_{n,k}<2c_0^*\), \(w_1<1\), and
\begin{equation}\label{eq:finite-interpolation-tests}
 c>\frac3{100},\qquad Z>\frac{19}{1000},\qquad
 \frac{V_0}{nr^{n-1}Z}<\frac9{100},\qquad Y_0<\frac{71}{100}.
\end{equation}
The program \href{https://github.com/zhangteng2000/ehp-lemniscate-verification-20261108/blob/ca8c451bfbc7312b9c57acc5e672c0b69cd601ae/verification/verify_interpolation_degrees.py}{\nolinkurl{verify_interpolation_degrees.py}} evaluates
the arctangent and logarithm in the expression for \(\chi\) by
rational series bounds. It computes the Taylor coefficients of
\(w/\sin w\) in \eqref{eq:finite-V} recursively.

The second implementation is \href{https://github.com/zhangteng2000/ehp-lemniscate-verification-20261108/blob/ca8c451bfbc7312b9c57acc5e672c0b69cd601ae/verification/independent_interpolation_check.py}{\nolinkurl{independent_interpolation_check.py}}.
For \eqref{eq:finite-count} it bounds the integral defining \(\chi\)
from above by midpoint sums on \(64\) equal intervals, using
concavity of the arctangent on \([0,\infty)\). It takes \(M=1\) and
\(|\alpha_l|\le6/2^l\), since on \(|w|=2\) one has
\(|1-\sin w/w|\le\sinh(2)/2-1<5/6\). Both implementations use
rational bounds with outward rounding. The strict inequalities
\eqref{eq:finite-interpolation-tests} establish the contradiction
in \eqref{eq:finite-value} for every degree in the table.
\end{proof}

\subsection{The circle comparison in degrees 28 through 255}
\label{sec:finite-circle}

Proposition~\ref{prop:finite-line} supplies the hypothesis of
Proposition~\ref{prop:uniform-circle-transfer} in this range.
Choose a centered maximizer whose finite critical values lie on
\(\mathbb T\), as permitted by Lemma~\ref{lem:extremizers}.
The circle comparison forces these values to coincide.
By Lemma~\ref{lem:critical-compactness}, the chosen maximizer is
\(z^n+v\) for some \(|v|=1\). Its lemniscate has length
\(\ell(\partial E_1(p_0))\), which proves the length inequality
for every monic polynomial of these degrees.
Section~\ref{sec:equality} determines all cases of equality.

\section{Critical values on a line}\label{sec:collinear-gram}

The boundary term in \eqref{eq:new-mass} involves polynomials
whose critical values all lie on the target line. We now use this
additional hypothesis. Successive differentiation will bound the
coefficients down to degree two. To control the linear coefficient,
we shall compare two moments of the critical points and their values.
Throughout this section \(n\ge3\).

\subsection{Bounds for the coefficients of a derivative}

For a centered monic polynomial \(G(z)=\sum_{j=0}^m g_jz^{m-j}\),
we seek constants \(e_{m,j}\) such that
\begin{equation}\label{eq:peel-caps}
 U(G)\le2c_*A
 \quad\Longrightarrow\quad |g_j|\le e_{m,j}A^j
 \quad(0\le j\le m,\ A>0).
\end{equation}
We take \(e_{m,0}=1\) and \(e_{m,1}=0\). In degree two the roots
are \(x,-x\), and \eqref{eq:pair-mass-sharp} is an equality for
equal masses. Thus \(U(G)=2c_*|x|\), and we may take \(e_{2,2}=1\).

We give the recursion for \(m\ge3\). All quantities in a single
step depend only on \(m\), the preceding row \(e_{m-1,i}\), and the
endpoints of one of the intervals
\([l,h]=[a/24,(a+1)/24]\), \(0\le a<24\), with \(a\) an integer.
Put \(\nu=1000/433\). Let \(\mathcal T_m\) consist of the members of
\[
 \{0,1/4,1/2,1,2,3,4,5,7,9,12,16\}
 \ \cup\ \{(m-1)/4,(m-1)/2,3(m-1)/4\}
\]
that are strictly less than \(m-1\). For the interval \([l,h]\),
define
\begin{equation}\label{eq:remaining-root-energy}
 A'=\min_{t\in\mathcal T_m}
       \frac{1-mtl/((m-1)(1+t))}{1-t/(m-1)} .
\end{equation}
For \(2\le j\le m\) put
\begin{equation}\label{eq:translated-coefficient-bound}
 T_j=\sum_{i=0}^{\min(j,m-1)}
 e_{m-1,i}(A')^i\left(\frac h{m-1}\right)^{j-i}
 \left|\binom{m-1-i}{j-i}-(m-1)\binom{m-1-i}{j-i-1}\right|.
\end{equation}
A binomial coefficient with lower index outside its usual range is
understood to be zero.

For a second bound, put \(a_0=0\). For \(1\le i\le m-1\), let
\(a_i\) be the minimum of \(ih\), \(\nu-l\), and the following
entries for which the indicated conditions hold.
\[
 1+\frac{2i}{m-i-1}-l\quad(i<m-1),\qquad
 1+\frac{i(\nu+1)}{m-1}-l\quad(m\ge5).
\]
Let \(\widehat a\) be the least concave majorant of the points
\((i,a_i)\) on \([0,m-1]\), taken piecewise linear, and put
\(\lambda_i=\widehat a(i)-\widehat a(i-1)\) for \(1\le i\le m-1\).
Since the \(a_i\) are nondecreasing, these slopes are nonnegative
and nonincreasing. Define \(p_q=\sum_{i=1}^{m-1}\lambda_i^q\)
and set
\[
 H_0=d_0=1,\quad H_1=d_1=0,\qquad
 H_j=\frac1j\sum_{q=2}^j p_qH_{j-q},\quad
 d_j=\frac{j-1}{j!}h^j\quad(2\le j\le m).
\]

\begin{lem}\label{lem:peel-recursion}
If \eqref{eq:peel-caps} holds in degree \(m-1\), it holds in
degree \(m\) with
\begin{equation}\label{eq:peel-envelope}
 e_{m,j}=\max_{[l,h]}
 \min\left\{T_j,\ \sum_{q=0}^j d_qH_{j-q},\
                 \binom mj\min(h,\nu/m)^j\right\},
 \qquad 2\le j\le m,
\end{equation}
where the maximum runs over the \(24\) closed intervals above.
\end{lem}

\begin{proof}
By scaling, take \(A=1\) in \eqref{eq:peel-caps}. Mark one root
\(x\) of largest modulus, and use \(y\) for the remaining roots,
with all labels retained except the one marked. The marked-root
bound used in \eqref{eq:finite-root-radius} gives \(|x|\le1\),
and \eqref{eq:central-absolute} gives
\(|x|+\sum_y|y|\le\nu\). Fix an interval with \(l\le|x|\le h\).

\emph{Step 1. Centering the remaining roots.}
Their mean is \(-x/(m-1)\). Thus
\(G(z)=(z-x)J(z+x/(m-1))\), where \(J\) is centered and monic
of degree \(m-1\). For \(0\le t<m-1\), the measure
\(\delta_x+\frac{t}{m-1}\sum_y\delta_y\) is the average of
\(\delta_x+t\delta_y\) over the remaining labels. By concavity and
\eqref{eq:pair-mass-sharp}, its energy is at least
\[
 \frac{2c_*t}{1+t}\frac1{m-1}\sum_y|x-y|
 \ge\frac{2c_*tm|x|}{(1+t)(m-1)}.
\]
The unused measure has energy \((1-t/(m-1))U(J)\), by translation
invariance and homogeneity. Superadditivity and \(U(G)\le2c_*\)
therefore give \(U(J)\le2c_*A'\), with \(A'\) as in
\eqref{eq:remaining-root-energy}. Expanding
\((z-x)J(z+x/(m-1))\) and using the bounds in degree \(m-1\)
gives precisely \(|g_j|\le T_j\) in
\eqref{eq:translated-coefficient-bound}.

\emph{Step 2. Moments of the remaining roots.}
In decreasing order of modulus, the first \(i\) remaining roots
have total modulus at most \(a_i\). The bounds \(ih\) and
\(\nu-l\) follow from the choice of \(x\). Marking \(i+1\) roots
in \eqref{eq:ordered-prefix} gives the bound with denominator
\(m-i-1\). The proof of that prefix inequality applies for every
\(m\ge2\). For \(m\ge5\), the proof of
Lemma~\ref{lem:ordered-radii} also bounds the sum of the first
\(i+1\) moduli by \(1+i(\nu+1)/(m-1)\). In each case, subtract
\(|x|\ge l\).

Consequently the partial sums of the \(\lambda_i\) dominate those
of the ordered moduli \(|y|\). For \(\tau\ge0\), the sum
\(\sum_y(|y|-\tau)_+\) is the maximum, over prefix lengths \(i\),
of the \(i\)-th partial sum minus \(i\tau\), including \(i=0\).
It follows that this sum is at most \(\sum_i(\lambda_i-\tau)_+\).
Integrating against \(q(q-1)\tau^{q-2}\,d\tau\) gives
\(\sum_y|y|^q\le p_q\) for every integer \(q\ge2\).

Since \(\sum_y y=-x\), the reversed polynomial has the expansion
\[
 w^mG(1/w)
 =(1-xw)e^{xw}
   \exp\!\left(-\sum_{q\ge2}\frac{w^q}{q}\sum_y y^q\right).
\]
The coefficient moduli of the first factor are bounded by \(d_j\).
Those of the exponential are bounded by the coefficients of
\(\exp(\sum_{q\ge2}p_q w^q/q)\), which satisfy the recurrence
defining \(H_j\). Hence \(|g_j|\le\sum_{q=0}^j d_qH_{j-q}\).
Finally, Maclaurin's inequality gives
\(|g_j|\le\binom mj\min(h,\nu/m)^j\).
Taking the minimum on each interval and the maximum over the
intervals proves \eqref{eq:peel-envelope}. The use of closed
intervals includes endpoints and ties for the largest root modulus.
\end{proof}

We next pass from these coefficient bounds for \(G\) to bounds
for \(F\). Let \(F(z)=\sum_{i=0}^n c_i z^i\) be centered and monic,
and suppose \(\Delta(F,\R)\le L_1\), where \(L_1>0\).
Fix \(2\le k\le n-2\), put \(m=n-k\), and let \(G=m!F^{(k)}/n!\).
For \(2\le j\le m\) and each available bound \(|c_{n-j}|\le B_{n-j}\),
put \(q_j=B_{n-j}\binom mj/\binom nj\). Set \(q_j=+\infty\) if
there is no previous bound for that coefficient.

\begin{lem}\label{lem:radial-descent}
All roots of \(G\) have modulus at most \(r\) if \(r>0\) satisfies
\(r(2c_0^*-E_{n,k})\ge L_1\) or
\begin{equation}\label{eq:radial-descent}
 \sum_{j=2}^m
 \min\left\{\frac{q_j}{r^j},\
       e_{m,j}\left(\frac{L_1/r+E_{n,k}}{2c_0^*}\right)^j
     \right\}<1 .
\end{equation}
\end{lem}

\begin{proof}
Let \(R\) be the largest root modulus. The first test follows
from \eqref{eq:finite-root-radius}. For the second, suppose that
\(R\ge r\) and take a root of modulus \(R\). By
\eqref{eq:finite-root-radius} and \eqref{eq:peel-caps}, each
\(|g_j|\) is at most
\(e_{m,j}((L_1+E_{n,k}R)/(2c_0^*))^j\), as well as \(q_j\).
Evaluating \(G\) at the chosen root and dividing by \(R^m\)
shows that the left side of \eqref{eq:radial-descent} at \(R\)
is at least one. Each summand is nonincreasing in the radius,
which contradicts \eqref{eq:radial-descent} at \(r\).
\end{proof}

With \(r\) verified, the coefficient identity for \(G=m!F^{(k)}/n!\)
gives
\begin{equation}\label{eq:descending-update}
 |c_{n-j}|\le
 \frac{\binom nj}{\binom mj}\,e_{m,j}
 \left(\frac{L_1+E_{n,k}r}{2c_0^*}\right)^j,
 \qquad 2\le j\le m.
\end{equation}
In particular, \(j=m\) gives a bound for \(c_k\). Starting with
\(c_n=1\) and \(c_{n-1}=0\), we apply
Lemma~\ref{lem:radial-descent} successively for
\(k=n-2,n-3,\ldots,2\), keeping the smaller bound whenever a
coefficient has already been estimated. A further pass through the
same orders uses all bounds established in the previous pass.

\subsection{Two moments of the critical values}

Write \(m=n-1\), and let \(v_1,\ldots,v_m\) be the finite critical
values of \(F\), listed with the multiplicities of the corresponding
zeros of \(F'\).
With \(v_{\rm av}=m^{-1}\sum_i v_i\), define the second centered
power sum
\[
 \mathcal V(F)=\sum_{i=1}^m(v_i-v_{\rm av})^2.
\]
If the critical values lie on a line and have diameter two, then
\(|\mathcal V(F)|\ge2\). Indeed, after a rotation their centered
values are real, and the squares of the two endpoint distances
from the mean sum to at least half the squared endpoint separation.

\begin{lem}\label{lem:companion-variance}
Let \(F(z)=z^n+\sum_{j=0}^{n-2}c_jz^j\). Given bounds
\(B_j\ge0\) for \(1\le j\le n-2\), the maximum of
\(|\mathcal V(F)|\) over \(|c_j|\le B_j\) is
\(\mathcal V(z^n-\sum_{j=1}^{n-2}B_jz^j)\).
\end{lem}

\begin{proof}
An additive constant has no effect on \(\mathcal V\). Set \(c_0=0\)
and write \(c_j=-x_j\), initially treating the \(x_j\) as formal
variables. Let \(C\) be the companion matrix of \(F'/n\), with
first row \(((n-i)x_{n-i}/n)_{i=1}^m\), where \(x_{n-1}=0\),
ones on the subdiagonal, and zeros elsewhere. Put
\(M=\sum_{j=1}^{n-2}(1-j/n)x_jC^j\). At a critical point \(\zeta\),
subtracting \(\zeta F'(\zeta)/n\) removes the leading term of \(F(\zeta)\).
It follows that the eigenvalues of \(-M\) are the \(v_i\), with
their algebraic multiplicities. The trace identity
\(\mathcal V(F)=\operatorname{tr}(M^2)-(\operatorname{tr}M)^2/m\)
therefore gives
\begin{equation}\label{eq:companion-variance}
 \mathcal V(F)=\sum_{i\ne j}M_{ij}M_{ji}
       +\frac1m\sum_{i<j}(M_{ii}-M_{jj})^2.
\end{equation}

Every entry of \(M\) has nonnegative coefficients. The same holds
for \(M_{ii}-M_{jj}\) whenever \(i<j\). To see this for \(C^q\)
and adjacent indices, regard \(C_{ab}\) as the weight of the
directed edge from \(a\) to \(b\). A closed walk of length \(q\)
starting at \(i+1\) begins with the edge \(i+1\to i\), of weight
one. Moving that edge to the end gives a closed walk starting at
\(i\). The map is injective and preserves the weight. The assertion
for \(M\) follows by summing over \(q\), and the general diagonal
difference follows by telescoping. Thus \eqref{eq:companion-variance}
is a polynomial with nonnegative coefficients in the \(x_j\).
Its value at \(x_j=B_j\) bounds its modulus on the entire polydisk
and is attained there.
\end{proof}

For a large linear coefficient we use collinearity a second time.
We first state the form of Gon\c{c}alves' inequality
\cite{Gon50} needed to control the moduli of the critical points.

\begin{lem}[{\cite[Eq.~(1.4)]{BMV07}}]\label{lem:gram-root-product}
Let \(P(z)=z^m+\sum_{j=0}^{m-1}s_jz^j\) have degree \(m\ge1\)
and \(|s_0|=1\), and put
\(J=\prod_{P(\xi)=0}\max\{1,|\xi|\}\), with roots counted with
multiplicity. Then
\[
 J^2+J^{-2}\le 2+\sum_{j=1}^{m-1}|s_j|^2.
\]
\end{lem}

\begin{proof}
In the factorization of \(P\), replace \(z-\xi\) by
\(1-\bar\xi z\) for each \(|\xi|>1\). The resulting polynomial
has the same modulus as \(P\) on \(\mathbb T\), leading coefficient
of modulus \(J\), and constant coefficient of modulus \(J^{-1}\).
The last assertion uses \(\prod|\xi|=1\). Parseval's identity,
retaining these two coefficients, gives the inequality.
\end{proof}

Fix \(u>0\), and suppose \(|c_j|\le B_j\) for \(2\le j\le n-2\).
Set \(t_0=(u/n)^{1/m}\), \(A_j=B_jt_0^{j-n}\),
\(d_{n-j}=jA_j/n\) for \(2\le j\le n-2\), and \(d_1=0\),
\(d_m=1\). Write \(P_l\) for the \(l\)-th power sum of the
roots of \(z^m-\sum_{i=1}^m d_i z^{m-i}\). Newton's identities
take the form
\begin{equation}\label{eq:gram-power-sums}
 \sum_{l\ge1}P_lw^l
   =\frac{\sum_{i=1}^m i d_iw^i}{1-\sum_{i=1}^m d_iw^i}.
\end{equation}
In particular \(P_l\ge0\). Put \(B=(\sum_{i=1}^{m-1}d_i^2)^{1/2}\)
and choose \(r>1\) with \(\sum_{i=1}^m d_i r^{-i}<1\).

\begin{lem}\label{lem:collinear-gram}
Under these assumptions, if the critical values of \(F\) lie on
an affine line and \(|c_1|\ge u\), then
\begin{equation}\label{eq:collinear-gram}
 m(m-P_2)\le
 \sum_{j=2}^{n-2}(1-j/n)A_j
 \left\{P_{j-1}+P_{j+1}
       +B\left(1+\frac{(r+1)^2r^{j-1}}4\right)\right\}.
\end{equation}
\end{lem}

\begin{proof}
\emph{Step 1. The rescaled critical points.}
Put \(t=(|c_1|/n)^{1/m}\). Denote the critical points of
\(t^{-n}F(tz)\) by \(\zeta_1,\ldots,\zeta_m\), indexed so that
\(v_i=F(t\zeta_i)\). The corresponding critical points of \(F\)
in the original coordinate are \(t\zeta_i\). The monic polynomial
\(F'(tz)/(nt^m)\) with roots \(\zeta_i\) has constant coefficient
of modulus one. Since
\(t\ge t_0\), the modulus of its coefficient of \(z^{m-i}\)
is at most \(d_i\). The choice of \(r\) then gives
\(|\zeta_i|<r\) for all \(i\). From \eqref{eq:gram-power-sums},
or directly by Newton's identities, we have
\(|\sum_i\zeta_i^l|\le P_l\). We also have
\(\sum_i\zeta_i=0\), \(\prod_i|\zeta_i|=1\), and
\(\sum_i|\zeta_i|^2\ge m\).

Let \(J=\prod_i\max\{1,|\zeta_i|\}\). By
Lemma~\ref{lem:gram-root-product}, \((J-J^{-1})^2\le B^2\).
Since \(J-J^{-1}=2\sinh(\log J)\ge2\log J\), this gives
\begin{equation}\label{eq:gram-log-moduli}
 \sum_{|\zeta_i|>1}\log|\zeta_i|
 =\sum_{|\zeta_i|<1}-\log|\zeta_i|
 =\log J\le B/2.
\end{equation}

\emph{Step 2. Comparing the two moments.}
For \(0<s\le1\) and \(j\ge2\),
\(s^{j-1}(1-s^2)\le-2\log s\). For \(1\le s<r\), the bound
\(\log s\ge2(s-1)/(s+1)\) gives
\(s^{j-1}(s^2-1)\le r^{j-1}(r+1)^2\log s/2\).
Apply these inequalities in \eqref{eq:gram-log-moduli} to estimate
the difference between \(\sum_i\overline{\zeta_i}\zeta_i^j\)
and \(\sum_i\zeta_i^{j-1}\). We obtain
\begin{equation}\label{eq:gram-mixed-moment}
 \left|\sum_i\overline{\zeta_i}\zeta_i^j\right|
 \le P_{j-1}+B\left(1+\frac{(r+1)^2r^{j-1}}4\right).
\end{equation}

Write \(w_i=t^{-n}(F(t\zeta_i)-c_0)\). Eliminating
\((t\zeta_i)^n\) with \(F'(t\zeta_i)=0\) gives
\[
 w_i=m\frac{c_1}{|c_1|}\zeta_i+
       \sum_{j=2}^{n-2}(1-j/n)c_jt^{j-n}\zeta_i^j .
\]
The \(w_i\) are collinear. Thus \(w_i=w_*+e^{i\psi}x_i\)
for real \(x_i\), and \(\sum_i\zeta_i=0\) gives
\(\big|\sum_i\overline{\zeta_i}w_i\big|
  =\big|\sum_i\zeta_iw_i\big|\).
Separating the linear term on the two sides yields
\[
 m\left(\sum_i|\zeta_i|^2-\left|\sum_i\zeta_i^2\right|\right)
 \le\sum_{j=2}^{n-2}(1-j/n)A_j
 \left(\left|\sum_i\overline{\zeta_i}\zeta_i^j\right|
       +\left|\sum_i\zeta_i^{j+1}\right|\right).
\]
Using \(\sum_i|\zeta_i|^2\ge m\), the power-sum bounds from
Step~1, and \eqref{eq:gram-mixed-moment} proves
\eqref{eq:collinear-gram}.
\end{proof}

\begin{prop}\label{prop:collinear-finite}
Let \(15\le n\le27\), and let \(F\) be monic of degree \(n\).
If all its finite critical values lie on a line \(\mathcal L\) and
have diameter \(d>0\), then
\[
 \Delta(F,\mathcal L)>\frac{49}{60}\left(\frac d2\right)^{1/n}.
\]
\end{prop}

\begin{proof}
Normalize \(d=2\), take \(\mathcal L=\R\), and center \(F\).
Suppose \(\Delta(F,\R)\le L_1=49/60\). Construct the array
\(e_{m,j}\) by Lemma~\ref{lem:peel-recursion}. Four passes through
the derivative orders in \eqref{eq:descending-update} give upper
bounds \(B_j\) for \(|c_j|\), \(2\le j\le n-2\). Each radius in
\href{https://github.com/zhangteng2000/ehp-lemniscate-verification-20261108/blob/ca8c451bfbc7312b9c57acc5e672c0b69cd601ae/verification/peeled_primary.json}{\texttt{peeled\_primary.json}} is verified by
Lemma~\ref{lem:radial-descent} using only bounds available before
that update.

Use \(u,r\) from Table~\ref{tab:collinear-gram-parameters}.
\par\nopagebreak[4]
\begin{table}[H]
\centering\small
\begin{tabular}{@{}crr@{\qquad}crr@{}}
\toprule
\(n\)&\(u\)&\(1000r\)&\(n\)&\(u\)&\(1000r\)\\
\midrule
15&\(7/15\)&1373&22&\(13/60\)&1153\\
16&\(9/20\)&1323&23&\(11/60\)&1137\\
17&\(5/12\)&1285&24&\(3/20\)&1127\\
18&\(23/60\)&1249&25&\(2/15\)&1113\\
19&\(1/3\)&1219&26&\(7/60\)&1103\\
20&\(17/60\)&1195&27&\(1/10\)&1095\\
21&\(1/4\)&1171&&&\\
\bottomrule
\end{tabular}
\par\smallskip
\caption{The cutoff \(u\) for \(|c_1|\) and the radius \(r\)
in Lemma~\ref{lem:collinear-gram}.}
\label{tab:collinear-gram-parameters}
\end{table}

If \(|c_1|\le u\), set \(B_1=u\). Lemma~\ref{lem:companion-variance}
then gives \(|\mathcal V(F)|/2<94/100\), contradicting the lower
bound \(2\) for \(|\mathcal V(F)|\). If \(|c_1|\ge u\), the
computed \(d_i\) satisfy \(\sum_i d_i r^{-i}<1\) and \(m-P_2>0\).
The right-hand side of \eqref{eq:collinear-gram}, divided by
\(m(m-P_2)\), is less than \(84/100\), contradicting the lemma.
The two tests include the common endpoint \(|c_1|=u\).

For each of the thirteen degrees,
\href{https://github.com/zhangteng2000/ehp-lemniscate-verification-20261108/blob/ca8c451bfbc7312b9c57acc5e672c0b69cd601ae/verification/verify_peeled_line.py}{\nolinkurl{verify_peeled_line.py}} evaluates the first test by Newton
power sums. The independent implementation
\href{https://github.com/zhangteng2000/ehp-lemniscate-verification-20261108/blob/ca8c451bfbc7312b9c57acc5e672c0b69cd601ae/verification/audit_peeled_line.py}{\nolinkurl{audit_peeled_line.py}} evaluates
\eqref{eq:companion-variance} by matrix multiplication, recomputes
the bounds in \eqref{eq:peel-envelope} by multiplying power
series, and checks each of the proposed radii before using it.
In both implementations the largest upper bounds for the two tests
occur in degree \(15\).
\end{proof}

\section{Line and circle estimates in degrees nine through twenty-seven}
\label{sec:conditional-low}

We continue the argument of Section~\ref{sec:collinear-gram}.
In degrees nine through fourteen, the bounds obtained by successive
differentiation depend appreciably on how far the largest root of
each derivative lies from the origin. We therefore keep track of
that radius relative to the mean line deficit until the coefficients
of the derivative have been estimated.

\subsection{The angular mean outside a disk}

Let \(F\) be centered and monic of degree \(n\). Fix
\(2\le k\le n-2\), put \(m=n-k\), and write \(G=m!F^{(k)}/n!\).
Let \(\xi_1,\ldots,\xi_m\) be the roots of \(G\), counted with
multiplicity, and put \(R=\max_i|\xi_i|\). If \(R=0\), then
\(G(z)=z^m\) and \(U(G)=0\). For \(R>0\), write \(y_i=\xi_i/R\)
and set
\[
 B(t)=\prod_{i=1}^m(1-y_it),\qquad
 A(t)=\frac1m\sum_{i=1}^m\frac1{1-y_it}.
\]
We also put \(v=\pi/(2k)\) and \(\beta=v/\sin v\).

\begin{lem}\label{lem:conditional-contour}
Suppose \(R>0\). Let \(T>1\) and \(\eta=T^{-2}\). Suppose that
\(|A(t)|\le a\) and \(|B(t)|\ge z_0>0\) on \(|t|\le\eta\), and
that \(|B(t)|\le w\) on \(|t|\le1\). If
\(\rho=(2a)^kw/(T^{2n}z_0)<1\), then, for every affine line \(\mathcal L\),
\begin{equation}\label{eq:conditional-contour}
 \Delta(F,\mathcal L)\ge\beta U(G)-ER,\qquad
 E=2mT\left(\beta-1+
       \sqrt{a/2}\,(1-v\cot v)\frac{\rho}{2n-1}\right).
\end{equation}
\end{lem}

\begin{proof}
\emph{Step 1. The integral over the disk.}
By a source dilation, take \(R=1\), and write \(V=G'/G\).
Averaging \eqref{eq:sampling-integral} in the target angle gives
\[
 U(G)=\lim_{r\to\infty}J(r),\qquad
 J(r)=2mr-\frac1\pi\int_{|z|\le r}|V(z)|\,dA(z).
\]
For \(r>1\), the circular mean of \(zV(z)=\sum_i z/(z-y_i)\)
is \(m\). Thus the circular mean of \(|zV(z)|\) is at least \(m\),
and differentiation gives \(J'(r)\le0\). In particular,
\begin{equation}\label{eq:conditional-inner-integral}
 \int_{|z|\le T}|V(z)|\,dA(z)\le 2\pi mT-\pi U(G).
\end{equation}

Let \(\psi\) be the angle of \(\mathcal L\), and write \(\Lambda_k(\psi)\)
as in \eqref{eq:sampling-integral}. In the notation
of Lemma~\ref{lem:sampled-derivative}, the sampling function \(P_k\)
has mean \(2/\pi\) and maximum \(2/\pi+\alpha_k\).
By \eqref{eq:conditional-inner-integral}, the disk \(|z|\le T\)
contributes at most \(\alpha_k(2\pi mT-\pi U(G))/2\) to
\(U(G)-\Lambda_k(\psi)\).

\emph{Step 2. The exterior integral.}
Since \(\Re A(t)>1/2\) on \(|t|<1\), we may take analytic powers
of \(A\) normalized at zero. For \(j\ge1\), put
\[
 H_+=A^{kj+1/2}B^{-j},\qquad H_-=A^{-kj+1/2}B^j,\qquad
 H_-^\#(t)=\overline{H_-(\bar t)}.
\]
As in the proof of Lemma~\ref{lem:sampled-derivative}, the
\(j\)-th angular moment on \(|z|=r\), with its constant phase
omitted, is
\[
 \frac mr\,\frac1{2\pi i}\int_{|\omega|=1}
 \omega^{-2nj}H_+\!\left(\frac1{r\omega}\right)
 H_-^\#\!\left(\frac{\omega}{r}\right)\frac{d\omega}{\omega}.
\]
For \(r>T\), move the contour to \(|\omega|=cr\), with \(c<1\)
sufficiently close to one that \(1/(cr^2)\le\eta\). The two
arguments then have moduli \(1/(cr^2)\) and \(c\), respectively.
The first factor is bounded by \(a^{kj+1/2}z_0^{-j}\).
For the reflected factor, \(\Re A>1/2\) gives the bound
\(2^{kj-1/2}w^j\). Letting \(c\uparrow1\) in this upper bound,
we find that the angular moment has modulus at most
\begin{equation}\label{eq:conditional-exterior-mode}
 \frac mr\sqrt{\frac a2}
 \left[\rho\left(\frac Tr\right)^{2n}\right]^j.
\end{equation}
All contours used here lie strictly inside the annulus of
analyticity \(r^{-1}<|\omega|<r\).

After radial integration, the contribution of
\eqref{eq:conditional-exterior-mode}, before multiplication by
its Fourier coefficient, is at most
\(\pi mT\sqrt{a/2}\,\rho^j/(2nj-1)\). The absolute values of the
nonconstant cosine coefficients of \(P_k\) sum to \(\beta_k\).
Since \(\rho<1\), the exterior contribution to
\(U(G)-\Lambda_k(\psi)\) is at most
\(\pi mT\sqrt{a/2}\,\beta_k\rho/(2n-1)\).
Together with Step~1, the identities
\(\pi\alpha_k=2(\beta-1)\) and
\(\pi\beta_k=2(1-v\cot v)\) give \eqref{eq:conditional-contour}
for \(R=1\), since \(\Delta(F,\mathcal L)\ge\Lambda_k(\psi)\).
Restoring the source scale proves the result.
\end{proof}

\subsection{Coefficient bounds at a prescribed root radius}

Put \(D=U(G)/(2c_*)\). The bound \(2c_*R\le U(G)\) from the
marked-root argument gives \(R\le D\). For \(D>0\) put
\(\lambda=R/D\) and write \(G(z)=\sum_{j=0}^m g_jz^{m-j}\).
For \(m\ge3\) divide \([0,1]\) into \(128\) equal closed intervals.
For each such interval \([l,h]\) we construct constants \(e_j(l,h)\)
with
\begin{equation}\label{eq:conditional-coefficient-band}
 |g_j|\le D^j e_j(l,h)\qquad\text{whenever }l\le\lambda\le h.
\end{equation}
In degree two, \(U(G)=2c_*R\), so only \([l,h]=[1,1]\) is needed,
with \(e_0=e_2=1\) and \(e_1=0\).

Use the proof of Lemma~\ref{lem:peel-recursion} on the new partition.
At degree \(m\), retain the minimum in \eqref{eq:peel-envelope} on
each interval. In the translated coefficient bound \(T_j\), use the
maximum of the preceding degree's bounds over all its intervals.
One further entry in the minimum is useful. To derive it, normalize
\(D=1\), mark a root \(x\) of largest modulus, and let \(A'\) be
the bound for \(U(J)/(2c_*)\) obtained in
Lemma~\ref{lem:peel-recursion} after centering the remaining
\(m-1\) roots. Their mean modulus before centering is at most
\[
 u=\min\left\{h,\frac{\nu-l}{m-1},\frac{\nu A'+h}{m-1}\right\},
 \qquad \nu=\frac{1000}{433}.
\]
Indeed their total modulus is at most \(\nu-|x|\). Centering
gives total modulus at most \(\nu A'\), and translating back
adds at most \(|x|\). Applying Maclaurin's inequality to the
remaining roots and multiplying by \(z-x\) gives the additional
bound \(\binom{m-1}{j}u^j+h\binom{m-1}{j-1}u^{j-1}\)
for \(2\le j\le m\).

We next specify the three disk bounds in
Lemma~\ref{lem:conditional-contour}. Write \(e_j=e_j(l,h)\), and
set \(\sigma=\min\{1,\nu/(ml)\}\) if \(l>0\), and \(\sigma=1\)
if \(l=0\). Since \(\sum_i|\xi_i|^2\le R\sum_i|\xi_i|\le R\nu D\),
we have \(m^{-1}\sum_i|y_i|^2\le\sigma\). For \(\eta=T^{-2}\)
we may take
\begin{equation}\label{eq:conditional-disk}
\begin{aligned}
 z_0&=\max\left\{(1-\eta^2)^{m/2},\
                   1-\sum_{j=2}^m e_j(\eta/l)^j\right\},\\
 w&=\min\left\{(1+\sigma)^{m/2},\
                   1+\sum_{j=2}^m e_j/l^j\right\},\\
 a&=\min\left\{(1-\eta^2)^{-1},\
       1+\frac{\sigma\eta^2}{1-\eta},\
       1+\frac{\sum_{j=2}^m j e_j(\eta/l)^j}{m z_0}\right\}.
\end{aligned}
\end{equation}
At \(l=0\) omit all entries containing \(1/l\).

For the first bound on \(|B|\) from below, use
\(|1-y_it|\ge1-\Re(y_it)\). On \([-\eta,\eta]\),
\(\log(1-s)\) lies above the chord joining its endpoint values.
Since \(\sum_i y_i=0\), summing the chord inequality gives
\(\log|B(t)|\ge(m/2)\log(1-\eta^2)\) for \(|t|\le\eta\).
For \(|t|\le1\), the arithmetic--geometric mean inequality applied
to \(|1-y_it|^2\) gives \(|B(t)|\le(1+\sigma)^{m/2}\).

Since \(\Re A(t)>1/2\) for \(|t|<1\), the function
\(1-1/A(t)\) maps the unit disk into itself. Its value and first
derivative vanish at zero. Applying Schwarz's lemma twice gives
\(|1-1/A(t)|\le |t|^2\). For the second bound, subtract
\(1+y_it\) from each summand in \(A(t)\) and use the bound on
\(m^{-1}\sum_i|y_i|^2\). The remaining entries in
\eqref{eq:conditional-disk} follow from
\(B(t)=\sum_{j=0}^m g_jR^{-j}t^j\),
\eqref{eq:conditional-coefficient-band}, and
\(A=1-tB'/(mB)\).

Suppose now that \(\Delta(F,\mathcal L)\le L\), where \(L>0\), and that
\(|g_j|\le q_j\) is already known, with \(q_j=+\infty\) allowed.
For a proposed \(s>0\), either of the following conditions
suffices to give \(D\le s\) on \([l,h]\).
\begin{equation}\label{eq:conditional-tests}
 s(2c_*\beta-Eh)\ge L
 \qquad\text{or}\qquad
 l>0\ \text{ and }\
 \sum_{j=2}^m\frac{\min\{q_j/s^j,e_j(l,h)\}}{l^j}<1 .
\end{equation}
The first follows from \eqref{eq:conditional-contour} and
\(R\le hD\). For the second, if \(D\ge s\), evaluate \(G\) at
a root of modulus \(R\ge lD\). Since \(g_1=0\), division by the
leading term gives
\(1\le\sum_{j=2}^m|g_j|/R^j\). By
\eqref{eq:conditional-coefficient-band} the right-hand side is
at most the sum in \eqref{eq:conditional-tests}, a contradiction.

Write \(F(z)=\sum_{i=0}^n c_i z^i\). At each derivative order,
an existing bound \(|c_{n-j}|\le B_{n-j}\) gives
\(q_j=B_{n-j}\binom mj/\binom nj\). Once a value \(s(l,h)\) has
been verified on every interval, we obtain
\begin{equation}\label{eq:conditional-update}
 |c_{n-j}|\le\frac{\binom nj}{\binom mj}
           \max_{[l,h]} e_j(l,h)s(l,h)^j,\qquad 2\le j\le m.
\end{equation}
If \(D=0\), then \(G=z^m\) and these coefficients vanish. We apply
\eqref{eq:conditional-update} once at each order
\(k=n-2,n-3,\ldots,2\), retaining the smaller bound whenever a
coefficient has already been estimated.

\begin{prop}\label{prop:conditional-line}
Let \(F\) be centered and monic of degree \(9\le n\le14\).
Suppose that all its finite critical values lie on a line \(\mathcal L\)
and have diameter two. Then
\[
 \Delta(F,\mathcal L)>
 \begin{cases}
   41/50,& n=9,\\
   7/8,& 10\le n\le14 .
 \end{cases}
\]
\end{prop}

\begin{proof}
\emph{Step 1. The coefficient estimates.}
After a target translation and simultaneous rotations of source and
target, we may take \(\mathcal L=\R\) and keep \(F\) monic. Suppose that
\(\Delta(F,\R)\le L\), where \(L\) is the constant in the statement.
The file \href{https://github.com/zhangteng2000/ehp-lemniscate-verification-20261108/blob/ca8c451bfbc7312b9c57acc5e672c0b69cd601ae/verification/threshold9_contours.json}{\nolinkurl{threshold9_contours.json}} gives \(100T\) as an
integer for each degree, derivative order, and interval \([l,h]\).
With the bounds \eqref{eq:conditional-disk}, every choice satisfies
\(\rho<1\). The program \href{https://github.com/zhangteng2000/ehp-lemniscate-verification-20261108/blob/ca8c451bfbc7312b9c57acc5e672c0b69cd601ae/verification/verify_conditional_line.py}{\nolinkurl{verify_conditional_line.py}}
applies \eqref{eq:conditional-tests} and \eqref{eq:conditional-update}
to obtain bounds on \(|c_2|,\ldots,|c_{n-2}|\). The independent
implementation \href{https://github.com/zhangteng2000/ehp-lemniscate-verification-20261108/blob/ca8c451bfbc7312b9c57acc5e672c0b69cd601ae/verification/audit_conditional_line.py}{\nolinkurl{audit_conditional_line.py}} checks the
stored bounds for \(D\) and recomputes the coefficient bounds by
multiplying the power series in Lemma~\ref{lem:peel-recursion}.
There are \(5\,766\) interval tests over the six degrees.

For \(10\le n\le14\), take \(u,r\) from
Table~\ref{tab:conditional-gram}.
\par\nopagebreak[4]
\begin{table}[H]
\centering
\begin{tabular}{@{}lrrrrr@{}}
\toprule
\(n\)    & 10 & 11 & 12 & 13 & 14\\
\midrule
\(100u\)  & 66 & 58 & 52 & 45 & 39\\
\(1000r\) & 1445 & 1385 & 1340 & 1290 & 1260\\
\bottomrule
\end{tabular}
\par\smallskip
\caption{Cutoffs for \(|c_1|\) and radii for
Lemma~\ref{lem:collinear-gram}.}
\label{tab:conditional-gram}
\end{table}

If \(|c_1|\le u\),
Lemma~\ref{lem:companion-variance} gives \(|\mathcal V(F)|/2<24/25\).
This contradicts \(|\mathcal V(F)|\ge2\), since the critical
values are collinear and have diameter two. If \(|c_1|\ge u\),
the conditions on \(r\) in Lemma~\ref{lem:collinear-gram} hold,
\(m'-P_2>0\) for \(m'=n-1\), and the right-hand side of
\eqref{eq:collinear-gram} divided by \(m'(m'-P_2)\) is less than
\(24/25\). This violates \eqref{eq:collinear-gram}. Both comparisons use
rational enclosures. The two implementations evaluate the variance
by Newton power sums and by matrix multiplication, respectively.

\emph{Step 2. A test for nonreal roots in degree nine.}
The same coefficient calculation, now with \(L=41/50\), gives
\begin{equation}\label{eq:nonic-source-coefficients}
 |c_j|\le M_j\quad(2\le j\le7),\qquad
 (M_2,\ldots,M_7)=10^{-3}(103,14,84,311,521,785).
\end{equation}
For \(|c_1|\le12/25\), these bounds and
Lemma~\ref{lem:companion-variance} again give
\(|\mathcal V(F)|/2<47/50\). It remains to consider
\(|c_1|>12/25\). Write \(\gamma=\arg c_1\) and \(m=8\).
For fixed \(\theta,h\), let \(q(t)=\sum_{j=0}^m q_jt^j\) be the
real polynomial given on the real axis by
\[
 q(t)=\frac19\frac{d}{dt}\Im F(e^{i\theta}(t+ih)).
\]
Set \(B_m=1\), \(B_{m-1}=0\), and \(B_j=(j+1)M_{j+1}/9\) for
\(1\le j\le m-2\). Put \(V(H)=\sum_{j=1}^m B_jH^j\).
Suppose that \(0<X\le|\sin(9\theta)|\),
\(0\le Y\le|\sin(\theta+\gamma)|\), and \(|h|\le H\).
The contribution of \(c_1\) to \(q_0\) is
\(|c_1|\sin(\theta+\gamma)/9\). Expanding the remaining terms
in powers of \(t+ih\) and using \eqref{eq:nonic-source-coefficients}
gives
\[
 |q_0|\ge D_0:=\frac{4Y}{75}-V(H),\qquad
 |q_1|\le V'(H),\qquad 2|q_2|\le V''(H).
\]
At the other end of the polynomial, centering gives
\[
 q_m=\sin(9\theta),\qquad q_{m-1}=mh\cos(9\theta),\qquad
 q_{m-2}=-\binom m2 h^2\sin(9\theta)
                 +\tfrac79\Im(c_7e^{7i\theta}).
\]

Suppose \(D_0>0\) and all roots \(t_1,\ldots,t_m\) of \(q\) are
real, counted with multiplicity. All are nonzero. Applying the
second Newton identity to \(q\) and to its reversed polynomial
gives
\begin{equation}\label{eq:nonic-two-moments}
 \begin{aligned}
 \sum_i t_i^2
 &= (q_{m-1}/q_m)^2-2q_{m-2}/q_m
 \le T_0:=mH^2(m/X^2-1)+2B_{m-2}/X,\\
 \sum_i t_i^{-2}
 &= (q_1/q_0)^2-2q_2/q_0
 \le \frac{V'(H)^2}{D_0^2}+\frac{V''(H)}{D_0}.
 \end{aligned}
\end{equation}
For the first estimate, the terms in \(h^2\) combine to
\(mh^2(m/\sin^2(9\theta)-1)\).
Cauchy--Schwarz gives
\(m^2\le(\sum_i t_i^2)(\sum_i t_i^{-2})\).
It follows from \eqref{eq:nonic-two-moments} that the conditions
\begin{equation}\label{eq:nonic-strip}
 D_0>0,\qquad m^2D_0^2>T_0\bigl(V'(H)^2+D_0V''(H)\bigr)
\end{equation}
force a nonreal conjugate pair for every \(|h|\le H\). The
derivative slice therefore has at most \(m-2\) distinct real
zeros. By Rolle's theorem the original slice has at most \(n-2\).

\emph{Step 3. The area of the strips.}
The angular dependence in \eqref{eq:nonic-strip} enters only
through \(|\sin(9\theta)|\) and \(|\sin(\theta+\gamma)|\).
The change \((\theta,\gamma)\mapsto(\theta-\pi/9,\gamma+\pi/9)\)
and simultaneous reflection of both angles preserve the area
certified by this test. Thus it suffices to cover
\(0\le\gamma\le\pi/18\).

The program \href{https://github.com/zhangteng2000/ehp-lemniscate-verification-20261108/blob/ca8c451bfbc7312b9c57acc5e672c0b69cd601ae/verification/verify_reciprocal_strips.py}{\nolinkurl{verify_reciprocal_strips.py}} divides this
phase interval into six equal parts and divides \([0,\pi]\) into
\(360\) equal intervals in \(\theta\). On each closed rectangle it
uses lower bounds \(X,Y\) for the two absolute sines. It assigns
height zero or a height \(H=b/512\), \(b\) a positive integer,
for which \eqref{eq:nonic-strip} holds. In order of increasing
phase, the sums of the height numerators over the angular intervals
are
\[
 28463,\quad 28431,\quad 28346,\quad 28258,\quad 28144,\quad 28020.
\]
For each fixed phase the certified area in \((\theta,h)\) is
therefore at least \(2\pi\cdot28020/(360\cdot512)>19/20\).
An independent calculation in \href{https://github.com/zhangteng2000/ehp-lemniscate-verification-20261108/blob/ca8c451bfbc7312b9c57acc5e672c0b69cd601ae/verification/audit_low_degree_transfer.py}{\nolinkurl{audit_low_degree_transfer.py}}
covers the full phase range modulo \(\pi\) in \(36\) intervals,
with \(336\) angular intervals and heights in multiples of
\(1/384\). It uses a cubic lower bound for the sine and gives a
minimum height sum of \(18963\), hence an area greater than
\(23/25\). On every certified strip the integrand
\((n-N(\theta,h))/2\) in \eqref{eq:new-crofton} is at least one.
Either area bound contradicts \(\Delta(F,\R)\le41/50\).
\end{proof}

\subsection{The ninth exterior coefficient on a short arc}

Write \(x=1/n\), \(b=(1-x)/2\), and recall from
\eqref{eq:exterior-constant-normalizations} that
\[
 I_0=\frac{\Gamma(1+2b)}{\Gamma(1+b)^2},\qquad
 H^2=\frac{U_n}{2\pi}-1=\frac{n-1}{4}I_0.
\]
In degree nine we need a one-sided bound on \(h_9\). We first
estimate its nonlinear part on the full set of Laurent coefficients
allowed by the area theorem.

\begin{lem}\label{lem:nonic-area-remainder}
Let \(Z(t)=\sum_{j\ge1}j a_jt^{j+1}\), where
\(\sum_{j\ge1}j|a_j|^2\le1\). Then
\[
 \left|[t^9]\left(\sqrt{1-Z(t)}-1
                   -\tfrac12\log(1-Z(t))\right)\right|<\frac9{10}.
\]
\end{lem}

\begin{proof}
Denote the coefficient by \(R_9\), and put \(y_j=\sqrt j\,|a_j|\).
Since \(Z(t)=O(t^2)\), only \(Z^2/8\), \(5Z^3/48\), and
\(11Z^4/128\) contribute. Collecting the terms of degree nine gives
\begin{equation}\label{eq:nonic-area-expansion}
\begin{aligned}
 |R_9|\le{}&
 \tfrac{\sqrt6}{4}y_1y_6+\tfrac{\sqrt{10}}4y_2y_5
                         +\tfrac{\sqrt{12}}4y_3y_4\\
 &+\tfrac58\bigl(y_1^2y_4+\sqrt6\,y_1y_2y_3\bigr)
   +\tfrac{5\sqrt2}{24}y_2^3
   +\tfrac{11\sqrt2}{32}y_1^3y_2.
\end{aligned}
\end{equation}
Let \(s=y_1^2+y_2^2\). Young's inequality bounds the first
line of \eqref{eq:nonic-area-expansion} by
\((11/40)s+(23/40)(1-s)\). Indeed,
\(2\sqrt{(11/40)(23/40)}>\sqrt{10}/4\) for the first two
products, and \(\sqrt{12}/8<23/40\) for the third.
Also \(y_2^3\le y_2^2\) and \(y_1^3y_2\le y_1^2/2\), since
\(2y_1y_2\le s\le1\). Thus the last two terms are at most
\(3s/10\). Together these bounds account for at most \(23/40\).

For the remaining two terms, put \(u=y_1^2\). Cauchy--Schwarz
and \(y_3^2+y_4^2\le1-s\) give
\[
 \tfrac58\bigl(y_1^2y_4+\sqrt6\,y_1y_2y_3\bigr)
 \le\tfrac58\sqrt{u(6s-5u)(1-s)}
 \le\tfrac{3\sqrt5}{8}s\sqrt{1-s}
 \le\tfrac{\sqrt{15}}{12}.
\]
Here \(u(6s-5u)\le9s^2/5\), and \(s\sqrt{1-s}\) attains its
maximum \(2/(3\sqrt3)\) at \(s=2/3\). Since
\(\sqrt{15}/12<13/40\), \eqref{eq:nonic-area-expansion} gives
\(|R_9|<23/40+13/40=9/10\).
\end{proof}

Suppose now that \(p\) is centered and monic, with finite critical
values on \(\mathbb T\) of diameter \(2\delta\), where
\(0<\delta\le1/2\). By simultaneous rotations of the source and
target, arrange that a farthest pair is \(c\pm i\delta\), with
\(c=\sqrt{1-\delta^2}>0\). Write
\(p(z)=c+\delta Q(\delta^{-x}z)\). The critical values of \(p\)
lie on the shorter arc between the endpoints. In particular, each
such value \(v\) satisfies \(\Re v\ge c\) and
\(|v-c|^2=1-2c\Re v+c^2\le\delta^2\). Hence the critical values
\(q_i\) of \(Q\) satisfy \(|q_i|\le1\) and \(\Re q_i\ge0\).

\begin{lem}\label{lem:nonic-chord-sign}
In this normalization, if \(n=9\), then
\[
 \Re h_9(p)\ge bc-\tfrac9{10}\delta .
\]
\end{lem}

\begin{proof}
The exterior inverse of \(Q\) is univalent on \(|w|>1\), and its
Laurent coefficients satisfy \(\sum j|a_j|^2\le1\) by the area
theorem. Lemma~\ref{lem:inverse-log-identity} gives
\(h_9(Q)=\frac1{18}\sum_{i=1}^8 q_i+R_9(Q)\).
The real part of the sum is nonnegative, and
Lemma~\ref{lem:nonic-area-remainder} gives \(|R_9(Q)|<9/10\).
Now apply \(h_9(p)=bc+\delta h_9(Q)\) from
\eqref{eq:affine-inverse}.
\end{proof}

\subsection{The circle comparison in degrees nine through twenty-seven}
\label{sec:collinear-circle}

\begin{prop}\label{prop:circle-nine-to-twenty-seven}
Let \(p\) be a monic polynomial of degree \(9\le n\le27\), and
put \(p_0(z)=z^n-1\). Then
\[
 \ell(\partial E_1(p))\le\ell(\partial E_1(p_0)).
\]
\end{prop}

\begin{proof}
By Lemma~\ref{lem:extremizers}, there is a maximizer whose finite
critical values lie on \(\mathbb T\). Let \(p\) be any such maximizer,
centered by a translation of the source, and let \(2\delta\) be the
diameter of its critical values.
We shall exclude \(\delta>0\). Use \(L,K,\tau_0\) from
Table~\ref{tab:circle-nine-twenty-seven}.
\par\nopagebreak[4]
\begin{table}[H]
\centering
\begin{tabular}{@{}lccc@{}}
\toprule
\(n\) & \(L\) & \(K\) & \(\tau_0\)\\
\midrule
\(9\) & \(41/50\) & \(9/10\) & \(1/320\)\\
\(10\)--\(14\) & \(7/8\) & \(H\) & \(1/(8n^2)\)\\
\(15\)--\(27\) & \(49/60\) & \(H\) & \(1/(128n^2)\)\\
\bottomrule
\end{tabular}
\par\smallskip
\caption{Constants for the circle comparison, with
\(H=(U_n/(2\pi)-1)^{1/2}\).}
\label{tab:circle-nine-twenty-seven}
\end{table}

The line bounds \(L\)
are supplied by Propositions~\ref{prop:conditional-line}
and~\ref{prop:collinear-finite}. In degree nine the value \(K=9/10\)
comes from Lemma~\ref{lem:nonic-area-remainder}.

\emph{Step 1. Critical values on a short arc.}
Put \(\delta_0=2\tau_0/(1+\tau_0^2)\). When
\(0<\delta\le\delta_0\), the normalization
\eqref{eq:short-arc-normalization} gives a parameter
\(0<\tau\le\tau_0\) and a centered monic polynomial \(F\) whose
critical values lie in \([-1,1]\) and include both endpoints.
Thus \(\Delta(F,\R)>L\).

For \(9\le n\le14\) take \(R=3\) and \(b_0=2/5\) in
Lemma~\ref{lem:curvature-general}. The restrictions
\(\tau<b_0<1-1/R\) and \(\tau R\le1/2\) hold throughout the
stated ranges. In degree nine, Lemma~\ref{lem:nonic-area-remainder}
and \eqref{eq:signed-short-coefficient} allow the upper bound
\(9/20+4\tau/9\) for the signed coefficient \(h\). In the
other degrees take \(h\le H\).

With these choices it suffices to bound the right-hand side of
\eqref{eq:curvature-general} at \(\tau=\tau_0\). Indeed, after
the indicated substitution for \(h\), the elementary terms are
nondecreasing in \(\tau\). In the two series terms, incorporating
the outside factor \(s^{1+x}\) leaves powers \(s^{1-v}\) with
\(v>1\), namely \(v=(j+k-1)/n\) in the row sums and \(v=q-x\)
in the linear sum. These powers increase as
\(s=R/\sqrt{1+\tau^2R^2}\) decreases. The two rational computations
in \href{https://github.com/zhangteng2000/ehp-lemniscate-verification-20261108/blob/ca8c451bfbc7312b9c57acc5e672c0b69cd601ae/verification/single_coefficient_transfer.py}{\nolinkurl{single_coefficient_transfer.py}} and
\href{https://github.com/zhangteng2000/ehp-lemniscate-verification-20261108/blob/ca8c451bfbc7312b9c57acc5e672c0b69cd601ae/verification/independent_single_coefficient_transfer.py}{\nolinkurl{independent_single_coefficient_transfer.py}} bound the
error by less than \(4/5\) for \(n=9\) and less than \(5/8\)
for \(10\le n\le14\), including the tails of both series.
For \(15\le n\le27\), \eqref{eq:curvature256} bounds the error
by \(1/2+8/128^{5/6}<4/5\).
In every case the error is smaller than \(L\).
The curvature comparison therefore gives
\(\ell(\partial E_1(p))<\ell(\partial E_1(p_0))\), contrary to
maximality.

\emph{Step 2. The remaining diameters.}
Choose two critical values \(v_1,v_2\) at distance \(2\delta\).
For every nonexceptional \(\xi\in\mathbb T\), their images under
\(M_\xi\) are separated by
\(2\delta/(|1-\xi v_1|\,|1-\xi v_2|)\), and all transformed
critical values lie on \(\mathcal L_\xi\). Scaling the appropriate
line bound in Table~\ref{tab:circle-nine-twenty-seven} and
averaging gives
\begin{equation}\label{eq:circle-nine-mean-line}
 \big\langle\Delta(F_\xi,\mathcal L_\xi)\big\rangle
 \ge L\delta^{1/n}
       \big\langle |1-\xi v_1|^{-1/n}|1-\xi v_2|^{-1/n}\big\rangle
 \ge L\delta^{1/n}.
\end{equation}
The last step is Jensen's inequality, since
\(\langle\log|1-\xi v_i|\rangle=0\).

For \(0<\delta\le1/2\), rotate the source and target, preserving
monicity, so that a farthest pair is \(c\pm i\delta\), where
\(c=\sqrt{1-\delta^2}\). Put \(Q(z)=\delta^{-1}(p(\delta^{1/n}z)-c)\).
All critical values of \(p\) lie on the shorter arc between the
endpoints, hence in \(\overline{D(c,\delta)}\). The critical values
of \(Q\) therefore lie in the closed unit disk.
Equations~\eqref{eq:exterior-parseval} and~\eqref{eq:universal-length-bound}
give \(|h_n(Q)|\le H\), and \eqref{eq:affine-inverse} gives
\(h_n(p)=bc+\delta h_n(Q)\).
In degree nine,
Lemma~\ref{lem:nonic-chord-sign} gives the sharper lower bound
for \(\Re h_9(p)\). In all cases,
\(|h_n(p)|\ge[b\sqrt{1-\delta^2}-K\delta]_+\), where
\([y]_+=\max(y,0)\). For a closed interval \([l,h]\) contained
in \([\delta_0,1/2]\), put
\[
 e(l,h)=2\pi(I_0-1)[b\sqrt{1-h^2}-Kh]_+^2.
\]
On intervals contained in \([1/2,1]\), use \(e(l,h)=0\).
Then \eqref{eq:full-exterior-coefficient} gives
\(\sum_{k\ge1}\mathcal E_k\ge e(l,h)\) on each interval.

By \eqref{eq:new-mass} and \eqref{eq:circle-nine-mean-line},
it suffices to verify
\begin{equation}\label{eq:circle-nine-interval-test}
 U_n-\ell(\partial E_1(p_0))<e(l,h)+Ll^{1/n}
\end{equation}
on a partition of \([\delta_0,1]\). The two programs in Step~1
verify \eqref{eq:circle-nine-interval-test} for all nineteen
degrees. Bisection uses at most \(19\) closed intervals per degree,
and the independent trisection uses at most \(26\). Each program
checks that consecutive intervals share an endpoint and that the
first and last endpoints are \(\delta_0\) and \(1\).

For the independent evaluation of \(I_0=\sum_{j\ge0}\binom bj^2\),
the partial sum includes \(0\le j\le s\) with \(s=240\).
For \(j>s\) the recurrence for the binomial coefficients gives
\(\big|\binom bj\big|\le(s/j)\big|\binom bs\big|\), so the omitted
tail is at most \(s\binom bs^2\). The beta values are enclosed
using \eqref{eq:split-beta-series}. Thus
\eqref{eq:circle-nine-interval-test} contradicts maximality for
the remaining positive diameters as well.

Thus every centered maximizer whose finite critical values lie on
\(\mathbb T\) has a single critical value \(v\in\mathbb T\).
By Lemma~\ref{lem:critical-compactness}, it has the form \(z^n+v\),
so its lemniscate is a rotation of \(\partial E_1(p_0)\).
Since such a maximizer exists, this proves the length comparison.
\end{proof}
\section{The comparison in degree eight}\label{sec:octic}

On the line \(z=e^{i\theta}(t+ih)\), intersections with
\(F^{-1}(\R)\) are the real zeros of
\(t\mapsto\Im F(e^{i\theta}(t+ih))\). If a derivative of order \(d\)
has degree \(n-d\) and has a nonreal zero, it has at most
\(n-d-2\) distinct real zeros. Rolle's theorem then bounds the
number of distinct real zeros of the original slice by \(n-2\).
By \eqref{eq:new-crofton}, a set of such lines of
\((\theta,h)\)-measure \(A\) contributes at least \(A\) to
\(\Delta(F,\R)\). We shall use this observation in all three degrees.

In degree eight we use the first derivative. Its constant term at
\(h=0\) contains the linear coefficient of \(F\). To exploit this
coefficient we apply the arithmetic--geometric mean inequality to
the reciprocals of the roots.

\begin{lem}\label{lem:octic-product}
Let \(q(t)=\sum_{j=0}^m q_jt^j\) be a real polynomial of degree
\(m\ge2\). If all its zeros are real, then
\begin{equation}\label{eq:reciprocal-root-product}
 m|q_m|^{2/m}|q_0|^{2-2/m}\le q_1^2-2q_0q_2 .
\end{equation}
\end{lem}

\begin{proof}
If \(q_0=0\), the right-hand side is \(q_1^2\). Otherwise let
\(r_1,\ldots,r_m\) be the zeros, counted with multiplicity.
Comparing the first two coefficients in
\(q(t)/q_0=\prod_{j=1}^m(1-t/r_j)\) gives
\(\sum_j r_j^{-2}=(q_1/q_0)^2-2q_2/q_0\).
The geometric mean of the \(r_j^{-2}\) is \(|q_m/q_0|^{2/m}\).
Multiplication by \(q_0^2\) now gives \eqref{eq:reciprocal-root-product}.
\end{proof}

\begin{prop}\label{prop:octic-line}
Let \(F\) be a centered monic polynomial of degree eight.
If its finite critical values lie on a line \(\mathcal L\) and have
diameter two, then \(\Delta(F,\mathcal L)>39/50\).
\end{prop}

\begin{proof}
\emph{Step 1. The linear coefficient.}
A rigid motion in the target and a rotation in the source take
\(\mathcal L\) to \(\R\) and preserve monicity and the deficit.
Write \(F(z)=\sum_{j=0}^8c_jz^j\) in these coordinates, and suppose
that \(\Delta(F,\R)\le39/50\). Applying
\eqref{eq:conditional-disk} and \eqref{eq:conditional-tests}
at derivative orders \(6,5,4,3,2\) gives \(|c_j|\le M_j\), where
\begin{equation}\label{eq:octic-caps}
 (M_2,M_3,M_4,M_5,M_6,M_7,M_8)
   =10^{-3}(197,62,319,512,774,0,1000).
\end{equation}
The contour and energy data are given in
\href{https://github.com/zhangteng2000/ehp-lemniscate-verification-20261108/blob/ca8c451bfbc7312b9c57acc5e672c0b69cd601ae/verification/octic_contours.json}{\nolinkurl{octic_contours.json}} and \href{https://github.com/zhangteng2000/ehp-lemniscate-verification-20261108/blob/ca8c451bfbc7312b9c57acc5e672c0b69cd601ae/verification/octic_energy.json}{\nolinkurl{octic_energy.json}}.
There is one normalized radius for the quadratic derivative and
there are 128 radius intervals for each of the other four derivatives.

If \(|c_1|\le13/25\), Lemma~\ref{lem:companion-variance} and
\eqref{eq:octic-caps} give \(|\mathcal V(F)|/2<97/100\).
If \(\bar v\) is the mean of all seven critical values, counted
with multiplicity, the two extreme values are distance two apart
and their squared distances from \(\bar v\) sum to at least \(2\).
Thus \(|\mathcal V(F)|\ge2\), and \(|c_1|>13/25\).

\emph{Step 2. Counting intersections.}
Put \(\gamma=\arg c_1\) and
\(q(t)=\tfrac18\frac{d}{dt}\Im F(e^{i\theta}(t+ih))\).
Its leading coefficient is \(\sin(8\theta)\), so \(q\) has degree
\(m=7\) except in finitely many directions. Set
\(V(H)=\sum_{j=1}^m (j+1)M_{j+1}H^j/8\).
On an angular interval choose \(0<X\le|\sin(8\theta)|\) and
\(0\le Y\le|\sin(\theta+\gamma)|\). For \(|h|\le H\), expansion
in powers of \(t+ih\) gives
\[
 |q_0|\ge D:=13Y/200-V(H),\qquad
 |q_1|\le V'(H),\qquad 2|q_2|\le V''(H).
\]
If \(q\) had only real zeros and \(D>0\), then
\eqref{eq:reciprocal-root-product}, divided by
\(|q_0|^{2-2/m}\), would give
\(mX^{2/m}\le V'(H)^2D^{-2+2/m}+V''(H)D^{-1+2/m}\).
A nonreal pair is therefore guaranteed throughout the strip whenever
\begin{equation}\label{eq:octic-strip}
 D>0,\qquad
 m^mX^2D^{2m-2}>\bigl(V'(H)^2+D V''(H)\bigr)^m .
\end{equation}

The two absolute sines entering \eqref{eq:octic-strip} are unchanged
under \((\theta,\gamma)\mapsto(\theta-\pi/8,\gamma+\pi/8)\) and
under simultaneous reflection. Since the angular integral has
period \(\pi\) in \(\theta\), it suffices to take
\(0\le\gamma\le\pi/16\). Divide this phase interval into eight equal
parts and \(0\le\theta\le\pi\) into 512 equal parts. In each rectangle
the bounds \(X,Y\) hold throughout, and the accepted height is
\(H=k/640\), with \(k=0\) if no strip is used.
Writing \(k_j\) for the chosen integer on the \(j\)-th angular
interval, the strip certificate gives
\(\min_{\text{phase intervals}}\sum_{j=1}^{512}k_j=42820\)
(\href{https://github.com/zhangteng2000/ehp-lemniscate-verification-20261108/blob/ca8c451bfbc7312b9c57acc5e672c0b69cd601ae/verification/octic_root_product_strips.json}{\nolinkurl{octic_root_product_strips.json}}). Counting both signs of
\(h\) gives
\[
 \Delta(F,\R)\ge\frac{2\pi\cdot42820}{512\cdot640}
      >\frac{821}{1000}>\frac{39}{50}.
\]
This contradicts the assumed deficit. \href{https://github.com/zhangteng2000/ehp-lemniscate-verification-20261108/blob/ca8c451bfbc7312b9c57acc5e672c0b69cd601ae/verification/octic_line_independent.json}{An independent check} over the
full phase interval \(0\le\gamma\le\pi\), with heights in multiples
of \(1/480\), gives the lower bound \(819/1000\).
\end{proof}

By Propositions~\ref{prop:octic-line} and~\ref{prop:small-circle-transfer},
a degree-eight maximizer whose critical values lie on the unit circle
has just one critical value. Lemmas~\ref{lem:extremizers}
and~\ref{lem:critical-compactness} then give the length comparison
with \(z^8-1\). Section~\ref{sec:equality} treats arbitrary maximizers.

\section{The comparison in degrees six and seven}\label{sec:last-degrees}

Set \(L_6=73/100\) and \(L_7=367/500\), as in
Proposition~\ref{prop:small-circle-transfer}. The derivative of
order \(n-3\) is cubic and depends, apart from its leading term,
only on the coefficients of \(z^{n-2}\) and \(z^{n-3}\). It gives
a bound on their common scale for each fixed ratio of their moduli.
We keep this information when we return to the first two derivatives.

\subsection{Bounds from the cubic derivative}

Let \(n\in\{6,7\}\), and write
\(F(z)=z^n+\sum_{j=0}^{n-2}c_jz^j\). Assume that its finite
critical values lie on \(\mathcal L\) and have diameter two. Put
\(f=3|c_{n-2}|/\binom n2\) and \(g=|c_{n-3}|/\binom n3\).
With \(s\ge0\), we write \((f,g)\) in one of the two forms
\begin{equation}\label{eq:sextic-septic-charts}
\begin{array}{c@{\qquad}c@{\qquad}c@{\qquad}l}
 \mathrm A & f=s^2  & g=ys^3 & 0\le y\le 5/16,\\
 \mathrm B & f=ys^2 & g=s^3  & 0\le y\le 9/4 .
\end{array}
\end{equation}
For \(f,g>0\), A applies if \(g/f^{3/2}\le5/16\). Otherwise
\(f/g^{2/3}<(16/5)^{2/3}<9/4\), so B applies.
The axes are included, and at \(f=g=0\) we take \(s=y=0\) in A.
In each row of \eqref{eq:sextic-septic-charts}, divide the
\(y\)-interval into 32 equal closed subintervals, numbered from
\(0\) to \(31\) in increasing order.

\Needspace{24\baselineskip}
\begin{lem}\label{lem:last-cubic-profile}
If \(y\) belongs to subinterval \(i\) in parametrization A or B,
then \(\Delta(F,\mathcal L)\ge s q_i\), where \(q_i\) is taken from the
corresponding column of Table~\ref{tab:last-profile}.
\end{lem}
\par\nopagebreak[4]
\begin{table}[H]
\centering\small
\begin{tabular}{r@{\quad}rrrr@{\qquad}r@{\quad}rrrr}
\toprule
\(i\)&6A&6B&7A&7B&\(i\)&6A&6B&7A&7B\\
\midrule
 0&181&201&101&201&16&200&263&213&285\\
 1&182&201&107&201&17&201&267&214&288\\
 2&184&201&122&201&18&202&271&215&292\\
 3&185&203&134&201&19&203&275&216&296\\
 4&187&207&145&201&20&204&279&217&300\\
 5&188&211&154&201&21&205&283&218&303\\
 6&189&216&162&201&22&206&286&219&286\\
 7&191&221&169&201&23&207&290&220&286\\
 8&192&226&176&212&24&208&275&221&286\\
 9&193&231&182&224&25&208&275&222&290\\
10&194&235&188&235&26&209&275&223&262\\
11&195&240&194&245&27&210&275&224&262\\
12&196&244&199&255&28&211&275&225&286\\
13&198&249&204&265&29&212&275&221&300\\
14&198&254&208&274&30&212&314&221&300\\
15&199&259&212&281&31&213&318&221&300\\
\bottomrule
\end{tabular}
\par\smallskip
\caption{The integers \(100q_i\) in Lemma~\ref{lem:last-cubic-profile}.
Each column gives the degree and the parametrization in
\eqref{eq:sextic-septic-charts}.}
\label{tab:last-profile}
\end{table}

\begin{proof}
For \(s=0\) the assertion follows from \eqref{eq:new-crofton}.
For \(s>0\) set \(F_s(z)=s^{-n}F(sz)\) and
\(\mathcal L_s=s^{-n}\mathcal L\). Then
\(\Delta(F,\mathcal L)=s\,\Delta(F_s,\mathcal L_s)\), and the coefficient pair
for \(F_s\) is \((1,y)\) in A and \((y,1)\) in B.
By Lemma~\ref{lem:cubic-slice}, it suffices to find an area of at
least \(q_i\) on which the expression \(D\) in \eqref{eq:new-D}
is positive. Indeed, the discriminant of the cubic slice is
\(-27D\).

Use the phase coordinates \(\alpha,\beta\) introduced before
\eqref{eq:cubic-phase-expansion}.
Adding \(\pi\) to \(\alpha\) leaves \(D\) unchanged. Adding \(\pi\)
to \(\beta\) replaces \(h\) by \(-h\). Conjugation reverses both
phases and preserves area. We may therefore restrict
\((\alpha,\beta)\) to \([0,\pi/2]\times[0,\pi]\).
For each fixed choice of \(y,\alpha,\beta\), we count rectangles in
\([0,\pi]\times[-R,R]\), with \(R=5/4\) in A and \(R=7/4\) in B.

Here is the criterion used in that count. Given a box in
\((y,\alpha,\beta)\) and a rectangle in \((\theta,h)\), let \(w_j\)
be the five coordinate widths and let \(U_j\) be the corresponding
bounds from \eqref{eq:cubic-hessian}. Multilinear interpolation gives
\begin{equation}\label{eq:small-cubic-cell}
 \min_{\text{vertices}}D>\frac18\sum_j U_jw_j^2
 \quad\Longrightarrow\quad D>0\ \text{throughout the cell}.
\end{equation}
The coefficient of \(f^2r^2\) in \(U_\theta\) is bounded above by
\(64\) for \(n=6\) and by \(96\) for \(n=7\). In \(U_h\) the
negative term \(-120v^4\) uses \(v=\inf|h|\) on the height interval.
Thus a single application of \eqref{eq:small-cubic-cell} certifies
the rectangle for every parameter in the box.

Initially the three parameter intervals have \(8,16,32\) subdivisions,
respectively. Bisecting each coordinate of an unresolved box gives
eight children. At depths \(0,1,2,3\), the numbers \((T,J)\) of
angular and height intervals are, respectively,
\((256,128),(512,256),(1024,256),(1024,512)\).
If \(b\) rectangles pass, their interiors are disjoint and their
total area is \(2\pi Rb/(TJ)\). Using \(\pi>333/106\), we accept
the box once \(2R(333/106)b/(TJ)\ge q_i\) for every \(y\)-subinterval
whose interior meets the interior of its \(y\)-projection.
The tests hold on closed boxes and therefore include shared faces.

Table~\ref{tab:last-leaves} gives the complete subdivision.
\par\nopagebreak[4]
\begin{table}[H]
\centering\small
\begin{tabular}{lrrrr}
\toprule
 &Depth 0&Depth 1&Depth 2&Depth 3\\
\midrule
6A&1378&21139&4840&0\\
6B&1087&21885&17112&3072\\
7A&1155&20823&21640&0\\
7B&1024&20960&28359&4552\\
\bottomrule
\end{tabular}
\par\smallskip
\caption{Number of terminal parameter boxes at each depth.}
\label{tab:last-leaves}
\end{table}

The supplement's \href{https://github.com/zhangteng2000/ehp-lemniscate-verification-20261108/blob/ca8c451bfbc7312b9c57acc5e672c0b69cd601ae/verification/verify_terminal_cubic_tree.py}{\nolinkurl{verify_terminal_cubic_tree.py}} checks
that every nonterminal box has all eight children and recomputes
the accepted area of every terminal box. One evaluation uses
outward rounding on the \(2^{-24}\) lattice. The second uses the
factored expression \eqref{eq:new-D} and sine enclosures of width
at most \(2^{-46}\). The combined error in a vertex value and its
interpolation allowance is less than \(2^{-28}\). The additional
margin \(2^{-25}\) thus gives the strict inequality in
\eqref{eq:small-cubic-cell} on every accepted cell. These checks
give \(\Delta(F_s,\mathcal L_s)\ge q_i\), as required.
\end{proof}

\subsection{The first two derivatives of a real slice}

\begin{prop}\label{prop:last-line}
Let \(F\) be a centered monic polynomial of degree \(n\in\{6,7\}\).
If its finite critical values lie on a line \(\mathcal L\) and have
diameter two, then \(\Delta(F,\mathcal L)>L_n\).
\end{prop}

\begin{proof}
\emph{Step 1. Bounds on the coefficients.}
Take \(\mathcal L=\R\) as in the proof of Proposition~\ref{prop:octic-line}
and suppose that \(\Delta(F,\R)\le L_n\).
Fix a row of \eqref{eq:sextic-septic-charts} and a subinterval
\([y_-,y_+]\) from Table~\ref{tab:last-profile}.
Lemma~\ref{lem:last-cubic-profile} gives \(0\le s\le S:=L_n/q_i\).
On a smaller interval \(s_-\le s\le s_+\), the two monomials in
\eqref{eq:sextic-septic-charts} give upper and lower bounds for
both coefficient moduli. For example, in B they give
\[
 \frac{\binom n2}{3}y_-s_-^2\le |c_{n-2}|
       \le\frac{\binom n2}{3}y_+s_+^2,\qquad
 \binom n3 s_-^3\le |c_{n-3}|\le\binom n3 s_+^3 .
\]

Starting with these upper bounds, apply \eqref{eq:conditional-disk}
and \eqref{eq:conditional-tests} twice in the order \(k=n-2,\ldots,2\).
For \(m=n-k\), let \(A\) be the upper bound for \(U(G)/(2c_*)\) on a
normalized radius interval \([a,b]\), where \(G=m!F^{(k)}/n!\).
On this interval the resulting bound is
\begin{equation}\label{eq:small-degree-descent}
 |c_{n-j}|\le
 \frac{\binom nj}{\binom mj}\,e_j(a,b)A^j,\qquad 2\le j\le m .
\end{equation}
If a known lower bound exceeds the right-hand side, discard that
radius interval. At radius zero \(G(z)=z^m\), so this case is
incompatible with a positive lower bound for any of the same
coefficients. At each update take the maximum over the remaining
possibilities and retain the smaller of this bound and the old one.
Write \(M_j\) for the resulting upper bounds, with \(M_{n-1}=0\)
and \(M_n=1\).

When the first derivative is used below, \href{https://github.com/zhangteng2000/ehp-lemniscate-verification-20261108/blob/ca8c451bfbc7312b9c57acc5e672c0b69cd601ae/verification/verify_small_scalar_primary.py}{the supplement} supplies
a rational \(u>0\) for which Lemma~\ref{lem:companion-variance}
gives \(|\mathcal V(F)|<2\) if \(|c_1|\le u\).
The critical values have diameter two on \(\R\), so
\(|\mathcal V(F)|\ge2\). We can then use \(|c_1|\ge u\).
The second derivative imposes no condition on \(c_1\).

\emph{Step 2. Two moments of the roots.}
For \(d=1,2\) put \(m=n-d\), \(C=\binom nd\), and write
\[
 Q(t)=\frac{1}{d!\,C}\left(\frac{\mathrm d}{\mathrm dt}\right)^d
          \Im F\bigl(e^{i\theta}(t+ih)\bigr)
     =\sum_{j=0}^m Q_jt^j .
\]
Put \(a=\sin(n\theta)\), \(B=\binom{n-2}d/C\), and let
\(\phi=\arg c_{n-2}\), with any choice if \(c_{n-2}=0\).
The leading coefficients are \(Q_m=a\), \(Q_{m-1}=mh\cos(n\theta)\),
and
\(Q_{m-2}=-\binom m2h^2a+B|c_{n-2}|\sin((n-2)\theta+\phi)\).
Fix a direction with \(a\ne0\) and put \(x=|a|\).
If \(r_1,\ldots,r_m\) are the roots of \(Q\), counted with
multiplicity, Newton's identity gives
\begin{equation}\label{eq:small-forward-moment}
 x^2\sum_j r_j^2
   =mh^2(m-x^2)-2B a|c_{n-2}|\sin((n-2)\theta+\phi).
\end{equation}

For \(|h|\le H\) put
\(V(H)=C^{-1}\sum_{j=1}^m\binom{j+d}{d}M_{j+d}H^j\),
\(W=V'(H)\), and \(Z=V''(H)\).
If \(|c_d|\ge v\) and \(Y\le|\sin(d\theta+\arg c_d)|\), with
\(Y\ge0\), put \(D=vY/C-V(H)\). When \(c_d=0\) we take \(v=0\)
and any phase. Expansion in powers of \(t+ih\) gives
\(|Q_0|\ge D\), \(|Q_1|\le W\), and \(2|Q_2|\le Z\).
Choose \(0<X\le x\), and put
\(T(X)=mH^2(m-X^2)+2BM_{n-2}X\).
When \(D>0\), any one of the following inequalities forces \(Q\)
to have a nonreal pair.
\begin{equation}\label{eq:last-ordinary}
\begin{aligned}
 m^m X^2D^{2m-2}&>(W^2+DZ)^m,\\
 m^m X^{2m-2}D^2&>T(X)^m,\\
 m^2 X^2D^2&>T(X)(W^2+DZ).
\end{aligned}
\end{equation}

To see this, suppose all the roots are real. They are nonzero since
\(D>0\). The first test is Lemma~\ref{lem:octic-product} applied
as in \eqref{eq:octic-strip}. For the second, the geometric mean of
the \(r_j^2\) is \(|Q_0/a|^{2/m}\), while
\eqref{eq:small-forward-moment} gives \(x^2\sum_jr_j^2\le T(x)\).
Thus \(mD^{2/m}\le T(x)/x^{2-2/m}\). Both
\(T(x)/x^{2-2/m}\) and \(T(x)/x^2\) are nonincreasing on
\(0<x\le1\), since \(T(x)/x=mH^2(m/x-x)+2BM_{n-2}\)
is nonnegative and nonincreasing and \(m\ge2\).
This justifies substituting \(X\) in the second test.
For the third, Newton's identity gives
\(\sum_jr_j^{-2}\le W^2/D^2+Z/D\). Now use
\(m^2\le(\sum_jr_j^2)(\sum_jr_j^{-2})\) and
\(\sum_jr_j^2\le T(X)/X^2\).

\emph{Step 3. Retaining the sign of the \(c_{n-2}\) term.}
Suppose \(b_0\le|c_{n-2}|\le b_1\) and \(0<X\le x\le X_1\).
Let \(j_0\) be a lower bound for
\(\operatorname{sgn}(a)\sin((n-2)\theta+\phi)\) on the angular
interval under consideration. Define
\[
 A=\begin{cases}
       -2Bb_0Xj_0,&j_0\ge0,\\
       -2Bb_1X_1j_0,&j_0<0,
    \end{cases}
 \qquad T=mH^2(m-X^2)+A,\qquad T_c=m(m-1)H^2+A .
\]
The last term of \eqref{eq:small-forward-moment} is at most \(A\).
If all roots are real and \(\bar r=m^{-1}\sum_j r_j\), then
\begin{equation}\label{eq:small-signed-moments}
 x^2\sum_jr_j^2\le T,\qquad
 x^2\sum_j(r_j-\bar r)^2\le T_c .
\end{equation}
For the second inequality subtract
\(x^2m\bar r^2=mh^2\cos^2(n\theta)\) in
\eqref{eq:small-forward-moment}. In particular \(T_c<0\) rules out
real-rootedness even when \(D\le0\). If \(T_c\ge0\), then \(T\ge0\),
and the first inequality in \eqref{eq:small-signed-moments} gives
the forward product test \(m^mX^{2m-2}D^2>T^m\) whenever \(D>0\).

In degree seven we also use the two heights \(h\) and \(-h\)
together. Write \(Q_0(h)\) for the constant term at height \(h\)
and let \(V_{\rm even}\) be the sum of the even-degree terms of
\(V\). The odd powers of \(h\) cancel in
\((Q_0(h)+Q_0(-h))/2\), so
\[
 \max\{|Q_0(h)|,|Q_0(-h)|\}\ge vY/C-V_{\rm even}(H)
 \qquad (|h|\le H).
\]
The bounds in \eqref{eq:small-signed-moments} hold at both heights.
If both signs pass for \(0<h\le H_2\) and at least one passes for
\(H_2<h\le H_1\), the certified vertical measure is at least
\(2H_2+(H_1-H_2)=H_1+H_2\). The successful sign may depend on \(h\).

\emph{Step 4. Integrating the certified strips.}
On the intervals where \eqref{eq:last-ordinary} suffices, take
\(d=1\) and \(v=u\). As in degree eight, the phase of \(c_1\)
can be restricted to \([0,\pi/(2n)]\). The primary computation
uses 432 angular intervals, eight phase intervals and heights
\(k/512\). Bounds for \(X,Y\) hold on each closed angular
rectangle. A test at \(H\) applies to every \(|h|\le H\).

For the tests of Step~3, let \(\gamma=\arg c_d\).
The change \((\theta,\phi,\gamma)\mapsto
(\theta+\pi/n,\phi+2\pi/n,\gamma-d\pi/n)\) preserves both
\(a\sin((n-2)\theta+\phi)\) and \(|\sin(d\theta+\gamma)|\).
Simultaneous reflection also preserves these quantities.
It therefore suffices to take \(\phi\in[0,\pi/n]\) and
\(\gamma\in[0,\pi]\). The intervals requiring this test are
\(24\)--\(29\) in 6B, \(29\)--\(31\) in 7A, and \(22\)--\(31\)
in 7B. Ordinarily we subdivide \(s/S\) at
\(0,4/5,9/10,19/20,1\) and bisect \([0,M_2]\) for \(|c_2|\).
In 7B, intervals \(26,27\), the subdivision points for \(s/S\)
are \(0,80,88,92,95,97,99,100\), divided by \(100\), and
\([0,M_2]\) is divided into quarters. In these two subintervals,
parametrization B in \eqref{eq:sextic-septic-charts}, evaluated at
the lower endpoints of the \(s\)- and \(y\)-intervals, gives lower
bounds for \(|c_{n-2}|\) and \(|c_{n-3}|\). We discard a radius
interval if either lower bound exceeds the corresponding upper
bound in \eqref{eq:small-degree-descent}. After subdividing the
interval for \(|c_2|\), we repeat the coefficient estimates using
both endpoints of the chosen subinterval.
On each remaining box we use either \(d=1\) after the variance
test or \(d=2\) with \(v\) the lower endpoint for \(|c_2|\).

The basic mesh for Step~3 uses \(420,8,96\) intervals in
\(\theta,\phi,\gamma\), respectively, and heights \(k/400\).
In 7B, interval \(28\) uses \(840,12,144\) intervals and heights
\(k/600\). Intervals \(26,27\) use \(980,12,196\) intervals and
heights \(k/700\). For subintervals \(26\) and \(27\) of case \(7\mathrm B\),
Table~\ref{tab:last-profile} gives \(q_i=131/50\).
After returning to the original scale, the certified lower bounds
are strictly greater than \(7483/10000\) and \(7511/10000\),
respectively, both exceeding \(L_7=367/500\).

Taking the minimum over all subintervals and all phase boxes gives
strict lower bounds of \(7309/10000\), \(7302/10000\), \(7377/10000\)
and \(7352/10000\) in 6A, 6B, 7A and 7B, respectively. The first
two exceed \(L_6\), and the last two exceed \(L_7\).
The primary verification is \href{https://github.com/zhangteng2000/ehp-lemniscate-verification-20261108/blob/ca8c451bfbc7312b9c57acc5e672c0b69cd601ae/verification/verify_small_scalar_primary.py}{\nolinkurl{verify_small_scalar_primary.py}}.
The independent program \href{https://github.com/zhangteng2000/ehp-lemniscate-verification-20261108/blob/ca8c451bfbc7312b9c57acc5e672c0b69cd601ae/verification/audit_small_scalar_profiles.py}{\nolinkurl{audit_small_scalar_profiles.py}}
reconstructs the coefficient bounds and repeats the angular integration.
By \eqref{eq:new-crofton}, the certified area is a lower bound
for \(\Delta(F,\R)\), contrary to the assumption in Step~1.
\end{proof}

Propositions~\ref{prop:last-line} and~\ref{prop:small-circle-transfer}
now give the circle comparison in degrees six and seven.
As in degree eight, the passage from a selected maximizer to all
equality cases is made in Section~\ref{sec:equality}.

\section{Inverse images of lines in degree five}\label{sec:new-five-line}

\begin{prop}\label{prop:new-five-line}
Let \(F\) be a monic quintic whose finite critical values lie on a
line \(\mathcal L\) and have diameter \(d>0\). Then
\begin{equation}\label{eq:new-line-loss}
 \Delta(F,\mathcal L)>\frac{17}{16}(d/2)^{1/5}.
\end{equation}
\end{prop}

Two estimates enter the proof. The first bounds the deficit on any
target line in terms of the cubic and quadratic coefficients of the
centered quintic. The second uses collinearity to bound the same
coefficient expression from below in terms of \(d\).

We recall the scaling rule. If \(F\) is monic of degree \(n\),
\(P(z)=\lambda^{-n}(F(z_*+\lambda z)-w_*)\), and
\(\mathcal L'=\lambda^{-n}(\mathcal L-w_*)\), where \(\lambda\ne0\), then
\begin{equation}\label{eq:low-degree-line-scaling}
 \Delta(P,\mathcal L')=|\lambda|^{-1}\Delta(F,\mathcal L).
\end{equation}
Indeed, a target cutoff of size \(R^n\) for \(P\) corresponds to
one of size \((|\lambda|R)^n\) for \(F\), and the source lengths
differ by the factor \(|\lambda|\).

\begin{lem}\label{lem:new-cubic}
Let \(F(z)=z^5+Az^3+Bz^2+Cz+D_0\) be a centered monic quintic,
and set \(f=3|A|/10\) and \(g=|B|/10\). If \(f+g>0\), then every
target line \(\mathcal L\) satisfies
\begin{equation}\label{eq:new-cubic}
 \Delta(F,\mathcal L)>\max\{(7/4)\sqrt f,\ (5/2)g^{1/3}\}.
\end{equation}
\end{lem}

\begin{proof}
\emph{Step 1. Normalizing the coefficients and phases.}
Write \(\mathcal L=w_0+e^{i\psi}\R\). By Lemma~\ref{lem:cubic-slice}, the
deficit is at least the area of the set of \((\theta,h)\) on which
\(D\) in \eqref{eq:new-D} is positive. The second derivative of the
real slice is cubic, so neither \(C\) nor \(D_0\) enters \(D\).

Put \(\varphi=\arg A\), taking any phase for a zero coefficient.
Replace the old angular variable by \(\theta+\varphi/2\), and set
\(\alpha=5\varphi/2-\psi\), \(\beta=\arg B-3\varphi/2\).
The quantities in \eqref{eq:new-D} become
\[
 E=e^{i(5\theta+\alpha)},\qquad J=f e^{i(3\theta+\alpha)},\qquad
 Z=f e^{2i\theta},\qquad c=g\sin(2\theta+\alpha+\beta).
\]
The change of angular origin preserves area because
\((\theta,h)\mapsto(\theta+\pi,-h)\) represents the same source
line. Adding \(\pi\) to \(\alpha\) leaves \(D\) unchanged, and
adding \(\pi\) to \(\beta\) replaces \(h\) by \(-h\). Simultaneous
negation of the two phases has the effect of
\((\theta,h)\mapsto(-\theta,-h)\). It therefore suffices to take
\(0\le\alpha\le\pi/2\) and \(0\le\beta\le\pi\).

Under \((f,g,h)\mapsto(s^2f,s^3g,sh)\), with \(s>0\), the
expression \(D\) is multiplied by \(s^6\) and the area by \(s\).
If \((7/4)\sqrt f\ge(5/2)g^{1/3}\), divide out the scale
\(s=\sqrt f\). The normalized pair is \((1,q)\), with
\(q=g/f^{3/2}\le(7/10)^3<11/32\). In the other case take
\(s=g^{1/3}\). The pair is then \((q,1)\), with
\(q=f/g^{2/3}\le(10/7)^2<33/16\). Thus we need areas greater
than \(7/4\) and \(5/2\), respectively, in the two normalized cases.

\emph{Step 2. A lower bound for the area.}
Use \(0\le q\le11/32\), \(|h|\le5/4\) in case A, where \(f=1\),
and \(0\le q\le33/16\), \(|h|\le3/2\) in case B, where \(g=1\).
Write \(H_0\) for the corresponding upper bound on \(|h|\).
For each box in \((q,\alpha,\beta)\), divide
\([0,\pi]\times[-H_0,H_0]\) into \(K\times M\) rectangles.
The product of a box and a rectangle has 32 vertices. If \(w_j\)
are its five coordinate widths and \(U_j\) are the corresponding
upper bounds in \eqref{eq:cubic-hessian}, multilinear interpolation
gives
\begin{equation}\label{eq:quintic-positive-cell}
 \min_{\rm vertices}D>\frac18\sum_j U_jw_j^2
 \quad\Longrightarrow\quad D>0\ \text{throughout the cell}.
\end{equation}
Here \(U_q=U_g\) in A and \(U_q=U_f\) in B. In \(U_\theta\) the
coefficient of \(f^2r^2\) is bounded above by \(64\).

If \(N\) rectangles satisfy \eqref{eq:quintic-positive-cell},
their total area is \(2H_0\pi N/(KM)\) for every parameter in the
box. Since \(\pi>333/106\), the two sufficient integer tests are
\begin{equation}\label{eq:quintic-area-tests}
 3330N>742KM\quad\text{in A},\qquad
 1998N>530KM\quad\text{in B}.
\end{equation}
Table~\ref{tab:quintic-phase-cover} gives the counts.
\par\nopagebreak[4]
\begin{table}[H]
\centering\small
\begin{tabular}{@{}lccrrr@{}}
\toprule
Case & Divisions in \((q,\alpha,\beta)\) & \((K,M)\) &
Tested & Accepted & Least \(N\)\\
\midrule
A & \((8,16,32)\)  & \((256,128)\) & 4096  & 3766  & 7302\\
A & \((16,32,64)\) & \((512,256)\) & 2640  & 2640  & 31959\\
B & \((8,16,32)\)  & \((256,128)\) & 4096  & 923   & 8693\\
B & \((16,32,64)\) & \((512,256)\) & 25384 & 25356 & 34784\\
B & \((32,64,128)\)& \((512,256)\) & 224   & 224   & 36625\\
\bottomrule
\end{tabular}
\par\smallskip
\caption{Subdivision of the quintic phase parameters.
The last column is the least \(N\) among the boxes accepted in that row.}
\label{tab:quintic-phase-cover}
\end{table}

A box failing
\eqref{eq:quintic-area-tests} is bisected in all three parameter
coordinates, and all eight children are tested at the next depth.
For instance, the 28 remaining boxes in the second row for B give
exactly the 224 boxes in its last row.

The supplement's \href{https://github.com/zhangteng2000/ehp-lemniscate-verification-20261108/blob/ca8c451bfbc7312b9c57acc5e672c0b69cd601ae/prior/quintic/QuinticPhaseGrid.cs}{\nolinkurl{QuinticPhaseGrid.cs}} evaluates the
expanded sextic in \(h\) by Horner's rule, with outward rounding
on the \(2^{-24}\) lattice. At a shared vertex the interpolation
allowance covers every incident cell. All trigonometric arguments
are multiples of \(\pi/4096\). Half-angle identities and interval
multiplication give enclosures for the sine at these nodes. \href{https://github.com/zhangteng2000/ehp-lemniscate-verification-20261108/blob/ca8c451bfbc7312b9c57acc5e672c0b69cd601ae/prior/quintic/QuinticPhaseAudit.cs}{A
second program} evaluates the factored expression \eqref{eq:new-D}
using sine enclosures of width at most \(2^{-46}\). Its combined
arithmetic and trigonometric error is less than \(2^{-28}\), within
the allowance \(2^{-25}\). It recomputes the area on every terminal
box and checks that the boxes cover both parameter domains.
Multiplying the resulting area bounds by \(s\) gives
\eqref{eq:new-cubic}.
\end{proof}

\begin{lem}\label{lem:new-shapes}
Suppose that the critical points of a monic quintic \(F\) have
diameter one and that its finite critical values are collinear
with diameter \(d>0\). Define \(f,g\) from the centered translate
of \(F\) as in Lemma~\ref{lem:new-cubic}. Then
\begin{equation}\label{eq:new-source}
 \max\{(7/4)\sqrt f,\ (5/2)g^{1/3}\}
       >\frac{17}{16}(d/2)^{1/5}.
\end{equation}
\end{lem}

\begin{proof}
\emph{Step 1. Coordinates for the four critical points.}
By a rigid motion of the source, put a farthest pair at
\(\zeta_1=0\) and \(\zeta_2=1\).
Rotate the target to restore monicity and call the resulting
polynomial \(P\). The moduli of its centered coefficients agree
with those for \(F\). Write \(\zeta_3,\zeta_4\) for the other two
critical points, counted with multiplicity, and order them by
\(\Re\zeta_3\le\Re\zeta_4\). Interchanging \(0,1\) replaces the
sum of their real parts by \(2-\Re\zeta_3-\Re\zeta_4\), so we
may assume \(\Re\zeta_3+\Re\zeta_4\le1\). Complex conjugation
also allows \(\Im\zeta_3\ge0\). Both points are within distance
one of \(0\) and \(1\). Hence
\((\Re\zeta_3,\Im\zeta_3,\Re\zeta_4,\Im\zeta_4)
 \in[0,1]^3\times[-1,1]\).

Put \(\mu=(1+\zeta_3+\zeta_4)/4\), \(S=\zeta_3+\zeta_4-1\),
and \(T=\zeta_4-\zeta_3\).
Write \(P(\mu+z)=z^5+A'z^3+B'z^2+\cdots\).
The four centered critical points \(\zeta_j-\mu\) are
\(-S/4\pm1/2\) and \(S/4\pm T/2\). Comparing coefficients in
\(P'(\mu+z)=5\prod_{j=1}^4(z+\mu-\zeta_j)\), with
\(\sum_j(\zeta_j-\mu)=0\), gives
\begin{equation}\label{eq:quintic-centered-coefficients}
 \frac{3A'}{10}=-\frac14\sum_j(\zeta_j-\mu)^2
        =-\frac{2+2T^2+S^2}{16},\qquad
 \frac{B'}{10}=-\frac1{12}\sum_j(\zeta_j-\mu)^3
        =\frac{S(1-T^2)}{32}.
\end{equation}
The moduli of these two expressions are \(f\) and \(g\).

For a permutation \((i,j,k,l)\) of \((1,2,3,4)\), mark the pair
\(\zeta_i,\zeta_j\) and put \(m=(\zeta_i+\zeta_j)/2\).
Integration of \(P'(z)=5\prod_{\nu=1}^4(z-\zeta_\nu)\) along
\([\zeta_i,\zeta_j]\) gives
\begin{equation}\label{eq:new-chords}
 P(\zeta_j)-P(\zeta_i)=-\frac56(\zeta_j-\zeta_i)^3
       \left((m-\zeta_k)(m-\zeta_l)+\frac{(\zeta_j-\zeta_i)^2}{20}\right).
\end{equation}
Indeed, on writing \(z=\zeta_i+t(\zeta_j-\zeta_i)\), the term odd
in \(t-1/2\) integrates to zero against \(t(1-t)\). The integrals of
\(1\) and \((t-1/2)^2\) against this weight are \(1/6\) and \(1/120\).
Also \(|z-\zeta_k|,|z-\zeta_l|\le1\) throughout the segment, by
convexity of the two unit disks. The same integral gives
\(|P(\zeta_j)-P(\zeta_i)|\le5|\zeta_j-\zeta_i|^3/6\), and hence \(d\le5/6\).
Formula \eqref{eq:new-chords} also holds for coincident marked points.

\emph{Step 2. Verification on the parameter box.}
Use \eqref{eq:quintic-centered-coefficients} to obtain lower bounds
\(P_2\le f^2\) and \(Q_2\le g^2\), both nonnegative. Use the six
differences in \eqref{eq:new-chords} and \(d\le5/6\) to obtain
an upper bound \(D_2\ge d^2\). Either of the strict inequalities
\begin{equation}\label{eq:quintic-shape-tests}
 16\cdot28^{20}P_2^5>17^{20}D_2^2,\qquad
 64\cdot40^{30}Q_2^5>17^{30}D_2^3
\end{equation}
proves \eqref{eq:new-source} throughout the box. These are the
twentieth and thirtieth powers of the respective comparisons,
with the denominators cleared.

A box can also be discarded if its coordinate bounds violate
the ordering or a distance constraint. For any two critical-value
differences \(w_1,w_2\), collinearity requires
\(\Im(w_1\overline{w_2})=0\). Thus a box can be discarded if
even one enclosure for this imaginary part excludes zero.
Every box left unresolved is bisected in all four coordinates.
At depths \(0\) through \(8\), the respective numbers still
requiring subdivision are
\[
 1,\quad16,\quad135,\quad1043,\quad5095,\quad3475,\quad
 1217,\quad316,\quad0 .
\]
The program \href{https://github.com/zhangteng2000/ehp-lemniscate-verification-20261108/blob/ca8c451bfbc7312b9c57acc5e672c0b69cd601ae/prior/quintic/quintic_shape_integer_cover.py}{\nolinkurl{quintic_shape_integer_cover.py}} evaluates
\eqref{eq:quintic-shape-tests} with outward integer rectangles in
units of \(2^{-23}\). \href{https://github.com/zhangteng2000/ehp-lemniscate-verification-20261108/blob/ca8c451bfbc7312b9c57acc5e672c0b69cd601ae/prior/quintic/quintic_shape_arbitrary_audit.json}{A replay with unbounded integers} gives the
same classifications at every depth. The independent program
\href{https://github.com/zhangteng2000/ehp-lemniscate-verification-20261108/blob/ca8c451bfbc7312b9c57acc5e672c0b69cd601ae/prior/quintic/independent_quintic_shapes.py}{\nolinkurl{independent_quintic_shapes.py}} encloses the quantities
in complex disks, evaluates the primitive directly at its critical
points and bisects the longest coordinate interval. Its tree has
\(347691\) internal nodes and \(347692\) terminal nodes, with
maximum depth \(38\). Both checks cover the original closed box,
including configurations with repeated critical points or values.
\end{proof}

We now prove Proposition~\ref{prop:new-five-line}.
\begin{proof}
The diameter \(\rho\) of the critical points is positive since
\(d>0\). After centering, take \(P(z)=\rho^{-5}F(\rho z)\) and
\(\mathcal L'=\rho^{-5}\mathcal L\). The critical points of \(P\) have diameter
one and their images have diameter \(d/\rho^5\).
Lemma~\ref{lem:new-shapes} gives \(f+g>0\) for \(P\). Applying
\eqref{eq:new-cubic} and then \eqref{eq:new-source} gives
\(\Delta(P,\mathcal L')>(17/16)(d/(2\rho^5))^{1/5}\).
Equation \eqref{eq:low-degree-line-scaling} now gives
\eqref{eq:new-line-loss}.
\end{proof}

\section{The quintic comparison}\label{sec:new-five}

By Lemma~\ref{lem:extremizers} we may choose a quintic maximizer
whose finite critical values lie on \(\mathbb T\).
Proposition~\ref{prop:new-five-line} verifies the hypothesis of
Proposition~\ref{prop:small-circle-transfer} with \(L_5=17/16\).
Hence any maximizer with this property has only one critical value
\(v\). By Lemma~\ref{lem:critical-compactness} it has the form
\(p(z)=(z-z_0)^5+v\), with \(|v|=1\). Each such polynomial has
the same lemniscate length as \(z^5-1\), by translation and rotation.
This proves the length comparison for every monic quintic.
Section~\ref{sec:equality} treats the equality assertion for an
arbitrary maximizer.

\section{The quartic comparison}\label{sec:quartic-global}

For the quartic estimates below, write \(p_0(z)=z^4-1\) and let
\(d_p\) denote the diameter of the finite critical values of \(p\).

\begin{prop}\label{prop:quartic-boundary-gap}
If \(p\) is a monic quartic, all its critical values lie on
\(\mathbb T\), and \(d_p>0\), then
\[
 \ell(\partial E_1(p_0))-\ell(\partial E_1(p))
       >2^{-5}d_p^{1/4}.
\]
\end{prop}

The research notes of Plasma AI~\cite[entries L179--L181]{Pla26}
give a quartic line-deficit estimate, a length comparison at
small curvature, and estimates for the two families with
coincident critical values. Here we prove \eqref{eq:quartic-line}
for every monic quartic with collinear critical values of positive
diameter.
A quartic has three finite critical values, counted with
multiplicity. When these lie on a line, their positions up to
similarity are determined by the position of the third value
between the two extremes. We first use this parameter to bound
the line deficit.

\begin{lem}\label{lem:quartic-line}
Let \(F\) be a monic quartic whose finite critical values lie on
a line \(\mathcal L\) and have diameter \(d>0\). Then
\begin{equation}\label{eq:quartic-line}
 \Delta(F,\mathcal L)>\frac{29}{20}(d/2)^{1/4}.
\end{equation}
\end{lem}

\begin{proof}
\emph{Step 1. Parametrizing the critical points.}
Choose critical points over the two extreme critical values.
By an affine change of source put them at \(-1,1\), and denote the
third critical point by \(\zeta\). Restore monicity in the target and
choose the constant of integration so that the resulting polynomial
is \(G_\zeta(z)=z^4-\tfrac43\zeta z^3-2z^2+4\zeta z\). Indeed,
\(G_\zeta'(z)=4(z^2-1)(z-\zeta)\). The difference of the two extreme
values is \(d_\zeta=G_\zeta(1)-G_\zeta(-1)=16\zeta/3\), so \(\zeta\ne0\).
The ratio \(\nu=(G_\zeta(\zeta)-G_\zeta(-1))/d_\zeta\) is real and lies in
\([0,1]\). Swapping the endpoints replaces \(\nu\) by \(1-\nu\),
so we may take \(0\le\nu\le1/2\). Direct evaluation gives
\begin{equation}\label{eq:quartic-line-parameter}
 \nu=\frac{(3-\zeta)(1+\zeta)^3}{16\zeta},\qquad
 \zeta^4-6\zeta^2+8(2\nu-1)\zeta-3=0 .
\end{equation}

For \(0\le t\le1\) put \(s=1-t^2\), \(r=\sqrt{1+s+s^2}\),
\(u=\sqrt{2r+2+s}\), and \(v=\sqrt3\,s/u\).
All four roots in \eqref{eq:quartic-line-parameter} are given by
\begin{equation}\label{eq:quartic-line-charts}
 \zeta=t+u,\quad t-u,\quad -t+iv,\quad -t-iv,
 \qquad \nu=(1-tr)/2 .
\end{equation}
Indeed \(u^2-v^2=4+2s\) and \(u^2v^2=3s^2\). Expanding
\(((z-t)^2-u^2)((z+t)^2+v^2)\) gives the quartic in
\eqref{eq:quartic-line-parameter}. Moreover,
\(t^2r^2=1-s^3\), so \(tr\) increases from zero to one and
\(\nu\) runs through the whole interval. At \(t=1\) the parameter
equation becomes \((\zeta-3)(\zeta+1)^3=0\). For \(\zeta=-1\),
\(G_\zeta\) has a double critical point at \(-1\). For \(\zeta=3\),
its critical points are distinct, but \(G_3(-1)=G_3(3)=-9\).
Conjugation preserves length and interchanges the last two branches
of \eqref{eq:quartic-line-charts}.
It is therefore enough to check the first three on their closed
parameter intervals.

\emph{Step 2. The real slices.}
The centered monic polynomial
\(P_\zeta(z)=G_\zeta(z+\zeta/3)-G_\zeta(-1)\) has critical values on \(\zeta\R\).
Its quadratic and linear coefficients are
\(U=-2-2\zeta^2/3\) and \(V=8\zeta/3-8\zeta^3/27\), respectively.
Put \(\varepsilon=\bar \zeta/|\zeta|\). A source point \(z\) lies over
\(\zeta\R\) exactly when \(\Im(\varepsilon P_\zeta(z))=0\).
Thus Lemma~\ref{lem:cubic-slice} applies to the monic polynomial
\(P_\zeta\) with target direction \(\zeta/|\zeta|\). In \eqref{eq:new-D} use
\[
 E=\varepsilon e^{4i\theta},\quad
 J=\tfrac12\varepsilon U e^{2i\theta},\quad
 c_1=\Im\bigl(\tfrac14\varepsilon V e^{i\theta}\bigr),\quad
 Z=\tfrac12\bar U e^{2i\theta},
\]
where \(c_1\) occupies the place of the real coefficient denoted
by \(c\) in that formula.

\emph{Step 3. The area estimate.}
Divide each of the three parameter intervals in
\eqref{eq:quartic-line-charts} into 128 equal closed intervals.
On each interval in \(t\), enclose \(\varepsilon,\ \varepsilon U/2,
\varepsilon V/4\), and \(\overline U/2\).
Divide \([0,\pi]\times[-3,3]\) into \(256\times256\) rectangles.
At each angular grid point, multiply the four enclosures by
\(e^{4i\theta}, e^{2i\theta}, e^{i\theta}\), and \(e^{2i\theta}\),
respectively. The imaginary part of the third product gives
\(c_1\). Only \(\theta\) and \(h\) are interpolated.
With \(f=|U|/2\) and \(g=|V|/4\), the bounds \(U_\theta,U_h\)
in \eqref{eq:cubic-hessian} apply uniformly on the parameter
interval. The coefficient of \(f^2r^2\) in \(U_\theta\) for \(n=4\)
is \(56/3\), where \(r\) here denotes the upper bound on \(|h|\)
in that estimate. A rectangle is counted if its four vertex lower
bounds for \(D\) all exceed
\(\bigl(U_\theta(\pi/256)^2+U_h(6/256)^2\bigr)/8\).
The interpolation remainder then gives \(D>0\) throughout the
rectangle for every \(t\) in the parameter interval.

For the interval shown in Figure~\ref{fig:ehp-quartic-cells},
the accepted rectangles contribute \(7{,}641\cdot6\pi/256^2\)
to the area bound.
\par\nopagebreak[4]
\begin{figure}[H]
\centering
\includegraphics[width=\linewidth]{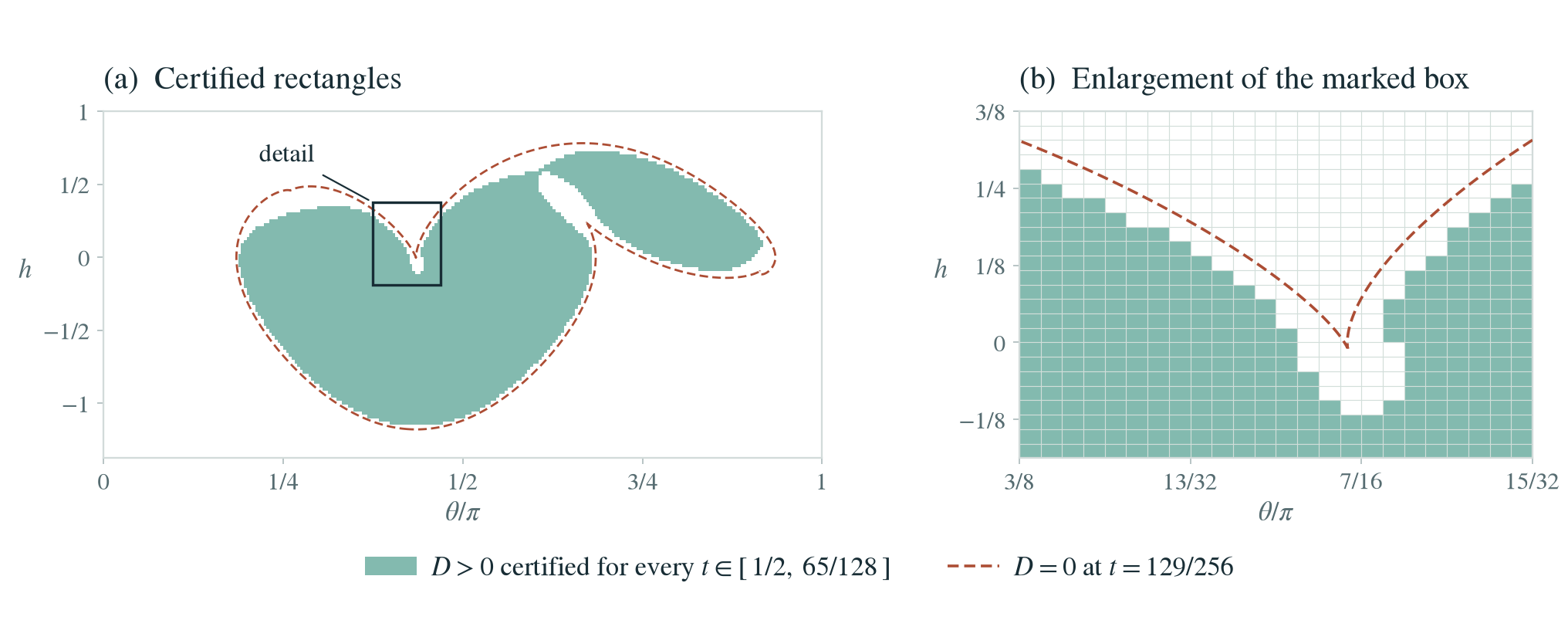}
\caption{One parameter interval in the proof of Lemma~\ref{lem:quartic-line},
on the branch \(\zeta=-t+iv\). The shaded cells are the \(7{,}641\)
rectangles on which the factored interval calculation certifies \(D>0\)
for every \(t\in[1/2,65/128]\). Each cell has widths \(\pi/256\) in
\(\theta\) and \(6/256\) in \(h\). The left panel contains all accepted
cells. The dashed curve shows the numerically plotted level \(D=0\)
at the midpoint \(t=129/256\).}
\label{fig:ehp-quartic-cells}
\end{figure}

For an interval \(I\), let \(N_I\) be the count and let
\(\gamma_I\) be a lower bound on
\(\gamma=(2/|d_\zeta|)^{1/4}=(3/(8|\zeta|))^{1/4}\) for \(t\in I\).
Writing \(\mathcal L_\zeta=G_\zeta(-1)+\zeta\R\), we obtain
\(\gamma\Delta(G_\zeta,\mathcal L_\zeta)\ge\gamma_I\,6\pi N_I/256^2\).
Table~\ref{tab:quartic-line-cover} gives strict lower bounds
for this quantity on each branch.
\par\nopagebreak[4]
\begin{table}[H]
\centering\small
\begin{tabular}{@{}lcc@{}}
\toprule
Branch & Expanded expression for \(D\) & Factored expression for \(D\)\\
\midrule
\(t+u\)      & \(61/25\) & \(49/20\)\\
\(t-u\)      & \(37/25\) & \(38/25\)\\
\(-t\pm iv\) & \(73/50\) & \(151/100\)\\
\bottomrule
\end{tabular}
\par\smallskip
\caption{Lower bounds for the quartic line deficit after the
critical-value diameter is normalized to two.}
\label{tab:quartic-line-cover}
\end{table}

The program \href{https://github.com/zhangteng2000/ehp-lemniscate-verification-20261108/blob/ca8c451bfbc7312b9c57acc5e672c0b69cd601ae/verification/quartic_line/verify_quartic_line.py}{\nolinkurl{quartic_line/verify_quartic_line.py}}
reconstructs all 384 intervals in two ways. The first uses
\(v=\sqrt3\,s/u\) with outward rounding on a \(2^{-60}\) lattice.
The second keeps rational operations exact and encloses square roots
on a \(2^{-100}\) lattice. It uses \(v^2=2r-2-s\), which is
nondecreasing in \(s\) since \((1+2s)/r\ge1\). Both computations round the coefficient bounds outwards to multiples
of \(2^{-24}\) and use \(\pi>333/106\) in the area estimate.
They evaluate the expanded and factored expressions for \(D\)
separately, using the same interval operations and interpolation
error routine.
Every entry of Table~\ref{tab:quartic-line-cover} exceeds \(29/20\).
Equation \eqref{eq:low-degree-line-scaling} now gives
\eqref{eq:quartic-line}.
\end{proof}

We return to Proposition~\ref{prop:quartic-boundary-gap}.
\begin{proof}
Translate the source to center \(p\), and put \(\delta=d_p/2\),
\(b=3/8\). We shall prove the slightly stronger estimate
\begin{equation}\label{eq:quartic-strengthened-gap}
 \ell(\partial E_1(p_0))-\ell(\partial E_1(p))
       >\frac{3}{40}\delta^{1/4}.
\end{equation}

We first state two consequences of the deformation estimates.
Choose critical values \(v_1,v_2\) with \(|v_1-v_2|=2\delta\).
For almost every \(\xi\in\mathbb T\), their images under \(M_\xi\)
are finite critical values of \(F_\xi\) on \(\mathcal L_\xi\), and their
distance is \(2\delta/(|1-\xi v_1||1-\xi v_2|)\). Apply
\eqref{eq:quartic-line} and average in \(\xi\). For \(v\in\mathbb T\)
the mean of \(\log|1-\xi v|\) is zero. Jensen's inequality for the
exponential therefore gives
\begin{equation}\label{eq:quartic-mean-line}
 \big\langle\Delta(F_\xi,\mathcal L_\xi)\big\rangle
       \ge \frac{29}{20}\delta^{1/4}.
\end{equation}
The logarithms are integrable, and \(v_1\ne v_2\), so the product
\(|1-\xi v_1|^{-1/4}|1-\xi v_2|^{-1/4}\) is integrable as well.

For a centered monic quartic \(q\), let \(a_1(q)=[w^{-1}]f_q(w)\).
Since \(\sqrt{1-s}-1-\tfrac12\log(1-s)=s^2/8+O(s^3)\) at zero,
Lemma~\ref{lem:inverse-log-identity} gives
\begin{equation}\label{eq:quartic-first-coefficient}
 h_4(q)=\frac18\sum_{q'(\zeta)=0}q(\zeta)+\frac18a_1(q)^2,
\end{equation}
where critical points are counted with multiplicity.

\emph{Step 1. The short arc.}
Put \(\tau_0=1/64\) and
\(\delta_0=2\tau_0/(1+\tau_0^2)=128/4097\).
If \(0<\delta\le\delta_0\), use the rotation and normalization in
\eqref{eq:short-arc-normalization}. The resulting polynomial \(F\)
has real critical values in \([-1,1]\), including both endpoints,
and
\[
 p(z)=1+2i\tau F_{i\tau}\bigl((2i\tau)^{-1/4}z\bigr),
 \qquad \delta=\frac{2\tau}{1+\tau^2},\qquad 0<\tau\le\tau_0 .
\]
By \eqref{eq:moved-segment-capacity} the critical values of
\(F_{i\tau}\) lie on an arc of capacity at most \(1/2\).
Lemma~\ref{lem:last-arc} gives \(|a_1(F_{i\tau})|^2\le1/2\).
For \(r\in[-1,1]\),
\(\Im M_{i\tau}(r)=\tau r^2/(1+\tau^2r^2)\le\tau\).
Hence \eqref{eq:quartic-first-coefficient} yields
\(\max(0,\Im h_4(F_{i\tau}))\le1/16+3\tau/8\).

Take \(R=5/2\) and \(b_0=2/5\) in
Lemma~\ref{lem:curvature-general}. Then
\(\tau<b_0<1-1/R\), \(\tau R\le5/128\), and
\(s=R/\sqrt{1+\tau^2R^2}>249/100\).
Uniform upper bounds for the four groups on the right of
\eqref{eq:curvature-general}, in their order there with the first
two summands grouped together, are
\[
 \frac{631}{1000},\qquad \frac{253}{1000},\qquad
 \frac{58}{1000},\qquad \frac{433}{1000}.
\]
For clarity, the scalar enclosures use \(H^2<7/8\), \(H<15/16\),
\(B(3/8,1/2)<39/10\), and \(B(1/8,1/2)<19/2\).
For the row sum \(\mathcal T\), retain \(2\le j,k<200\).
Writing \(y=(249/100)^{-1/4}\), the sum of all omitted entries is
bounded by \(\sum_{m\ge202}m(6+2m)y^m\), since an entry with
\(j+k=m\) has kernel at most \(6+2m\) and there are at most \(m\)
such entries. For \(\mathcal A\), retain \(2\le q<24\) and bound
the remaining \(\ell^2\) norm by
\(\sum_{q\ge24}(3+2q)(249/100)^{-q}\).
The
\href{https://github.com/zhangteng2000/ehp-lemniscate-verification-20261108/blob/ca8c451bfbc7312b9c57acc5e672c0b69cd601ae/verification/quartic_line/independent_circle_audit.py}{independent quartic verification}
checks these geometric tails, the beta enclosures from
\eqref{eq:split-beta-series}, and each of the four bounds above
by rational arithmetic.

The four displayed bounds sum to \(11/8\). Since
\(\Delta(F,\R)>29/20\) by \eqref{eq:quartic-line}, and
\(29/20-11/8=3/40\), formula \eqref{eq:curvature-general} gives
\[
 \ell(\partial E_1(p_0))-\ell(\partial E_1(p))
     >\frac{3}{40}(2\tau)^{1/4}.
\]
As \(2\tau=(1+\tau^2)\delta\), this proves
\eqref{eq:quartic-strengthened-gap} for \(0<\delta\le\delta_0\).

\emph{Step 2. The remaining diameters.}
The scalar enclosures also give \(U_4<2353/200\) and
\(\ell(\partial E_1(p_0))>1107/100\), so their difference is
less than \(139/200\). Keeping \(\mathcal E_1\) in \eqref{eq:new-mass}
and using \eqref{eq:quartic-mean-line} gives
\begin{equation}\label{eq:quartic-diameter-function}
 \ell(\partial E_1(p_0))-\ell(\partial E_1(p))
   > \frac{29}{20}\delta^{1/4}+\mathcal E_1-\frac{139}{200}.
\end{equation}
For \(\delta\ge1/15\), the line term alone suffices.
Indeed \(15\cdot139^4<275^4\), so
\(139/200<(11/8)\delta^{1/4}\). Since \(\mathcal E_1\ge0\),
\eqref{eq:quartic-diameter-function} gives
\eqref{eq:quartic-strengthened-gap}.

It remains to treat \(\delta_0\le\delta\le1/15\). Rotate the
marked pair \(v_1,v_2\) to \(e^{\pm i\alpha}\), with
\(\sin\alpha=\delta\) and \(0<\alpha<\pi/3\), restoring monicity
by a source rotation. Every critical value is within angular
distance \(2\alpha\) of both endpoints, so it lies on the shorter
arc between them. This arc lies in the closed disk of center
\(c=\cos\alpha\) and radius \(\delta\). Since the disk has capacity
\(\delta\), Lemma~\ref{lem:last-arc} gives
\(|a_1(p)|^2\le\delta\). The three critical values have real parts
at least \(c\), and \eqref{eq:quartic-first-coefficient} now gives
\[
 |h_4(p)|\ge b\sqrt{1-\delta^2}-\delta/8
       \ge b(1-\delta^2)-\delta/8
       \ge b(1-15^{-2})-\frac1{120}=\frac{73}{200}.
\]
By \eqref{eq:exterior-coefficients} with \(k=1\) and \(\pi>157/50\),
\(\mathcal E_1\ge2\pi b^2|h_4(p)|^2>47/400\).
The remaining arithmetic is
\[
 \frac{139}{200}-\frac{47}{400}
     =\frac{11}{8}\frac{21}{50},
 \qquad
 (21/50)^4<\delta_0 .
\]
Thus \(139/200-\mathcal E_1<(11/8)\delta^{1/4}\) throughout this
interval. Substitution in \eqref{eq:quartic-diameter-function}
proves \eqref{eq:quartic-strengthened-gap} in the last case.
Finally \(d_p=2\delta\) and \(2^{1/4}<6/5\), so the coefficient of
\(d_p^{1/4}\) in \eqref{eq:quartic-strengthened-gap} is
\(3/(40\,2^{1/4})>1/16\). This implies the stated bound.
\end{proof}

By Lemma~\ref{lem:extremizers}, at least one quartic maximizer has
all its critical values on \(\mathbb T\).
Proposition~\ref{prop:quartic-boundary-gap} forces \(d_p=0\) for
every maximizer with this property, since \(p_0\) is itself an
admissible quartic. Lemma~\ref{lem:critical-compactness} therefore
gives \(p(z)=(z-z_0)^4+v\), \(|v|=1\), for each such maximizer.

We can now finish the length comparison in Theorem~\ref{thm:main}.
\begin{proof}
In degree one the lemniscate is a circle of radius one. In degree
two the centered maximizer supplied by Lemma~\ref{lem:extremizers}
has the form \(z^2+v\) with \(|v|=1\). In degree three,
Lemma~\ref{lem:cubic-comparison} gives the bound and makes it
strict when two distinct critical values lie on \(\mathbb T\).
The quartic case follows from Proposition~\ref{prop:quartic-boundary-gap},
and Section~\ref{sec:new-five} treats the quintic case.
For \(6\le n\le8\), Propositions~\ref{prop:last-line}
and~\ref{prop:octic-line} give the line estimates required by
Proposition~\ref{prop:small-circle-transfer}.
Section~\ref{sec:collinear-circle} gives the comparison for
\(9\le n\le27\). For \(28\le n\le255\) and \(n\ge256\), respectively,
Proposition~\ref{prop:finite-line} and Theorem~\ref{thm:line256}
supply the hypothesis of Proposition~\ref{prop:uniform-circle-transfer}.
This proves the length inequality in every degree. The equality
assertion for arbitrary maximizers is proved in
Section~\ref{sec:equality}.
\end{proof}
\section{Equality cases}\label{sec:equality}

The degree comparisons show that a maximizer whose finite critical
values lie on \(\mathbb T\) has only one such value when \(n\ge2\).
It remains to show that every maximizer belongs to this class.
Eremenko and Hayman~\cite[proof of Lemma~5 in the corrected
version]{EH99} move one critical value and use subharmonicity
of the length along the resulting deformation. Differentiating
their deformation also yields the strict conclusion below.
We first give a polynomial deformation for which the second
variation can be computed explicitly.

\begin{lem}\label{lem:new-strict-levi}
Let \(n\ge2\). If \(p\) locally maximizes
\(q\mapsto\ell(\partial E_1(q))\) in the coefficient space of
centered monic polynomials of degree \(n\), then
\(|p(\zeta)|=1\) for every critical point \(\zeta\) of \(p\).
\end{lem}

\begin{proof}
\emph{Step 1. Transporting the lemniscate.}
Suppose that \(p'(\zeta_1)=0\) and \(|p(\zeta_1)|\ne1\). Let \(S\) be the set
of distinct critical values of \(p\) on \(\mathbb T\), and put
\(G(w)=\prod_{v\in S}(w-v)\), with \(G=1\) if \(S\) is empty.
The target flow \(\dot w=G(w)\) fixes every value in \(S\), while
its velocity at \(p(\zeta_1)\) is nonzero.
For a monic polynomial \(F\) of degree \(n\), perform the division
\begin{equation}\label{eq:equality-division}
 G(F)=F'Q_F+R_F,\qquad \deg R_F\le n-2.
\end{equation}
Let \(q_t\) solve \(\dot q_t=R_{q_t}\), \(q_0=p\), for a complex
parameter \(t\). Since \(F'\) has leading coefficient \(n\), the
coefficients of \(Q_F\) and \(R_F\) are polynomial functions of the
coefficients of \(F\). The solution is therefore holomorphic for
small \(t\), and the degree bound in \eqref{eq:equality-division}
ensures that \(q_t\) remains monic and centered.

The rational function \(W=G(p)/p'\) is holomorphic on a neighborhood
of \(\partial E_1(p)\). Indeed, if \(\zeta\) is a critical point on this
lemniscate and has local degree \(m\), then \(p-p(\zeta)\) vanishes to
order \(m\) at \(\zeta\), whereas \(p'\) vanishes to order \(m-1\).
The factor \(w-p(\zeta)\) in \(G\) thus cancels the apparent pole.
The remaining zeros of \(p'\) form a finite set disjoint from the
lemniscate.

Let \(\Phi_t\) be the source flow
\(\dot\Phi_t=Q_{q_t}\circ\Phi_t\), and let \(\chi_t\) be the flow
of the fixed vector field \(W\), both starting at the identity.
By compactness of \(\partial E_1(p)\) and local existence and
uniqueness for holomorphic differential equations, the flows and
the composition below are defined and univalent on a fixed
neighborhood of the lemniscate for all sufficiently small \(t\).
As in \eqref{eq:new-flow}, the chain rule and
\eqref{eq:equality-division} show that \(q_t(\Phi_t(z))\) and
\(p(\chi_t(z))\) both solve \(\dot w=G(w)\) with initial value
\(p(z)\). Since \(W\) is independent of \(t\), its inverse flow
is \(\chi_{-t}\). We obtain
\begin{equation}\label{eq:equality-fixed-level}
 \Psi_t=\Phi_t\circ\chi_{-t},\qquad q_t(\Psi_t(z))=p(z).
\end{equation}

For each \(w\in\mathbb T\), the map \(\Psi_t\) takes the distinct
points of \(p^{-1}(w)\) to distinct points of \(q_t^{-1}(w)\).
By \eqref{eq:equality-fixed-level} and local univalence, it preserves
the multiplicity at each point. These multiplicities sum to \(n\),
so the images exhaust the fibre of \(q_t\), even when \(w\) is a
critical value. Thus \(\Psi_t\) maps the original lemniscate onto
the whole of \(\partial E_1(q_t)\). Integrating over its regular
arcs gives
\begin{equation}\label{eq:equality-transported-length}
 \ell(\partial E_1(q_t))
   =\int_{\partial E_1(p)}|\Psi_t'(z)|\,ds .
\end{equation}
The finitely many critical points have arclength measure zero.

\emph{Step 2. The Laplacian of the length.}
The source velocity at \(t=0\) is
\(V=\left.\partial_t\Psi_t\right|_{t=0}=Q_p-W=-R_p/p'\).
For a holomorphic function \(g\) with \(g(0)=1\), one has
\(\left.\partial_t\partial_{\bar t}|g(t)|\right|_{t=0}
=|g'(0)|^2/4\). Apply this identity to \(g(t)=\Psi_t'(z)\).
On the compact lemniscate \(\Psi_t'\) stays uniformly close to \(1\)
for small \(t\), and Cauchy's estimates bound its parameter
derivatives uniformly. We may therefore differentiate
\eqref{eq:equality-transported-length} under the integral to obtain
\begin{equation}\label{eq:equality-second-variation}
 \left.\partial_t\partial_{\bar t}
       \ell(\partial E_1(q_t))\right|_{t=0}
   =\frac14\int_{\partial E_1(p)} |V'(z)|^2\,ds .
\end{equation}

At the off-circle critical point \(\zeta_1\), \eqref{eq:equality-division}
gives \(R_p(\zeta_1)=G(p(\zeta_1))\ne0\). Hence \(V\) has a pole at \(\zeta_1\) and
is a nonconstant rational function. Its derivative has only finitely
many zeros, while the preimage of a small arc of regular values
on \(\mathbb T\) consists of arcs of positive length. Thus the
integral in \eqref{eq:equality-second-variation} is strictly positive.
At a local maximum the Laplacian
\(4\partial_t\partial_{\bar t}\ell(\partial E_1(q_t))\) would be
nonpositive at \(t=0\). This contradiction proves the lemma.
\end{proof}

For comparison, consider the deformation of Eremenko and Hayman
described above. If \(v=p(\zeta)\notin\mathbb T\) is a critical
value, their family consists of monic polynomials \(q_t\), with
\(q_0=p\), and source maps \(\psi_t\), with \(\psi_0(z)=z\).
The maps depend holomorphically on \(t\), satisfy
\(q_t(\psi_t(\zeta))=v+t\), and obey \(q_t\circ\psi_t=p\) on a
fixed neighborhood of \(\partial E_1(p)\) where \(\psi_t\) is
conformal. Moreover, \(\psi_t\) maps \(\partial E_1(p)\) onto
\(\partial E_1(q_t)\). Centering \(q_t\) by translation preserves its length.
At a local maximum, the calculation in
\eqref{eq:equality-second-variation}, applied to \(\psi_t\),
forces \(V' = 0\) on every regular arc, where
\(V=\left.\partial_t\psi_t\right|_{t=0}\).
On one such arc write \(V=C\), and put
\(\dot q=\left.\partial_tq_t\right|_{t=0}\). Differentiating
\(q_t\circ\psi_t=p\) at zero gives \(\dot q=-Cp'\) on the arc.
Both sides are polynomials, so the identity holds throughout \(\C\).
Thus \(\dot q(\zeta)=0\), whereas differentiating
\(q_t(\psi_t(\zeta))=v+t\) gives \(\dot q(\zeta)=1\).

We now finish the proof of Theorem~\ref{thm:main}. Suppose \(n\ge2\)
and \(p\) attains the model length. After a source translation, it is
a centered maximizer. Lemma~\ref{lem:new-strict-levi} puts all its
critical values on \(\mathbb T\), and the degree comparisons force
them to coincide at a single value \(v\). By
Lemma~\ref{lem:critical-compactness}, in the original coordinates
\(p(z)=(z-z_0)^n+v\). Conversely, if \(|v|=1\), choose
\(\lambda\in\mathbb T\) with \(\lambda^n=-v\). Then
\(p(z_0+\lambda z)=\lambda^n p_0(z)\), where \(p_0(z)=z^n-1\).
Thus \(\partial E_1(p)=z_0+\lambda\,\partial E_1(p_0)\), and the
two lengths agree. Writing \(v=-e^{i\theta}\) gives precisely
the asserted equality cases. In degree one every monic polynomial
has this form and its lemniscate is a circle of length \(2\pi\).

\appendix
\section{Verification of the finite estimates}\label{sec:certificates}

The \href{https://github.com/zhangteng2000/ehp-lemniscate-verification-20261108/blob/ca8c451bfbc7312b9c57acc5e672c0b69cd601ae/README.md}{computational supplement}~\cite{Zha26} contains the rational input data and
the programs for the finite estimates used in the proof. From the
root of the supplement, run
\begin{verbatim}
python -B verification/replay_all_degrees.py --workers 2
\end{verbatim}
with Python assertions enabled. The reference environment uses
64-bit Python~3.12 and a 64-bit C\# compiler. The Python package
versions are recorded in \href{https://github.com/zhangteng2000/ehp-lemniscate-verification-20261108/blob/ca8c451bfbc7312b9c57acc5e672c0b69cd601ae/requirements.txt}{\nolinkurl{requirements.txt}}. \href{https://github.com/zhangteng2000/ehp-lemniscate-verification-20261108/blob/ca8c451bfbc7312b9c57acc5e672c0b69cd601ae/verification/replay_all_degrees.py}{The driver}
compiles the counters from source. One of the quintic checks uses
NumPy arrays of arbitrary-precision integers, and an independent
check uses disks in the complex plane with integer bounds for
their centers and radii after scaling. The symbolic cross-checks
use SymPy.

\subsection{Enclosures for the beta integral}

For \(0<t\le2\), splitting the beta integral at \(\theta=\pi/4\)
gives
\begin{equation}\label{eq:split-beta-series}
\begin{aligned}
 \frac12 B\!\left(\frac t2,\frac12\right)
 &=\int_0^{\pi/2}\sin^{t-1}\theta\,d\theta\\
 &=2^{-t/2}\sum_{k=0}^{\infty}
       \frac{(1/2)_k\,2^{-k}}{k!(2k+t)}
    +2^{-1/2}\sum_{k=0}^{\infty}
       \frac{(1-t/2)_k\,2^{-k}}{k!(2k+1)}.
\end{aligned}
\end{equation}
On the first interval substitute \(u=\sin\theta\), and on the
second substitute \(u=\cos\theta\). The resulting integrands on
\([0,2^{-1/2}]\) are \(u^{t-1}(1-u^2)^{-1/2}\) and
\((1-u^2)^{-(1-t/2)}\), respectively. Expanding the binomial
factors and integrating term by term gives \eqref{eq:split-beta-series}.
The expansions converge uniformly on this interval, and
\(u^{t-1}\) is integrable at zero.

Let \(I_N(t)\) denote the right-hand side of
\eqref{eq:split-beta-series} with both sums restricted to
\(0\le k<N\). For \(0\le q\le1\), the product formula gives
\(0\le(q)_k/k!\le1\). Summing the two geometric tails gives
\[
 0\le \frac12 B\!\left(\frac t2,\frac12\right)-I_N(t)
 \le 2^{1-N}
 \left(\frac{2^{-t/2}}{2N+t}+\frac{2^{-1/2}}{2N+1}\right).
\]
For \(N=100\) this bounds the error in
\(\tfrac12 B(t/2,1/2)\) by \(2^{-105}\), uniformly for \(0<t\le2\).
Rational powers are enclosed by comparing integer powers of rational
endpoints. For bounds on \(\pi\), one can use the alternating
arctangent series in Machin's identity
\(\pi=16\arctan(1/5)-4\arctan(1/239)\).

\subsection{Programs for the degree estimates}

Table~\ref{tab:degree-verifiers} lists the principal programs.
\par\nopagebreak[4]
\begin{table}[H]
\centering\small
\setlength{\tabcolsep}{4pt}
\begin{tabularx}{\textwidth}{@{}ll>{\raggedright\arraybackslash}X@{}}
\toprule
Degree & Result & Programs\\
\midrule
\multicolumn{3}{@{}l}{\emph{Estimates on lines}}\\
\(n\ge256\) & Theorem~\ref{thm:line256} &
  \href{https://github.com/zhangteng2000/ehp-lemniscate-verification-20261108/blob/ca8c451bfbc7312b9c57acc5e672c0b69cd601ae/verification/effective_line_constants.py}{\nolinkurl{effective_line_constants.py}},
  \href{https://github.com/zhangteng2000/ehp-lemniscate-verification-20261108/blob/ca8c451bfbc7312b9c57acc5e672c0b69cd601ae/verification/audit_effective_line_256.py}{\nolinkurl{audit_effective_line_256.py}},
  \href{https://github.com/zhangteng2000/ehp-lemniscate-verification-20261108/blob/ca8c451bfbc7312b9c57acc5e672c0b69cd601ae/verification/inverse_beta_moment_audit.py}{\nolinkurl{inverse_beta_moment_audit.py}}\\
\(28\le n\le255\) & Proposition~\ref{prop:finite-line} &
  \href{https://github.com/zhangteng2000/ehp-lemniscate-verification-20261108/blob/ca8c451bfbc7312b9c57acc5e672c0b69cd601ae/verification/verify_interpolation_degrees.py}{\nolinkurl{verify_interpolation_degrees.py}},
  \href{https://github.com/zhangteng2000/ehp-lemniscate-verification-20261108/blob/ca8c451bfbc7312b9c57acc5e672c0b69cd601ae/verification/independent_interpolation_check.py}{\nolinkurl{independent_interpolation_check.py}}\\
\(15\le n\le27\) & Proposition~\ref{prop:collinear-finite} &
  \href{https://github.com/zhangteng2000/ehp-lemniscate-verification-20261108/blob/ca8c451bfbc7312b9c57acc5e672c0b69cd601ae/verification/verify_peeled_line.py}{\nolinkurl{verify_peeled_line.py}},
  \href{https://github.com/zhangteng2000/ehp-lemniscate-verification-20261108/blob/ca8c451bfbc7312b9c57acc5e672c0b69cd601ae/verification/audit_peeled_line.py}{\nolinkurl{audit_peeled_line.py}}\\
\(9\le n\le14\) & Proposition~\ref{prop:conditional-line} &
  \href{https://github.com/zhangteng2000/ehp-lemniscate-verification-20261108/blob/ca8c451bfbc7312b9c57acc5e672c0b69cd601ae/verification/verify_conditional_line.py}{\nolinkurl{verify_conditional_line.py}},
  \href{https://github.com/zhangteng2000/ehp-lemniscate-verification-20261108/blob/ca8c451bfbc7312b9c57acc5e672c0b69cd601ae/verification/audit_conditional_line.py}{\nolinkurl{audit_conditional_line.py}},
  \href{https://github.com/zhangteng2000/ehp-lemniscate-verification-20261108/blob/ca8c451bfbc7312b9c57acc5e672c0b69cd601ae/verification/verify_reciprocal_strips.py}{\nolinkurl{verify_reciprocal_strips.py}}\\
\(n=8\) & Proposition~\ref{prop:octic-line} &
  \href{https://github.com/zhangteng2000/ehp-lemniscate-verification-20261108/blob/ca8c451bfbc7312b9c57acc5e672c0b69cd601ae/verification/verify_conditional_line.py}{\nolinkurl{verify_conditional_line.py}},
  \href{https://github.com/zhangteng2000/ehp-lemniscate-verification-20261108/blob/ca8c451bfbc7312b9c57acc5e672c0b69cd601ae/verification/verify_root_product_strip.py}{\nolinkurl{verify_root_product_strip.py}},
  \href{https://github.com/zhangteng2000/ehp-lemniscate-verification-20261108/blob/ca8c451bfbc7312b9c57acc5e672c0b69cd601ae/verification/audit_octic.py}{\nolinkurl{audit_octic.py}}\\
\(n=6,7\) &
  \shortstack[l]{Lemma~\ref{lem:last-cubic-profile}\\
                  Proposition~\ref{prop:last-line}} &
  \href{https://github.com/zhangteng2000/ehp-lemniscate-verification-20261108/blob/ca8c451bfbc7312b9c57acc5e672c0b69cd601ae/verification/verify_terminal_cubic_tree.py}{\nolinkurl{verify_terminal_cubic_tree.py}},
  \href{https://github.com/zhangteng2000/ehp-lemniscate-verification-20261108/blob/ca8c451bfbc7312b9c57acc5e672c0b69cd601ae/verification/verify_small_scalar_primary.py}{\nolinkurl{verify_small_scalar_primary.py}},
  \href{https://github.com/zhangteng2000/ehp-lemniscate-verification-20261108/blob/ca8c451bfbc7312b9c57acc5e672c0b69cd601ae/verification/audit_small_scalar_profiles.py}{\nolinkurl{audit_small_scalar_profiles.py}}\\
\(n=5\) & Proposition~\ref{prop:new-five-line} &
  \href{https://github.com/zhangteng2000/ehp-lemniscate-verification-20261108/blob/ca8c451bfbc7312b9c57acc5e672c0b69cd601ae/prior/quintic/replay_general_quintic.py}{\nolinkurl{prior/quintic/replay_general_quintic.py}}\\
\midrule
\multicolumn{3}{@{}l}{\emph{Passage to the circle}}\\
\(n\ge28\) & Proposition~\ref{prop:uniform-circle-transfer} &
  \href{https://github.com/zhangteng2000/ehp-lemniscate-verification-20261108/blob/ca8c451bfbc7312b9c57acc5e672c0b69cd601ae/verification/uniform_circle_transfer.py}{\nolinkurl{uniform_circle_transfer.py}}\\
\(9\le n\le27\) & Section~\ref{sec:collinear-circle} &
  \href{https://github.com/zhangteng2000/ehp-lemniscate-verification-20261108/blob/ca8c451bfbc7312b9c57acc5e672c0b69cd601ae/verification/single_coefficient_transfer.py}{\nolinkurl{single_coefficient_transfer.py}},
  \href{https://github.com/zhangteng2000/ehp-lemniscate-verification-20261108/blob/ca8c451bfbc7312b9c57acc5e672c0b69cd601ae/verification/independent_single_coefficient_transfer.py}{\nolinkurl{independent_single_coefficient_transfer.py}}\\
\(5\le n\le8\) & Proposition~\ref{prop:small-circle-transfer} &
  \href{https://github.com/zhangteng2000/ehp-lemniscate-verification-20261108/blob/ca8c451bfbc7312b9c57acc5e672c0b69cd601ae/verification/stationary_circle_transfer.py}{\nolinkurl{stationary_circle_transfer.py}},
  \href{https://github.com/zhangteng2000/ehp-lemniscate-verification-20261108/blob/ca8c451bfbc7312b9c57acc5e672c0b69cd601ae/verification/independent_stationary_circle_transfer.py}{\nolinkurl{independent_stationary_circle_transfer.py}}\\
\midrule
\(n=4\) &
  \shortstack[l]{Lemma~\ref{lem:quartic-line}\\
                  Proposition~\ref{prop:quartic-boundary-gap}} &
  \href{https://github.com/zhangteng2000/ehp-lemniscate-verification-20261108/blob/ca8c451bfbc7312b9c57acc5e672c0b69cd601ae/verification/quartic_line/verify_quartic_line.py}{\nolinkurl{quartic_line/verify_quartic_line.py}}\\
\bottomrule
\end{tabularx}
\par\smallskip
\caption{Principal verification programs. Paths are relative to
\texttt{verification/}, except for the path beginning with
\texttt{prior/}, which is relative to the root of the supplement.}
\label{tab:degree-verifiers}
\end{table}

For \(n\ge28\), Proposition~\ref{prop:uniform-circle-transfer}
proves the passage from lines to circles uniformly in \(n\).
The program \href{https://github.com/zhangteng2000/ehp-lemniscate-verification-20261108/blob/ca8c451bfbc7312b9c57acc5e672c0b69cd601ae/verification/uniform_circle_transfer.py}{\nolinkurl{uniform_circle_transfer.py}} checks the
rational endpoint inequalities in its proof. For \(9\le n\le27\),
the two circle programs reconstruct covers of the critical-value
diameter interval. For \(5\le n\le8\), the circle programs also
check the short-arc bounds. In degree four a single driver runs
both line computations and both sets of scalar bounds used in
Section~\ref{sec:quartic-global}.

For the parameter covers, the verifier starts with the stated
closed domain. At each subdivision it accounts for every child.
A terminal box must satisfy either an exclusion justified by the
hypotheses or the required inequality throughout the box. In the
cubic discriminant computations, each counted rectangle in
\((\theta,h)\) works for every parameter in the box, and the
rectangles have disjoint interiors. In the paired-height argument
of Proposition~\ref{prop:last-line}, at least one of \(h,-h\) passes
for each height in the stated range, so only one copy of that range
is counted unless both signs pass. The closed enclosures include
shared faces and admissible collisions of critical points or values.

The quartic verifiers use integer and rational enclosures throughout.
The alternate discriminant evaluations in degrees five through seven
use binary64 arithmetic with a separately proved error allowance.
The sextic and septic input trees and their error estimates are in
\nolinkurl{cubic_terminal_6/} (\href{https://github.com/zhangteng2000/ehp-lemniscate-verification-20261108/blob/ca8c451bfbc7312b9c57acc5e672c0b69cd601ae/verification/cubic_terminal_6/A0.txt}{\nolinkurl{A0.txt}}, \href{https://github.com/zhangteng2000/ehp-lemniscate-verification-20261108/blob/ca8c451bfbc7312b9c57acc5e672c0b69cd601ae/verification/cubic_terminal_6/A1.txt}{\nolinkurl{A1.txt}}, \href{https://github.com/zhangteng2000/ehp-lemniscate-verification-20261108/blob/ca8c451bfbc7312b9c57acc5e672c0b69cd601ae/verification/cubic_terminal_6/A2.txt}{\nolinkurl{A2.txt}}, \href{https://github.com/zhangteng2000/ehp-lemniscate-verification-20261108/blob/ca8c451bfbc7312b9c57acc5e672c0b69cd601ae/verification/cubic_terminal_6/B0.txt}{\nolinkurl{B0.txt}}, \href{https://github.com/zhangteng2000/ehp-lemniscate-verification-20261108/blob/ca8c451bfbc7312b9c57acc5e672c0b69cd601ae/verification/cubic_terminal_6/B1.txt}{\nolinkurl{B1.txt}}, \href{https://github.com/zhangteng2000/ehp-lemniscate-verification-20261108/blob/ca8c451bfbc7312b9c57acc5e672c0b69cd601ae/verification/cubic_terminal_6/B2.txt}{\nolinkurl{B2.txt}}, \href{https://github.com/zhangteng2000/ehp-lemniscate-verification-20261108/blob/ca8c451bfbc7312b9c57acc5e672c0b69cd601ae/verification/cubic_terminal_6/B3.txt}{\nolinkurl{B3.txt}}, \href{https://github.com/zhangteng2000/ehp-lemniscate-verification-20261108/blob/ca8c451bfbc7312b9c57acc5e672c0b69cd601ae/verification/cubic_terminal_6/roundoff.py}{\nolinkurl{roundoff.py}}) and \nolinkurl{cubic_terminal_7/} (\href{https://github.com/zhangteng2000/ehp-lemniscate-verification-20261108/blob/ca8c451bfbc7312b9c57acc5e672c0b69cd601ae/verification/cubic_terminal_7/A0.txt}{\nolinkurl{A0.txt}}, \href{https://github.com/zhangteng2000/ehp-lemniscate-verification-20261108/blob/ca8c451bfbc7312b9c57acc5e672c0b69cd601ae/verification/cubic_terminal_7/A1.txt}{\nolinkurl{A1.txt}}, \href{https://github.com/zhangteng2000/ehp-lemniscate-verification-20261108/blob/ca8c451bfbc7312b9c57acc5e672c0b69cd601ae/verification/cubic_terminal_7/A2.txt}{\nolinkurl{A2.txt}}, \href{https://github.com/zhangteng2000/ehp-lemniscate-verification-20261108/blob/ca8c451bfbc7312b9c57acc5e672c0b69cd601ae/verification/cubic_terminal_7/B0.txt}{\nolinkurl{B0.txt}}, \href{https://github.com/zhangteng2000/ehp-lemniscate-verification-20261108/blob/ca8c451bfbc7312b9c57acc5e672c0b69cd601ae/verification/cubic_terminal_7/B1.txt}{\nolinkurl{B1.txt}}, \href{https://github.com/zhangteng2000/ehp-lemniscate-verification-20261108/blob/ca8c451bfbc7312b9c57acc5e672c0b69cd601ae/verification/cubic_terminal_7/B2.txt}{\nolinkurl{B2.txt}}, \href{https://github.com/zhangteng2000/ehp-lemniscate-verification-20261108/blob/ca8c451bfbc7312b9c57acc5e672c0b69cd601ae/verification/cubic_terminal_7/B3.txt}{\nolinkurl{B3.txt}}, \href{https://github.com/zhangteng2000/ehp-lemniscate-verification-20261108/blob/ca8c451bfbc7312b9c57acc5e672c0b69cd601ae/verification/cubic_terminal_7/roundoff.py}{\nolinkurl{roundoff.py}}).
The \href{https://github.com/zhangteng2000/ehp-lemniscate-verification-20261108/blob/ca8c451bfbc7312b9c57acc5e672c0b69cd601ae/audit/code_and_data_sha256.json}{SHA-256 manifest}
lists the checksums of the archived code and data. The full replay
also checks the regenerated
\href{https://github.com/zhangteng2000/ehp-lemniscate-verification-20261108/blob/ca8c451bfbc7312b9c57acc5e672c0b69cd601ae/verification/joint_circle_intervals.json}{degree \(9\)--\(14\) circle certificate}
against its entry in the manifest. In degree nine, the \href{https://github.com/zhangteng2000/ehp-lemniscate-verification-20261108/blob/ca8c451bfbc7312b9c57acc5e672c0b69cd601ae/verification/audit_nonic_majorant_exact.py}{coefficient calculation}
gives \href{https://github.com/zhangteng2000/ehp-lemniscate-verification-20261108/blob/ca8c451bfbc7312b9c57acc5e672c0b69cd601ae/verification/audit_nonic_majorant_exact.json}{the exact bounds} for Lemma~\ref{lem:nonic-area-remainder}.
The \href{https://github.com/zhangteng2000/ehp-lemniscate-verification-20261108/blob/ca8c451bfbc7312b9c57acc5e672c0b69cd601ae/verification/joint_circle_certificate.py}{certificate generator}, its
\href{https://github.com/zhangteng2000/ehp-lemniscate-verification-20261108/blob/ca8c451bfbc7312b9c57acc5e672c0b69cd601ae/verification/audit_low_degree_transfer.py}{independent check} (\href{https://github.com/zhangteng2000/ehp-lemniscate-verification-20261108/blob/ca8c451bfbc7312b9c57acc5e672c0b69cd601ae/verification/low_degree_transfer_independent.json}{record}), and both verifications of
\eqref{eq:circle-nine-interval-test} all use \(9/10\).
The \href{https://github.com/zhangteng2000/ehp-lemniscate-verification-20261108/blob/ca8c451bfbc7312b9c57acc5e672c0b69cd601ae/verification/single_coefficient_transfer.json}{bisection record} and the
\href{https://github.com/zhangteng2000/ehp-lemniscate-verification-20261108/blob/ca8c451bfbc7312b9c57acc5e672c0b69cd601ae/verification/independent_single_coefficient_transfer.json}{trisection record} give the closed intervals
and the short-arc estimates used in Section~\ref{sec:collinear-circle}.

\end{document}